\documentclass[12pt]{amsart}
\usepackage[dvips]{graphicx}
\usepackage{amssymb}
\usepackage{hhline}
\usepackage{mathrsfs}

\newtheorem{defi}{Definition}[subsection]
\newtheorem{thm}[defi]{Theorem}
\newtheorem{lem}[defi]{Lemma}
\newtheorem{prop}[defi]{Proposition}
\newtheorem{cor}[defi]{Corollary}
\newtheorem{obs}[defi]{Observation}
\newtheorem{example}[defi]{Example}
\newtheorem{exercise}[defi]{Exercise}
\newtheorem{rem}[defi]{Remark}

\title{Lecture Notes on Generalized Heegaard Splittings}

\author{Toshio Saito, Martin Scharlemann and Jennifer Schultens}
\date{}

\begin{document}
\maketitle
\section{Introduction}

These notes grew out of a lecture series given at RIMS in the summer 
of 2001.  The authors were visiting RIMS in conjunction with the 
Research Project on Low-Dimensional Topology in the Twenty-First 
Century.  They had been invited by Professor Tsuyoshi Kobayashi.  The lecture 
series was first suggested by Professor Hitoshi Murakami. 

The lecture series was aimed at a broad audience that included many 
graduate students.  Its purpose lay in familiarizing the audience with 
the basics of $3$-manifold theory and introducing some topics of 
current research.  The first portion of the lecture series was devoted 
to standard topics in the theory of $3$-manifolds.  The middle portion 
was devoted to a brief study of Heegaaard splittings and generalized 
Heegaard splittings.  The latter portion touched on a brand new topic: 
fork complexes. 

During this time Professor Tsuyoshi Kobayashi had raised some interesting 
questions about the connectivity properties of generalized Heegaard 
splittings.  The latter portion of the lecture series was motivated by 
these questions.  And fork complexes were invented in an effort to 
illuminate some of the more subtle issues arising in the study of 
generalized Heegaard splittings. 

\begin{figure}\begin{center}
  {\unitlength=1cm
  \begin{picture}(4,4)
   \put(0,0){\includegraphics[keepaspectratio]{%
             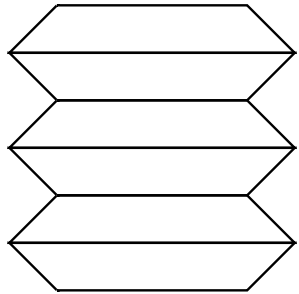}}
%
  \end{picture}}
  \caption{\sl The standard schematic diagram}
  \label{fig:fork1}
\end{center}\end{figure}

\begin{figure}[t]\begin{center}
  {\unitlength=1cm
  \begin{picture}(4,4)
   \put(0,0){\includegraphics[keepaspectratio]{%
             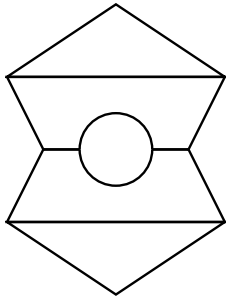}}
  \end{picture}}
\caption{\sl A more informative schematic diagram for a generalized
Heegaard splitting for a manifold homeomorphic to $($a surface$)\times S^1$}
  \label{fig:fork2}
\end{center}\end{figure}

In the standard schematic diagram for generalized Heegaard splittings, 
Heegaard splittings are stacked on top of each other in a linear 
fashion.  See Figure \ref{fig:fork1}.  This can cause confusion in those 
cases in which generalized Heegaaard splittings possess interesting 
connectivity properties.  In these cases, some of the topological 
features of the $3$-manifold are captured by the connectivity 
properties of the generalized Heegaard splitting rather than by the 
Heegaard splittings of submanifolds into which the generalized 
Heegaard splitting decomposes the $3$-manifold. See Figure \ref{fig:fork2}. 
Fork complexes provide a means of description in this context. 

The authors would like to express their appreciation of the 
hospitality extended to them during their stay at RIMS.  They would 
also like to thank the many people that made their stay at RIMS 
delightful, illuminating and productive, most notably Professor Hitoshi 
Murakami, Professor Tsuyoshi Kobayashi, Professor Jun Murakami, 
Professor Tomotada Ohtsuki, Professor Kyoji Saito, 
Professor Makoto Sakuma, Professor Kouki Taniyama and Dr. Yo'av Rieck. 
Finally, they would like to thank Dr. Ryosuke Yamamoto for drawing the fine 
pictures in these lecture notes. 

\section{Preliminaries}
\subsection{PL $3$-manifolds}

Let $M$ be a PL $3$-manifold, i.e., $M$ is a union of 
$3$-simplices $\sigma^3_{i}$ $(i=1,2,\dots, t)$ 
such that $\sigma^3_{i}\cap \sigma^3_{j}$ $(i\ne j)$ is emptyset, a vertex, an edge or a face 
and that for each vertex $v$, $\bigcup_{v\in \sigma^3_j} \sigma^3_j$ is a $3$-ball (\textit{cf}. \cite{Rourke}). 
Then the decomposition $\{\sigma^3_i\}_{1\le i\le t}$ of $M$ 
is called a \textit{triangulation} of $M$. 

\begin{example}
\rm
\begin{enumerate}
\item The $3$-ball $B^3$ is the simplest PL $3$-manifold in a sense that $B^3$ is homeomorphic to a 
$3$-simplex. 
\item The $3$-sphere $S^3$ is a $3$-manifold obtained from two $3$-balls by attaching their boundaries. 
Since $S^3$ is homeomorphic to the boundary of a $4$-simplex, we see that $S^3$ is a union of five $3$-simplices. 
It is easy to show that this gives a triangulation of $S^3$. 
\end{enumerate}
\end{example}

\begin{exercise}
\rm
Show that the following $3$-manifolds are PL $3$-manifolds. 
\begin{enumerate}
\item The solid torus $D^2\times S^1$. 
\item $S^2\times S^1$. 
\item The lens spaces. Note that a \textit{lens space} is obtained from two solid tori by attaching 
their boundaries. 
\end{enumerate}
\end{exercise}

\begin{figure}[htb]\begin{center}
  {\unitlength=1cm
  \begin{picture}(10.0,5.9)
   \put(2.9,3.5){\includegraphics[keepaspectratio,height=2.4cm]{%
   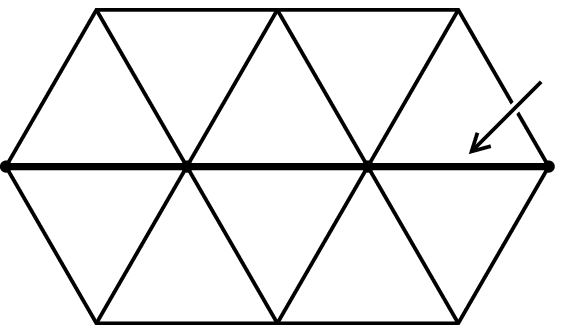}}
   \put(0,0){\includegraphics[keepaspectratio,height=2.4cm]{%
   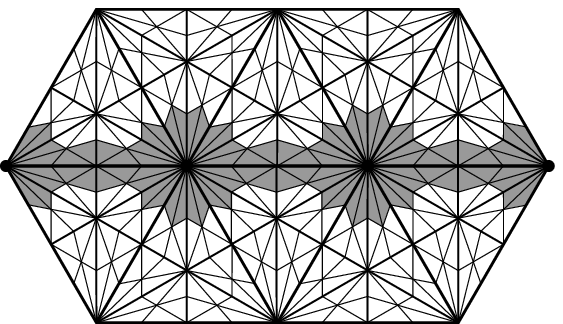}}
   \put(5.8,0){\includegraphics[keepaspectratio,height=2.4cm]{%
   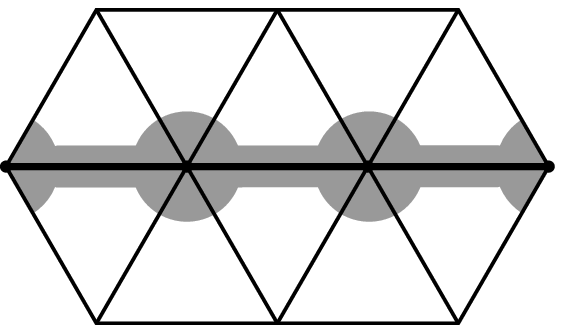}}
   \put(7.0,5.4){$X$}
   \put(2.5,3.5){$K$}
   \put(4.85,2.8){$\Downarrow$}
   \put(4.85,1.05){$\cong$}
   \put(4.1,0.6){$\nwarrow$}
   \put(5.4,0.6){$\nearrow$}
   \put(4.3,0.2){$\eta(X;K)$}
  \end{picture}}
  \caption{}
  \label{fig:0.1}
\end{center}\end{figure}

Let $K$ be a three dimensional simplicial complex and $X$ a sub-complex of $K$, that is, 
$X$ a union of vertices, edges, faces and $3$-simplices of $K$ such that 
$X$ is a simplicial complex. Let $K''$ be the second barycentric subdivision of $K$. 
A \textit{regular neighborhood} of $X$ in $K$, denoted by $\eta (X;K)$, is a union of the $3$-simplices 
of $K''$ intersecting $X$ (\textit{cf}. Figure \ref{fig:0.1}). 

\begin{figure}[htb]\begin{center}
  {\unitlength=1cm
  \begin{picture}(8.0,3.0)
   \put(0,0){\includegraphics[keepaspectratio,height=3cm]{%
   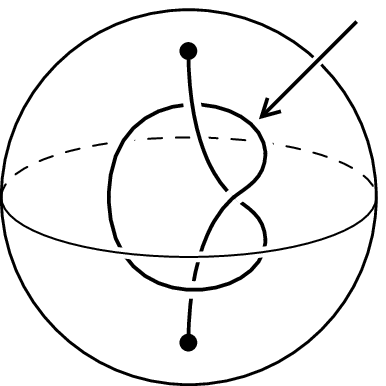}}
   \put(5,0){\includegraphics[keepaspectratio,height=3cm]{%
   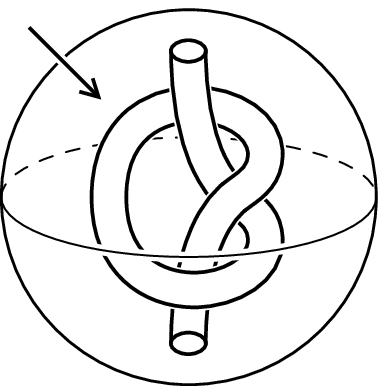}}
   \put(-.2,2.8){$M$}
   \put(3.6,2.7){$\eta(X;M)$}
   \put(2.9,2.9){$X$}
  \end{picture}}
  \caption{}
  \label{fig:0.2}
\end{center}\end{figure}

\begin{prop}
If $X$ is a PL $1$-manifold properly embedded in a PL $3$-manifold $M$ 
(namely, $X\cap \partial M=\partial X$), 
then $\eta (X;M)\cong X\times B^2$, where $X$ is identified with $X\times \{$a center of $B^2\}$ 
and $\eta (X;M)\cap \partial M$ is identified with $\partial X \times B^2$ 
(\textit{cf}. Figure \ref{fig:0.2}). 
\end{prop}

\begin{prop}
Suppose that a PL $3$-manifold $M$ is orientable. 
If $X$ is an orientable PL $2$-manifold properly embedded in $M$ (namely, $X\cap \partial M=\partial X$), 
then $\eta (X;M)\cong X\times [0,1]$, where $X$ is identified with $X\times \{1/2\}$ and 
$\eta (X;M)\cap \partial M$ is identified with $\partial X \times [0,1]$. 
\end{prop}

\begin{thm}[Moise \cite{Moise}]\label{Moise}
Every compact $3$-manifold is a PL $3$-manifold. 
\end{thm}

In the remainder of these notes, we work in the PL category unless otherwise specified. 

\subsection{Fundamental definitions}
By the term \textit{surface}, we will mean a connected compact $2$-manifold. 

Let $F$ be a surface.  A loop $\alpha$ in $F$ is said to be \textit{inessential} in $F$ if 
$\alpha$ bounds a disk in $F$, otherwise $\alpha$ is said to be \textit{essential} in $F$. 
An arc $\gamma$ properly embedded in $F$ is said to be \textit{inessential} in $F$ if 
$\gamma$ cuts off a disk from $F$, otherwise $\gamma$ is said to be \textit{essential} in $F$. 

Let $M$ be a compact orientable $3$-manifold. 
A disk $D$ properly embedded in $M$ is said to be \textit{inessential} in $M$ if 
$D$ cuts off a $3$-ball from $M$, otherwise $D$ is said to be \textit{essential} in $M$.
A $2$-sphere $P$ properly embedded in $M$ is said to be \textit{inessential} in $M$ if 
$P$ bounds a $3$-ball in $M$, otherwise $P$ is said to be \textit{essential} in $M$. 
Let $F$ be a surface properly embedded in $M$. We say that $F$ is $\partial$-\textit{parallel} in $M$ 
if $F$ cuts off a $3$-manifold homeomorphic to $F\times [0,1]$ from $M$. 
We say that $F$ is \textit{compressible} in $M$ if there is a disk $D\subset M$ such that 
$D\cap F=\partial D$ and $\partial D$ is an essential loop in $F$. Such a disk $D$ is called a 
\textit{compressing disk}. We say that $F$ is \textit{incompressible} in $M$ if
$F$ is not compressible in $M$. The surface $F$ is $\partial$-\textit{compressible} in $M$ if 
there is a disk $\delta\subset M$ such that $\delta\cap F$ is an arc which is essential in $F$, 
say $\gamma$, in $F$ and that 
$\delta\cap \partial M$ is an arc, say $\gamma'$, with $\gamma'\cup \gamma=\partial \delta$. 
Otherwise $F$ is said to be $\partial$-\textit{incompressible} in $M$. 
Suppose that $F$ is homeomorphic neither to a disk nor to a $2$-sphere. The surface 
$F$ is said to be \textit{essential} in $M$ if $F$ is incompressible in $M$ and is not $\partial$-parallel 
in $M$. 

\begin{defi}\label{D0}
\rm
Let $M$ be a connected compact orientable $3$-manifold. 
\begin{enumerate}
\item $M$ is said to be \textit{reducible} if there is a $2$-sphere in $M$  which does not 
bound a $3$-ball in $M$. Such a $2$-sphere is called a \textit{reducing $2$-sphere} of $M$. 
$M$ is said to be \textit{irreducible} if $M$ is not reducible. 
\item $M$ is said to be \textit{$\partial$-reducible} if there is a disk properly embedded in $M$ 
whose boundary is essential in $\partial M$. Such a disk is called a \textit{$\partial$-reducing disk}. 
\end{enumerate}
\end{defi}

\section{Heegaard splittings}

\subsection{Definitions and fundamental properties}
\begin{figure}[htb]\begin{center}
  {\unitlength=1cm
  \begin{picture}(10.5,10.5)
   \put(0,6.4){\includegraphics[keepaspectratio]{
               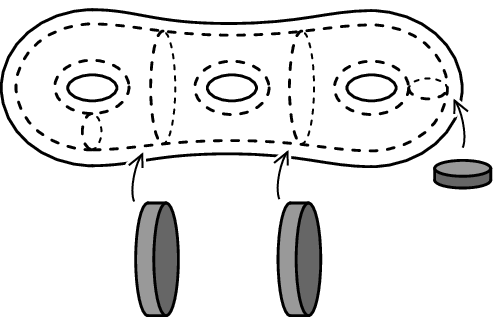}}
   \put(0,2.8){\includegraphics[keepaspectratio]{
               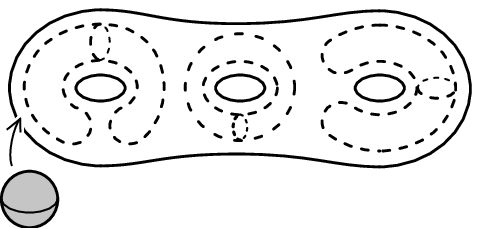}}
   \put(3,0){\includegraphics[keepaspectratio]{
             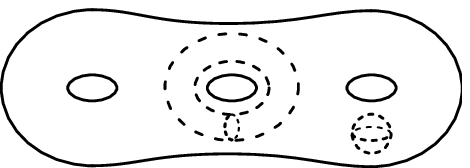}}
   \put(6,6){\includegraphics[keepaspectratio,width=4.5cm]{%
             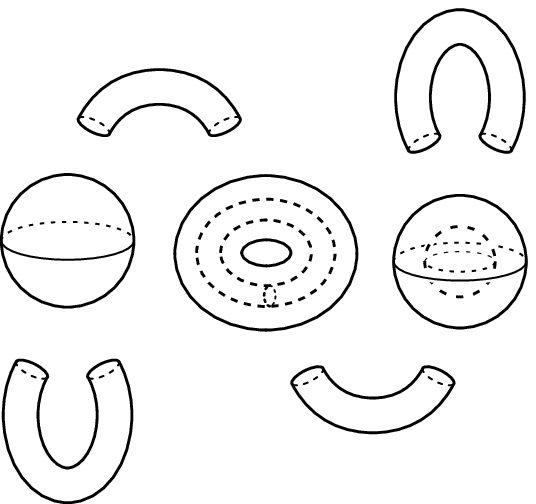}}
   \put(6,3){\includegraphics[keepaspectratio,width=4.5cm]{%
             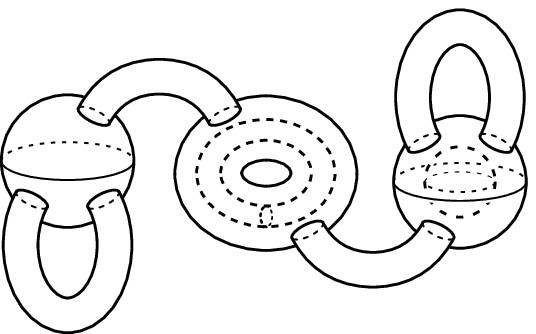}}
   \put(-.7,7.6){$F\times [0,1]$}
   \put(6.8,10.2){Dual discription}
   \put(2.2,5.8){$\downarrow$}
   \put(8.2,5.8){$\downarrow$}
   \put(3.2,2.5){$\searrow$}
   \put(7.2,2.5){$\swarrow$}
  \end{picture}}
  \caption{}
  \label{fig:0.3}
\end{center}\end{figure}

\begin{defi}\label{D1}
\rm
A $3$-manifold $C$ is called a \textit{compression body} if there exists a closed surface 
$F$ such that $C$ is obtained from $F\times [0,1]$ by attaching $2$-handles along mutually disjoint loops in 
$S\times \{1\}$ and filling in some resulting $2$-sphere boundary components with $3$-handles 
(\textit{cf}. Figure \ref{fig:0.3}). 
We denote $F\times \{0\}$ by $\partial_+ C$ and $\partial C\setminus \partial_+ C$ by $\partial_- C$. 
A compression body $C$ is called a \textit{handlebody} if $\partial_- C=\emptyset$. 
A compression body $C$ is said to be \textit{trivial} if $C\cong F\times [0,1]$. 
\end{defi}

\begin{defi}
\rm
For a compression body $C$, an essential disk in $C$ is called a \textit{meridian disk} of $C$. 
A union $\Delta$ of mutually disjoint meridian disks of $C$ 
is called a \textit{complete meridian system} if the manifold obtained from $C$ by cutting along 
$\Delta$ are the union of $\partial_- C\times [0,1]$ and (possibly empty) $3$-balls. 
A complete meridian system $\Delta$ of $C$ is \textit{minimal} if the number of the components of 
$\Delta$ is minimal among all complete meridian system of $C$. 
\end{defi}

\begin{rem}\label{R1}
\rm
The following properties are known for compression bodies. 
\begin{enumerate}
\item A compression body $C$ is reducible if and only if $\partial_- C$ contains a $2$-sphere component. 
\item A minimal complete meridian system $\Delta$ of a compression body $C$ 
cuts $C$ into $\partial_- C \times [0,1]$ if $\partial_- C\ne \emptyset$, and 
$\Delta$ cuts $C$ into a $3$-ball if $\partial_- C=\emptyset$ (hence $C$ is a handlebody). 
\item By extending the cores of the $2$-handles in the definition of the compression body $C$ vertically 
to $F\times [0,1]$, we obtain a complete meridian system $\Delta$ of $C$ such that the manifold obtained 
by cutting $C$ along $\Delta$ is homeomorphic to a union of $\partial_- C\times [0,1]$ and some 
(possibly empty) $3$-balls. This gives a \textit{dual description} of compression bodies. 
That is, a compression body $C$ is obtained from $\partial_- C\times [0,1]$ and some 
(possibly empty) $3$-balls by attaching some $1$-handles to $\partial_- C\times \{1\}$ and 
the boundary of the $3$-balls (\textit{cf}. Figure \ref{fig:0.3}). 
\item For any compression body $C$, $\partial_- C$ is incompressible in $C$. 
\item Let $C$ and $C'$ be compression bodies. Suppose that $C''$ is obtained from $C$ and $C'$ by 
identifying a component of $\partial_- C$ and $\partial_+ C'$. Then $C''$ is a compression body. 
\item Let $D$ be a meridian disk of a compression body $C$. Then there is a complete meridian system 
$\Delta$ of $C$ such that $D$ is a component of $\Delta$. Any component obtained by cutting $C$ along $D$ 
is a compression body. 
\end{enumerate}
\end{rem}

\begin{exercise}
\rm
Show Remark \ref{R1}. 
\end{exercise}

An annulus $A$ properly embedded in a compression body $C$ is called a \textit{spanning annulus} 
if $A$ is incompressible in $C$ and a component of $\partial A$ is contained in $\partial_+ C$ 
and the other is contained in $\partial_- C$.

\begin{lem}\label{spanning}
Let $C$ be a non-trivial compression body. 
Let $A$ be a spanning annulus in $C$. 
Then there is a meridian disk $D$ of $C$ with $D\cap A=\emptyset$. 
\end{lem}

\begin{proof}
Since $C$ is non-trivial, there is a meridian disk of $C$. 
We choose a meridian disk $D$ of $C$ such that $D$ intersects $A$ transversely and 
$|D\cap A|$ is minimal among all such meridian disks.  
Note that $A\cap \partial_- C$ is an essential loop in 
the component of $\partial_- C$ containing $A\cap \partial_- C$. 
We shall prove that $D\cap A=\emptyset$. To this end, we suppose $D\cap A\ne \emptyset$. 

\medskip
\noindent\textit{Claim 1.}\ \ \ \ 
There are no loop components of $D\cap A$. 

\medskip
\noindent\textit{Proof.}\ \ \ \ 
Suppose that $D\cap A$ has a loop component which is inessential 
in $A$. Let $\alpha$ be a loop component of $D\cap A$ which is \textit{innermost} in $A$, that is, 
$\alpha$ cuts off a disk $\delta_{\alpha}$ from $A$ such that the interior of $\delta_{\alpha}$ is disjoint 
from $D$. Such a disk $\delta_{\alpha}$ is called an \textit{innermost disk} for $\alpha$. 
We remark that $\alpha$ is not necessarily innermost in $D$. 
Note that $\alpha$ also bounds a disk in $D$, say $\delta'_{\alpha}$. 
Then we obtain a disk $D'$ by applying \textit{cut and paste operation} on $D$ with 
using $\delta_{\alpha}$ and $\delta'_{\alpha}$, i.e., $D'$ is obtained from $D$ by removing the interior of 
$\delta'_{\alpha}$ and then attaching $\delta_{\alpha}$ (\textit{cf}. Figure \ref{fig:04}). 
Note that $D'$ is a meridian disk of $C$. Moreover, 
we can isotope the interior of $D'$ slightly so that $|D'\cap A|<|D\cap A|$, a contradiction. 
(Such an argument as above is called an \textit{innermost disk argument}.) 

\begin{figure}[htb]\begin{center}
  {\unitlength=1cm
  \begin{picture}(12,4)
   \put(0,.5){\includegraphics[keepaspectratio]{%
   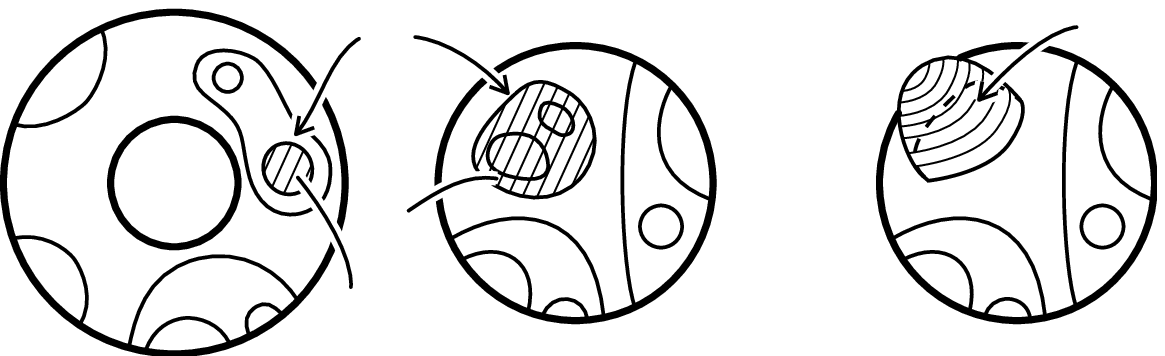}}
   \put(1.6,0){$A$}
   \put(3.82,3.7){$\alpha$}
   \put(3.82,1.7){$\delta'_{\alpha}$}
   \put(3.4,.9){$\delta_{\alpha}$}
   \put(5.8,0){$D$}
   \put(7.8,2.2){$\Longrightarrow$}
   \put(10.3,0){$D'$}
   \put(11,3.8){$\delta_{\alpha}$}
  \end{picture}}
  \caption{}
  \label{fig:04}
\end{center}\end{figure}

Hence if $D\cap A$ has a loop component, we may assume that the loop is 
essential in $A$. Let $\alpha'$ be a loop component of $D\cap A$ which is innermost in $D$, 
and let $\delta_{\alpha'}$ be the innermost disk in $D$ with $\partial \delta_{\alpha'}=\alpha'$. 
Then $\alpha'$ cuts $A$ into two annuli, and let $A'$ be the component obtained by cutting $A$ along 
$\alpha'$ such that $A'$ is adjacent to $\partial_- C$. Set $D''=A'\cup \delta_{\alpha'}$. 
Then $D''(\subset C)$ is a compressing disk of $\partial_- C$, contradicting (4) of Remark \ref{R1}. 
Hence we have Claim 1. 

\medskip
\noindent\textit{Claim 2.}\ \ \ \ 
There are no arc components of $D\cap A$. 

\medskip
\noindent\textit{Proof.}\ \ \ \ 
Suppose that there is an arc component of $D\cap A$. 
Note that $\partial D\subset \partial_+ C$. Hence we may assume that each component of $D\cap A$ 
is an inessential arc in $A$ whose endpoints are contained in $\partial_+ C$. 
Let $\gamma$ be an arc component of $D\cap A$ which is \textit{outermost} in $A$, that is, 
$\gamma$ cuts off a disk $\delta_{\gamma}$ from $A$ such that the interior of $\delta_{\gamma}$ is disjoint 
from $D$. Such a disk $\delta_{\gamma}$ is called an \textit{outermost disk} for $\gamma$. 
Note that $\gamma$ cuts $D$ into two disks $\bar{\delta}_{\gamma}$ and $\bar{\delta}'_{\gamma}$ 
(\textit{cf}. Figure \ref{fig:05}). 

\begin{figure}[htb]\begin{center}
  {\unitlength=1cm
  \begin{picture}(9.4,3.8)
   \put(0,.2){\includegraphics[keepaspectratio]{%
   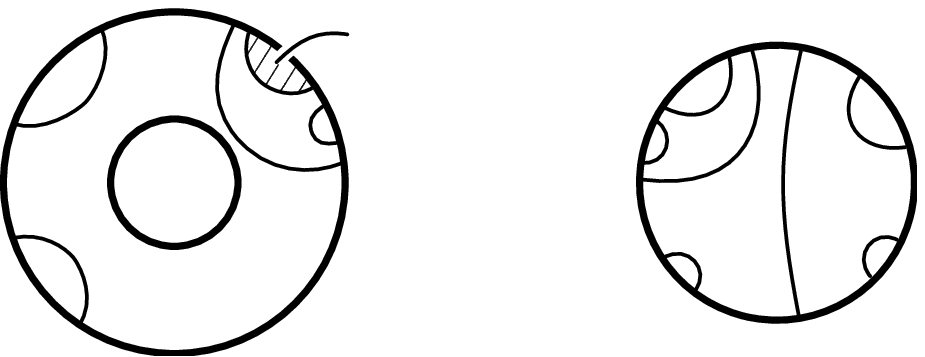}}
   \put(2.7,2.7){$\gamma$}
   \put(3.6,3.4){$\delta_{\gamma}$}
   \put(3.6,.8){$A$}
   \put(7.64,1.8){$\gamma$}
   \put(6.5,.5){$\underbrace{\hspace{1.5cm}}$}
   \put(7.1,-.3){$\bar{\delta}_{\gamma}$}
   \put(8.16,.5){$\underbrace{\hspace{1.2cm}}$}
   \put(8.6,-.3){$\bar{\delta}'_{\gamma}$}
   \put(9.5,.8){$D$} 
 \end{picture}}
  \caption{}
  \label{fig:05}
\end{center}\end{figure}

If both $\bar{\delta}_{\gamma}\cup \delta_{\gamma}$ and $\bar{\delta}'_{\gamma}\cup \delta_{\gamma}$ 
are inessential in $C$, then $D$ is also inessential in $C$, a contradiction. So we may assume that 
$\bar{D}=\bar{\delta}_{\gamma}\cup \delta_{\gamma}$ is essential in $C$. Then we can isotope 
$\bar{D}$ slightly so that $|\bar{D}\cap A|<|D\cap A|$, a contradiction. Hence we have Claim 2. 
(Such an argument as above is called an \textit{outermost disk argument}.) 

\medskip
Hence it follows from Claims 1 and 2 that $D\cap A=\emptyset$, 
and this completes the proof of Lemma \ref{spanning}. 
\end{proof}

\begin{rem}\label{R0}
\rm
Let $A$ be a spanning annulus in a non-trivial compression body $C$. 
\begin{enumerate}
\item By using the arguments of the proof of Lemma \ref{spanning}, we can show that there is a 
complete meridian system $\Delta$ of $C$ with $\Delta\cap A=\emptyset$. 
\item It follows from $(1)$ above that there is a meridian disk $E$ of $C$ such that $E\cap A=\emptyset$ 
and $E$ cuts off a $3$-manifold which is homeomorphic to $($a closed surface$)\times [0,1]$ containing $A$. 
\end{enumerate}
\end{rem}

\begin{exercise}
\rm
Show Remark \ref{R0}. 
\end{exercise}

Let $\bar{\alpha}=\alpha_1\cup \dots \cup \alpha_p$ be a union of mutually disjoint arcs in a 
compression body $C$. We say that $\bar{\alpha}$ is \textit{vertical} if there is a union of mutually disjoint 
spanning annuli $A_1\cup \dots \cup A_p$ in $C$ such that $\alpha_i\cap A_j=\emptyset$ $(i\ne j)$ and 
$\alpha_i$ is an essential arc properly embedded in $A_i$ $(i=1,2,\dots,p)$. 

\begin{lem}\label{vertical}
Suppose that $\bar{\alpha}=\alpha_1\cup \dots \cup \alpha_p$ is vertical in $C$. Let $D$ be a meridian disk of $C$. 
Then there is a meridian disk $D'$ of $C$ with $D'\cap \bar{\alpha}=\emptyset$ which is obtained by cut-and-paste 
operation on $D$. Particularly, 
if $C$ is irreducible, then $D$ is ambient isotopic such that $D\cap \bar{\alpha}=\emptyset$. 
\end{lem}

\begin{proof}
Let $\bar{A}=A_1\cup \dots \cup A_p$ be a union of annuli for $\bar{\alpha}$ as above. 
By using innermost disk arguments, we see that there is a meridian disk $D'$ such that no components 
of $D'\cap \bar{A}$ are loops which are inessential in $\bar{A}$. 
We remark that $D'$ is ambient isotopic to $D$ if $C$ is irreducible.
Note that each component of 
$\bar{A}$ is incompressible in $C$. Hence no components of $D'\cap \bar{A}$ are loops which 
are essential in $\bar{A}$. Hence each component of $D'\cap \bar{A}$ is an arc; moreover 
since $\partial D$ is contained in $\partial_+ C$, the endpoints of the arc components of 
$D'\cap \bar{A}$ are contained in $\partial_+ C\cap \bar{A}$. Then it is easy to see that 
there exists an arc $\beta_i(\subset A_i)$ such that $\beta_i$ is essential in $A_i$ and 
$\beta_i\cap D'=\emptyset$. Take an ambient isotopy $h_t$ $(0\le t\le 1)$ of $C$ such that 
$h_0(\beta_i)=\beta_i$, $h_t(\bar{A})=\bar{A}$ and $h_1(\beta_i)=\alpha_i$ 
$(i=1,2,\dots,p)$ (\textit{cf}. Figure \ref{fig:A01}). 
Then the ambient isotopy $h_t$ assures that $D'$ is isotoped 
so that $D'$ is disjoint from $\bar{\alpha}$. 

\begin{figure}[htb]\begin{center}
  {\unitlength=1cm
  \begin{picture}(4,4)
   \put(0,0){\includegraphics[keepaspectratio]{%
             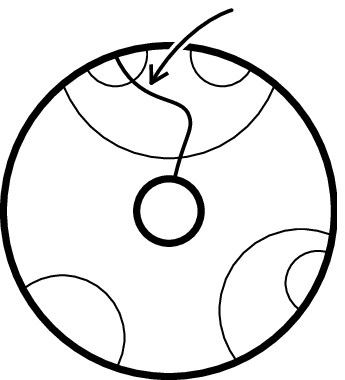}}
   \put(3.2,.2){$A_i$}
   \put(2.5,3.8){$\alpha_i$}
  \end{picture}}
  \caption{}
  \label{fig:A01}
\end{center}\end{figure}

\end{proof}

In the remainder of these notes, let $M$ be a connected compact orientable $3$-manifold. 

\begin{defi}
\rm
Let $(\partial_1 M,\partial_2 M)$ be a partition of $\partial$-components of $M$. 
A triplet $(C_1,C_2;S)$ is 
called a \textit{Heegaard splitting} of $(M;\partial_1 M,\partial_2 M)$ if $C_1$ and $C_2$ 
are compression bodies with $C_1\cup C_2=M$, $\partial_- C_1=\partial_1 M$, 
$\partial_- C_2=\partial_2 M$ and $C_1\cap C_2=\partial_+ C_1=\partial_+ C_2=S$. 
The surface $S$ is called a \textit{Heegaard surface} and 
the \textit{genus} of a Heegaard splitting is defined by the genus of the Heegaard surface. 
\end{defi}

\begin{thm}\label{S1}
For any partition $(\partial_1 M,\partial_2 M)$ of the boundary components of $M$, 
there is a Heegaard splitting of $(M;\partial_1 M,\partial_2 M)$. 
\end{thm}

\begin{proof}
It follows from Theorem \ref{Moise} that $M$ is triangulated, that is, there is a finite simplicial 
complex $K$ which is homeomorphic to $M$. Let $K'$ be a barycentric subdivision of $K$ and 
$K_1$ the $1$-skeleton of $K$. Here, a $1$\textit{-skeleton} of $K$ is a union of the vertices and edges of $K$. 
Let $K_2 \subset K'$ be the dual $1$-skeleton (see Figure \ref{fig:002}). 
Then each of $K_i$ $(i=1,2)$ is a finite graph in $M$. 

\begin{figure}\begin{center}
  {\unitlength=1cm
  \begin{picture}(4.0,4.0)
   \put(0,0){\includegraphics[keepaspectratio,height=4cm]{%
   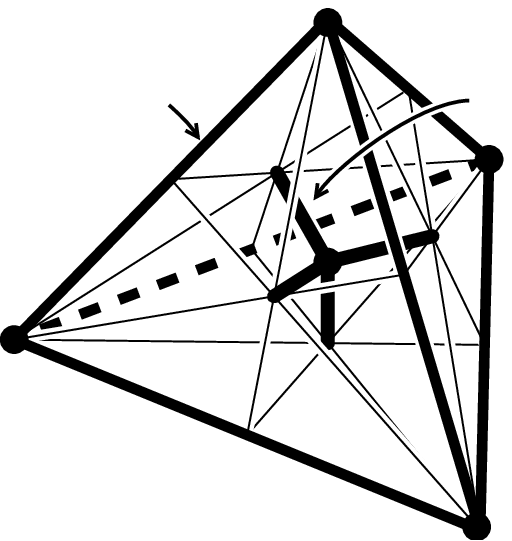}}
   \put(0.7,3.2){$K_1$}
   \put(3.5,3.2){$K_2$}
  \end{picture}}
  \caption{}
  \label{fig:002}
\end{center}\end{figure}

\medskip
\noindent\textit{Case 1.}\ \ \ \ 
$\partial M=\emptyset$. 

\medskip
Recall that $K_1$ consists of $0$-simplices and $1$-simplices. 
Set $C_1=\eta (K_1;M)$ and $C_2=\eta (K_2;M)$. 
Note that a regular neighborhood of a $0$-simplex corresponds to 
a $0$-handle and that a regular neighborhood of a $1$-simplex corresponds to 
a $1$-handle. Hence $C_1$ is a handlebody. Similarly, we see that $C_2$ is also a handlebody. 
Then we see that $C_1\cup C_2=M$ and $C_1\cap C_2=\partial C_1=\partial C_2$. 
Hence $(C_1,C_2;S)$ is a Heegaard splitting of $M$ with $S=C_1\cap C_2$. 

\medskip
\noindent\textit{Case 2.}\ \ \ \ 
$\partial M\ne \emptyset$. 

\medskip
In this case, we first take the barycentric subdivision of $K$ and use 
the same notation $K$. Recall that $K'$ is the barycentric subdivision of $K$. 
Note that no $3$-simplices of $K$ intersect both $\partial_1 M$ and $\partial_2 M$. 
Let $N(\partial_2 M)$ be a union of the $3$-simplices in $K'$ intersecting $\partial_2 M$. 
Then $N(\partial_2 M)$ is homeomorphic to $\partial_2 M\times [0,1]$, 
where $\partial_2 M\times \{0 \}$ is identified with $\partial_2 M$. Set 
$\partial'_2 M=\partial_2 M\times \{1 \}$. 
Let $\bar{K}_1$ ($\bar{K}_2$ resp.) be the maximal sub-complex of $K_1$ ($K_2$ resp.) 
such that $\bar{K}_1$ ($\bar{K}_2$ resp.) 
is disjoint from $\partial'_2 M$ ($\partial_1 M$ resp.) (\textit{cf}. Figure \ref{fig:02}). 

\begin{figure}\begin{center}
  {\unitlength=1cm
  \begin{picture}(8.8,4.0)
   \put(0.2,0){\includegraphics[keepaspectratio,height=4cm]{%
   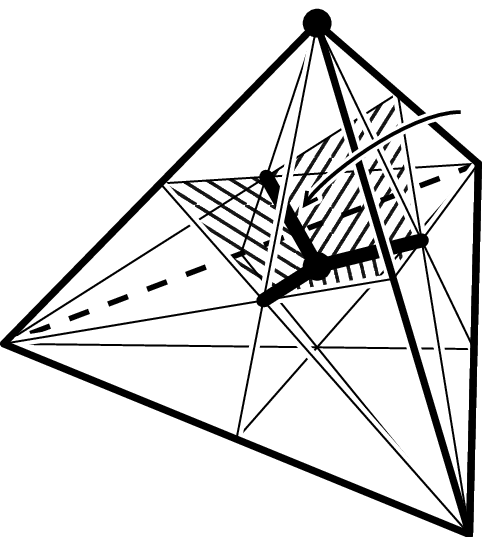}}
   \put(5.0,0){\includegraphics[keepaspectratio,height=4cm]{%
   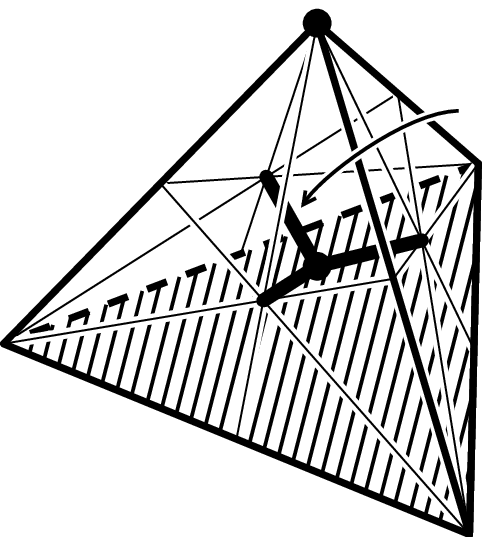}}
   \put(2.4,4.1){$\bar{K}_1$}
   \put(3.7,3.1){$\bar{K}_2$}
   \put(.4,2.6){${\partial'_2}M$}
   \put(-.3,1.0){${\partial_2}M$}
   \put(7.2,4.1){$\bar{K}_1$}
   \put(8.5,3.1){$\bar{K}_2$}
   \put(4.5,1.0){${\partial_1}M$}
  \end{picture}}
  \caption{}
  \label{fig:02}
\end{center}\end{figure}

Set $C_1=\eta(\partial_1 M\cup \bar{K}_1;M)$. Note that $C_1=\eta(\partial_1 M;M)\cup \eta(\bar{K}_1;M)$. 
Note again that a regular neighborhood of a $0$-simplex corresponds to 
a $0$-handle and that a regular neighborhood of a $1$-simplex corresponds to 
a $1$-handle. Hence $C_1$ is obtained from $\partial_1 M\times [0,1]$ by attaching $0$-handles and $1$-handles 
and therefore $C_1$ is a compression body with $\partial_- C_1=\partial_1 M$. 
Set $C_2=\eta(N(\partial_2 M)\cup \bar{K}_2;M)$. 
By the same argument, we can see that $C_2$ is a compression body with $\partial_- C_2=\partial_2 M$. 
Note that $C_1\cup C_2=M$ and $C_1\cap C_2=\partial C_1=\partial C_2$. 
Hence $(C_1,C_2;S)$ is a Heegaard splitting of $M$ with $S=C_1\cap C_2$. 
\end{proof}

We now introduce alternative viewpoints to Heegaard splittings as remarks below. 

\begin{figure}[htb]\begin{center}
  {\unitlength=1cm
  \begin{picture}(8.0,3.0)
   \put(0,0){\includegraphics[keepaspectratio]{
   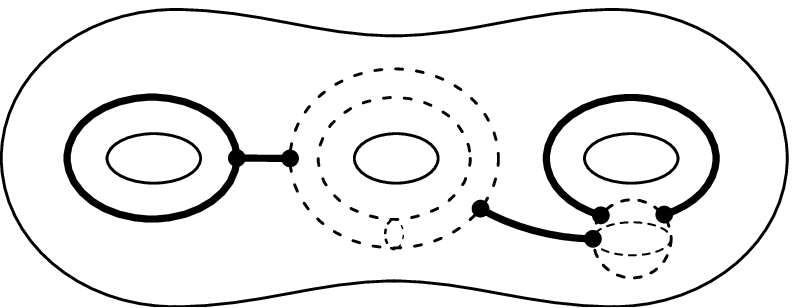}}
  \end{picture}}
  \caption{}
  \label{fig:02.1}
\end{center}\end{figure}

\begin{defi}
\rm
Let $C$ be a compression body. A finite graph $\Sigma$ in $C$ is called a \textit{spine} 
of $C$ if $C\setminus (\partial_- C\cup \Sigma)\cong \partial_+ C \times [0,1)$ and every vertex of 
valence one is in $\partial_- C$ (\textit{cf}. Figure \ref{fig:02.1}). 
\end{defi}

\begin{rem}
\rm
Let $(C_1,C_2;S)$ be a Heegaard splitting of $(M;\partial_1 M,\partial_2 M)$. 
Let $\Sigma_i$ be a spine of $C_i$, and set $\Sigma_i'=\partial_i M\cup \Sigma_i$ $(i=1,2)$. 
Then 

\bigskip
\begin{center}
$M\setminus (\Sigma_1'\cup \Sigma_2')=(C_1\setminus \Sigma_1')\cup_S (C_2\setminus \Sigma_2')\cong 
S\times (0,1)$. 
\end{center}

\medskip
Hence there is a continuous function $f:M\to [0,1]$ such that $f^{-1}(0)=\Sigma_1'$, 
$f^{-1}(1)=\Sigma_2'$ and $f^{-1}(t)\cong S$ $(0<t<1)$. This is called a 
\textit{sweep-out picture}. 
\end{rem}

\begin{rem}\label{R2}
\rm
Let $(C_1,C_2;S)$ be a Heegaard splitting of $(M;\partial_1 M,\partial_2 M)$. 
By a dual description of $C_1$, we see that $C_1$ is obtained from $\partial_1 M\times [0,1]$ and $0$-handles 
$\mathcal{H}^0$ by attaching $1$-handles $\mathcal{H}^1$. By Definition \ref{D1}, 
$C_2$ is obtained from $S\times [0,1]$ by attaching $2$-handles $\mathcal{H}^2$ and filling some $2$-sphere 
boundary components with $3$-handles $\mathcal{H}^3$. Hence we obtain the following decomposition of $M$: 

\bigskip
\begin{center}
$M=\partial_1 M\times [0,1]\cup \mathcal{H}^0\cup \mathcal{H}^1\cup S\times [0,1]\cup 
\mathcal{H}^2\cup \mathcal{H}^3$. 
\end{center}

\medskip
By collapsing $S\times [0,1]$ to $S$, we have: 

\bigskip
\begin{center}
$M=\partial_1 M\times [0,1]\cup \mathcal{H}^0\cup \mathcal{H}^1\cup_S 
\mathcal{H}^2\cup \mathcal{H}^3$. 
\end{center}

\medskip
This is called a \textit{handle decomposition of $M$ induced from $(C_1,C_2;S)$}. 
\end{rem}

\begin{defi}\label{p}
\rm
Let $(C_1,C_2;S)$ be a Heegaard splitting of $(M;\partial_1 M,\partial_2 M)$. 

\begin{enumerate}
\item The splitting $(C_1,C_2;S)$ is said to be \textit{reducible} if there are 
meridian disks $D_i$ $(i=1,2)$ of $C_i$ with $\partial D_1=\partial D_2$. 
The splitting $(C_1,C_2;S)$ is 
said to be \textit{irreducible} if $(C_1,C_2;S)$ is not reducible. 

\item The splitting $(C_1,C_2;S)$ is said to be \textit{weakly reducible} if there 
are meridian disks $D_i$ $(i=1,2)$ of $C_i$ with $\partial D_1\cap \partial D_2=\emptyset$. 
The splitting $(C_1,C_2;S)$ is said to be \textit{strongly irreducible} if $(C_1,C_2;S)$ 
is not \textit{weakly reducible}. 

\item The splitting $(C_1,C_2;S)$ is said to be \textit{$\partial$-reducible} if there is a disk $D$ 
properly embedded in $M$ such that $D\cap S$ is an essential loop in $S$. Such a disk $D$ 
is called a \textit{$\partial$-reducing disk} for $(C_1,C_2;S)$. 

\item The splitting $(C_1,C_2;S)$ is said to be \textit{stabilized} if there are 
meridian disks $D_i$ $(i=1,2)$ of $C_i$ such that $\partial D_1$ and $\partial D_2$ 
intersect transversely in a single point. Such a pair of disks is called a \textit{cancelling pair of 
disks} for $(C_1,C_2;S)$. 
\end{enumerate}
\end{defi}

\begin{example} 
\rm
Let $(C_1,C_2;S)$ be a Heegaard splitting such that each of $\partial_- C_i$ $(i=1,2)$ consists of 
two $2$-spheres and that $S$ is a $2$-sphere. 
Note that there does not exist an essential disk in $C_i$. Hence $(C_1,C_2;S)$ is strongly irreducible. 
\end{example}

Suppose that $(C_1,C_2;S)$ is stabilized, and let $D_i$ $(i=1,2)$ be disks as in 
$(4)$ of Definition \ref{p}. Note that since $\partial D_1$ intersects $\partial D_2$ 
transversely in a single point, we see that each of $\partial D_i$ $(i=1,2)$ is 
non-separating in $S$ and hence each of $D_i$ $(i=1,2)$ is non-separating in $C_i$. 
Set $C_1'=\mathrm{cl}(C_1\setminus \eta(D_1;C_1))$ and 
$C_2'=C_2\cup \eta(D_1;C_1)$. 
Then each of $C_i'$ $(i=1,2)$ is a compression body with $\partial_+ C_1'=\partial_+ C_2'$ 
(\textit{cf}. (6) of Remark \ref{R1}). Set $S'=\partial_+ C_1'(=\partial_+ C_2')$. 
Then we obtain the Heegaard splitting $(C_1',C_2';S')$ of $M$ with 
$\mathit{genus}(S')=\mathit{genus}(S)-1$. Conversely, $(C_1,C_2;S)$ is obtained 
from $(C_1',C_2';S')$ by adding a trivial handle. We say that $(C_1,C_2;S)$ is obtained 
from $(C_1',C_2';S')$ by \textit{stabilization}.

\begin{obs}\label{S3}
Every reducible Heegaard splitting is weakly reducible. 
\end{obs}

\begin{lem}\label{S2}
Let $(C_1,C_2;S)$ be a Heegaard splitting of $(M;\partial_1 M,\partial_2 M)$ with 
$\mathit{genus}(S)\ge 2$. If $(C_1,C_2;S)$ is stabilized, then $(C_1,C_2;S)$ is reducible. 
\end{lem}

\begin{proof}
Suppose that $(C_1,C_2;S)$ is stabilized, and let $D_i$ $(i=1,2)$ be meridian disks of $C_i$ 
such that $\partial D_1$ intersects $\partial D_2$ transversely in a single point. 
Then $\partial \eta(\partial D_1\cup \partial D_2;S)$ bounds a disk $D_i'$ in $C_i$ for each $i=1$ and $2$. 
In fact, $D'_1$ ($D'_2$ resp.) is obtained from two parallel copies 
of $D_1$ ($D_2$ resp.) by adding a band along 
$\partial D_2\setminus ($the product region between the parallel disks$)$ 
($\partial D_1\setminus ($the product region between the parallel disks$)$ resp.) 
(\textit{cf}. Figure \ref{fig:03}). 

\begin{figure}[htb]\begin{center}
  {\unitlength=1cm
  \begin{picture}(9.8,6)
   \put(0,3.6){\includegraphics[keepaspectratio]{%
   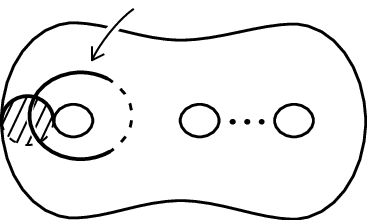}}
   \put(6,3.6){\includegraphics[keepaspectratio]{%
   Fig03p1.eps}}
   \put(0,0){\includegraphics[keepaspectratio]{%
   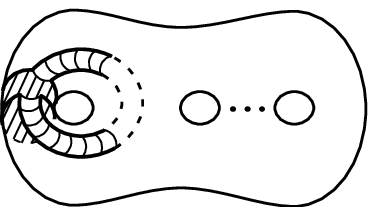}}
   \put(6,0){\includegraphics[keepaspectratio]{%
   Fig03p2.eps}}
   \put(-.6,4.5){$D_1$}
   \put(-.6,1.2){$D'_1$}
   \put(1.5,5.8){$\partial D_2$}
   \put(1.6,2.8){$C_1$}
   \put(4.6,2.8){{\Large $\Downarrow$}}
   \put(5.4,4.5){$D_2$}
   \put(5.4,1.2){$D'_2$}
   \put(7.5,5.8){$\partial D_1$}
   \put(7.7,2.8){$C_2$}
  \end{picture}}
  \caption{}
  \label{fig:03}
\end{center}\end{figure}

Note that $\partial D_1'=\partial D_2'$ cuts $S$ into a torus with a single hole and 
the other surface $S'$. Since $\mathit{genus}(S)\ge 2$, we see that $\mathit{genus}(S')\ge 1$. 
Hence $\partial D_1'=\partial D_2'$ is essential in $S$ and therefore $(C_1,C_2;S)$ is reducible. 
\end{proof}

\begin{defi}\label{trivial}
\rm
Let $(C_1,C_2;S)$ be a Heegaard splitting of $(M;\partial_1 M,\partial_2 M)$. 
\begin{enumerate}
\item 
Suppose that $M\cong S^3$. We call $(C_1,C_2;S)$ a \textit{trivial splitting} 
if both $C_1$ and $C_2$ are $3$-balls. 
\item 
Suppose that $M\not\cong S^3$. We call $(C_1,C_2;S)$ a \textit{trivial splitting} 
if $C_i$ is a trivial handlebody for $i=1$ or $2$. 
\end{enumerate}
\end{defi}

\begin{rem}
\rm
Suppose that $M\not\cong S^3$. 
If $(M;\partial_1 M,\partial_2 M)$ admits a trivial splitting $(C_1,C_2;S)$, 
then it is easy to see that $M$ is a compression body. 
Particularly, if $C_2$ ($C_1$ resp.) is trivial, then $\partial_- M=\partial_1 M$ and 
$\partial_+ M=\partial_2 M$ ($\partial_- M=\partial_2 M$ and 
$\partial_+ M=\partial_1 M$ resp.). 
\end{rem}

\begin{lem}\label{S4}
Let $(C_1,C_2;S)$ be a non-trivial Heegaard splitting of $(M;\partial_1 M,\partial_2 M)$. 
If $(C_1,C_2;S)$ is $\partial$-reducible, then $(C_1,C_2;S)$ is weakly reducible. 
\end{lem}

\begin{proof}
Let $D$ be a $\partial$-reducing disk for $(C_1,C_2;S)$. (Hence $D\cap S$ is an essential loop in $S$.) 
Set $D_1=D\cap C_1$ and $A_2=D\cap C_2$. 
By exchanging subscripts, if necessary, we may suppose that $D_1$ is a meridian disk of $C_1$ and 
$A_2$ is a spanning annulus in $C_2$. Note that $A_2\cap \partial_- C_2$ is an essential loop in 
the component of $\partial_- C_2$ containing $A_2\cap \partial_- C_2$. 
Since $C_2$ is non-trivial, there is a meridian disk of $C_2$. 
It follows from Lemma \ref{spanning} that we can choose a meridian disk $D_2$ of $C_2$ with 
$D_2\cap A_2=\emptyset$. This implies that $D_1\cap D_2=\emptyset$. Hence $(C_1,C_2;S)$ is weakly reducible. 
\end{proof}

\subsection{Haken's theorem}\label{Haken}

In this subsection, we prove the following. 

\begin{thm}\label{S5}
Let $(C_1,C_2;S)$ be a Heegaard splitting of $(M;\partial_1 M,\partial_2 M)$. 
\begin{enumerate}
\item If $M$ is reducible, then $(C_1,C_2;S)$ is reducible or 
$C_i$ is reducible for $i=1$ or $2$. 

\item If $M$ is $\partial$-reducible, then $(C_1,C_2;S)$ is $\partial$-reducible. 
\end{enumerate}
\end{thm}

Note that the statement $(1)$ of Theorem \ref{S5} is called Haken's theorem and proved by Haken 
\cite{Haken}, and the statement $(2)$ of Theorem \ref{S5} is 
proved by Casson and Gordon \cite{CG}. 

We first prove the following proposition, whose statement is weaker than that of 
Theorem \ref{S5}, after showing some lemmas. 

\begin{prop}\label{A}
If $M$ is reducible or $\partial$-reducible, then $(C_1,C_2;S)$ is reducible, 
$\partial$-reducible, or $C_i$ is reducible for $i=1$ or $2$. 
\end{prop}

We give a proof of Proposition \ref{A} by using Otal's idea 
(\textit{cf}. \cite{Otal}) of viewing the Heegaard splittings as a graph in the three dimensional space. 

\begin{figure}[htb]\begin{center}
  {\unitlength=1cm
  \begin{picture}(8.8,1.7)
   \put(0,0.2){\includegraphics[keepaspectratio]{%
   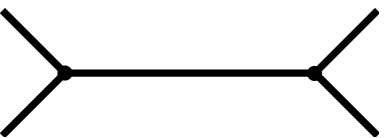}}
   \put(4.5,0){\includegraphics[keepaspectratio]{%
   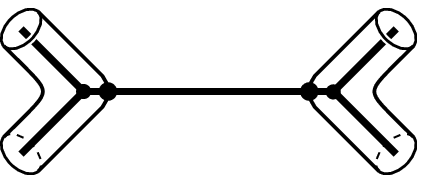}}
   \put(1.9,1.1){$\sigma$}
   \put(0.1,0.8){$p_1$}
   \put(3.5,0.8){$p_2$}
   \put(5.5,.3){$\bar{p}_1$}
   \put(7.5,.3){$\bar{p}_2$}
  \end{picture}}
  \caption{}
  \label{fig:B01}
\end{center}\end{figure}

\bigskip\noindent
\textbf{Edge slides of graphs.}
Let $\Gamma$ be a finite graph in a $3$-manifold $M$. 
Choose an edge $\sigma$ of $\Gamma$. Let $p_1$ and $p_2$ be the vertices of $\Gamma$ incident to $\sigma$. 
Set $\bar{\Gamma}=\Gamma\setminus \sigma$. 
Here, we may suppose that $\sigma\cap \partial \eta(\bar{\Gamma};M)$ consists of two points, say $\bar{p}_1$ 
and $\bar{p}_2$, and that $\mathrm{cl}(\sigma\setminus (p_1\cup p_2))$ consists of $\alpha_0$, $\alpha_1$ and 
$\alpha_2$ with $\partial \alpha_0=\bar{p}_1\cup \bar{p}_2$, $\partial \alpha_1= p_1\cup \bar{p}_1$ and 
$\partial \alpha_2= p_2\cup \bar{p}_2$ (\textit{cf}. Figure \ref{fig:B01}). 

\begin{figure}[htb]\begin{center}
  {\unitlength=1cm
  \begin{picture}(6.1,3.9)
   \put(0,0){\includegraphics[keepaspectratio]{%
   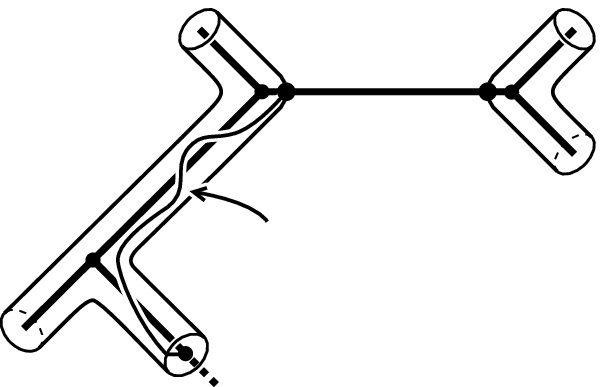}}
   \put(3.8,3.6){$\alpha_0$}
   \put(3.0,3.2){$\overbrace{\hspace{1.9cm}}$}
   \put(2.8,1.5){$\gamma$}
   \put(2.8,2.5){$\bar{p}_1$}
   \put(4.8,2.5){$\bar{p}_2$}
   \put(2.2,0.3){$p'_1$}
  \end{picture}}
  \caption{}
  \label{fig:B02}
\end{center}\end{figure}

Take a path $\gamma$ on $\partial \eta(\bar{\Gamma};M)$ with $\partial \gamma\ni \bar{p}_1$. 
Let $\bar{\sigma}$ be an arc obtained from $\gamma\cup \alpha_0\cup \alpha_2$ by adding a `straight short arc' 
in $\eta(\bar{\Gamma};M)$ connecting the endpoint of $\gamma$ other than $\bar{p}_1$ and a point $p'_1$ 
in the interior of an edge of $\bar{\Gamma}$ (\textit{cf}. Figure \ref{fig:B02}). 
Let $\Gamma'$ be a graph obtained from $\bar{\Gamma}\cup \bar{\sigma}$ by adding $p'_1$ as a vertex. 
Then we say that $\Gamma'$ is obtained from $\Gamma$ by an \textit{edge slide} on $\sigma$. 

\begin{figure}[htb]\begin{center}
  {\unitlength=1cm
  \begin{picture}(9.4,4.9)
   \put(0,0){\includegraphics[keepaspectratio]{%
             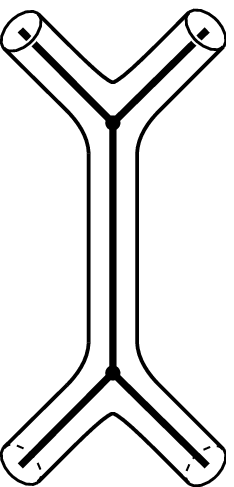}}
   \put(4.5,1.3){\includegraphics[keepaspectratio]{%
             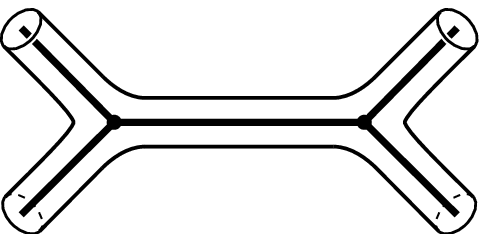}}
   \put(2.8,2.3){$\longleftrightarrow$}
  \end{picture}}
  \caption{}
  \label{fig:B03}
\end{center}\end{figure}

If $p_1$ is a trivalent vertex, then it is natural for us not to regard $p_1$ as a vertex of $\Gamma'$. 
Particularly, the deformation of $\Gamma$ which is depicted as in Figure \ref{fig:B03} is realized by an edge 
slide and an isotopy. This deformation is called a \textit{Whitehead move}. 

\bigskip\noindent
\textbf{A Proof of Proposition \ref{A}.}
Let $\Sigma$ be a spine of $C_1$. Note that $\eta(\partial_- C_1\cup \Sigma;M)$ is obtained from 
regular neighborhoods of $\partial_- C_1$ and the vertices of $\Sigma$ by attaching $1$-handles 
corresponding to the edges of $\Sigma$. Set $\Sigma_{\eta}=\eta(\Sigma;M)$. 
The notation $h^0_v$, called a \textit{vertex} of $\Sigma_{\eta}$, means a regular neighborhood 
of a vertex $v$ of $\Sigma$. 
Also, the notation $h^1_{\sigma}$, called an \textit{edge} of $\Sigma_{\eta}$, means a $1$-handle 
corresponding to an edge $\sigma$ of $\Sigma$. 
Let $\Delta=D_1\cup \dots \cup D_k$ be a minimal complete meridian system of $C_2$. 

Let $P$ be a reducing $2$-sphere or a $\partial$-reducing disk of $M$. 
If $P$ is a $\partial$-reducing disk, we may assume that $\partial P\subset \partial_- C_2$ 
by changing subscripts. 
We may assume that $P$ intersects $\Sigma$ and $\Delta$ 
transversely. 
Set $\Gamma=P\cap (\Sigma_{\eta}\cup \Delta)$. We note that 
$\Gamma$ is a union of disks $P\cap \Sigma_{\eta}$ and a union of arcs and loops $P\cap \Delta$ in $P$. 
We choose $P$, $\Sigma$ and $\Delta$ so that the pair 
$(|P\cap \Sigma|,|P\cap \Delta|)$ is minimal with respect to lexicographic order. 

\begin{lem}\label{A1}
Each component of $P\cap \Delta$ is an arc. 
\end{lem}

\begin{proof}
For some disk component, say $D_1$, of $\Delta$, suppose that $P\cap D_1$ has a loop 
component. Let $\alpha$ be a loop component of $P\cap D_1$ which is innermost in 
$D_1$, and let $\delta_{\alpha}$ be an innermost disk for $\alpha$. 
Let $\delta'_{\alpha}$ be a disk in $P$ with $\partial \delta'_{\alpha}=\alpha$. Set 
$P'=(P\setminus \delta'_{\alpha})\cup \delta_{\alpha}$ if $P$ is a $\partial$-reducing disk, or set 
$P'=(P\setminus \delta'_{\alpha})\cup \delta_{\alpha}$ and 
$P''=\delta'_{\alpha}\cup \delta_{\alpha}$ if $P$ is a reducing $2$-sphere. If $P$ is a $\partial$-reducing disk, 
then $P'$ is also a $\partial$-reducing disk. If $P$ is a reducing $2$-sphere, 
then either $P'$ or $P''$, say $P'$, is a reducing $2$-sphere. 
Moreover, we can isotope $P'$ so that $(|P'\cap \Sigma|,|P'\cap \Delta|)<
(|P\cap \Sigma|,|P\cap \Delta|)$. This contradicts the minimality of 
$(|P\cap \Sigma|,|P\cap \Delta|)$. 
\end{proof}

By Lemma \ref{A1}, we can regard $\Gamma$ as a graph in $P$ which consists of 
\textit{fat-vertices} $P\cap \Sigma_{\eta}$ and edges $P\cap \Delta$. 
An edge of the graph $\Gamma$ is called a \textit{loop} if the edge joins a fat-vertex of $\Gamma$ to itself, 
and a loop is said to be \textit{inessential} if the loop cuts off a 
disk from $\mathrm{cl}(P\setminus \Sigma_{\eta})$ whose interior is disjoint from 
$\Gamma \cap \Sigma_{\eta}$. 

\begin{lem}\label{A2}
$\Gamma$ does not contain an inessential loop. 
\end{lem}

\begin{proof}
Suppose that $\Gamma$ contains an inessential loop $\mu$. Then $\mu$ cuts off a disk $\delta_{\mu}$ 
from $\mathrm{cl}(P\setminus \Sigma_{\eta})$ such that the interior of $\delta_{\mu}$ is disjoint from 
$\Gamma \cap \Sigma_{\eta}$ (\textit{cf}. Figure \ref{fig:07}). 

\begin{figure}[htb]\begin{center}
  {\unitlength=1cm
  \begin{picture}(5,3.6)
   \put(0,0){\includegraphics[keepaspectratio]{%
   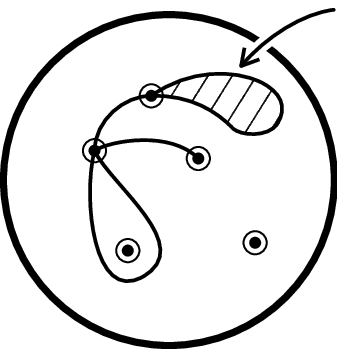}}
  \put(3.6,3.4){an inessential monogon}
  \end{picture}}
  \caption{}
  \label{fig:07}
\end{center}\end{figure}

We may assume that 
$\delta_{\mu}\cap \Delta=\delta_{\mu}\cap D_1$. Then $\mu$ cuts $D_1$ into two disks $D_1'$ and 
$D_1''$ (\textit{cf}. Figure \ref{fig:08}). 

\begin{figure}[htb]\begin{center}
  {\unitlength=1cm
  \begin{picture}(12,5)
   \put(0,0){\includegraphics[keepaspectratio]{%
   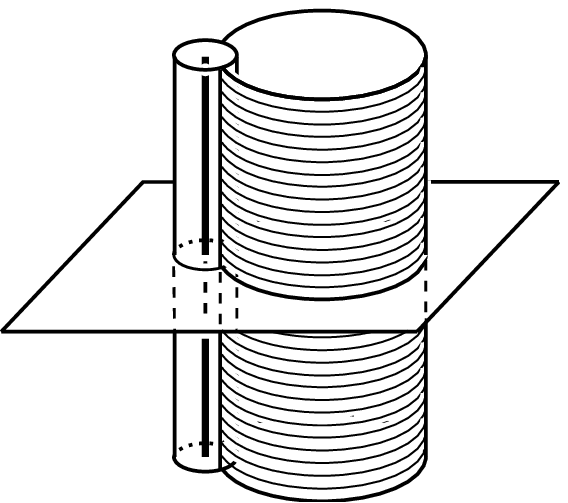}}
   \put(6.3,0){\includegraphics[keepaspectratio]{%
   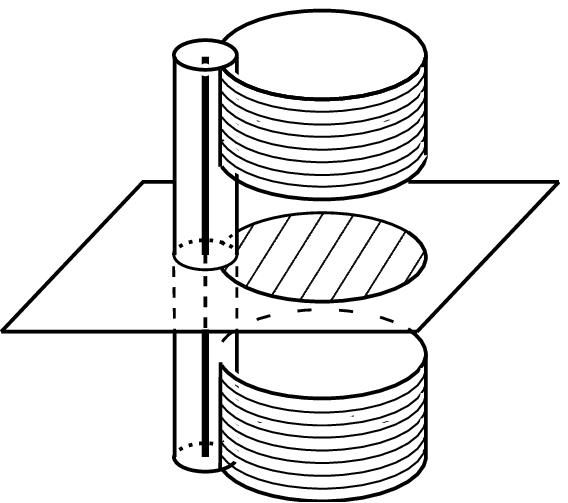}}
   \put(-.3,1.7){$P$}
   \put(6,1.7){$P$}
   \put(5.7,2.4){$\Longrightarrow$}
   \put(10.7,3.8){$D'_1$}
   \put(10.7,2.6){$\delta_\mu$}   
   \put(10.7,1.1){$D''_1$}
  \end{picture}}
  \caption{}
  \label{fig:08}
\end{center}\end{figure}

Let $C_2'$ be 
the component, which is obtained by cutting $C_2$ along $\Delta$, such that $C_2'$ 
contains $\delta_{\mu}$. Let $D_1^+$ be the copy of $D_1$ in $C_2'$ with $D_1^+\cap \delta_{\mu}\ne 
\emptyset$ and $D_1^-$ the other copy of $D_1$. Note that 
$C_2'$ is a $3$-ball or a $($a component of $\partial_- C_2$)$\times [0,1]$. This shows that there is a disk 
$\delta'_{\mu}$ in $\partial_+ C_2'$ such that $\partial \delta_{\mu}=\partial \delta'_{\mu}$ and 
$\partial \delta_{\mu}\cup \partial \delta'_{\mu}$ bounds a $3$-ball in $C_2'$. 
Note that $\delta'_{\mu}\cap D_1^+\ne \emptyset$. By changing superscripts, if necessary, we may assume that 
$\delta'_{\mu}\supset D_1'$ (\textit{cf}. Figure \ref{fig:09}). 

\begin{figure}[htb]\begin{center}
  {\unitlength=1cm
  \begin{picture}(6.7,5)
   \put(0,0){\includegraphics[height=4cm,keepaspectratio]{%
   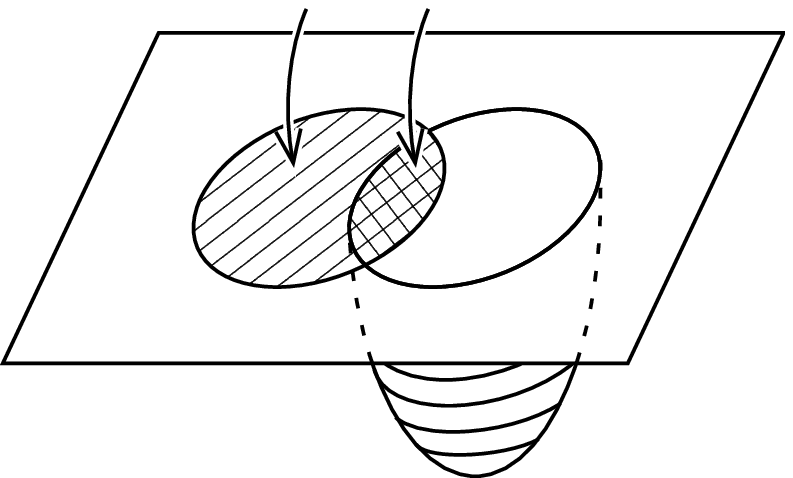}}
   \put(.3,1.2){$\partial_+ C'_2$}
   \put(2.9,4.9){$D^+_1$}
   \put(2.4,4.5){$\overbrace{\hspace*{40pt}}$}
   \put(2.4,4.2){$D''_1$}
   \put(3.4,4.2){$D'_1$}
   \put(2.9,1.66){$\underbrace{\hspace*{68pt}}$}
   \put(4.1,1.16){$\delta'_\mu$}
   \put(4.7,.3){$\delta_\mu$}
  \end{picture}}
  \caption{}
  \label{fig:09}
\end{center}\end{figure}

Set $D_0=\delta_{\mu}\cup D_1'$ if 
$\delta'_{\mu}\cap D_1^-\ne \emptyset$, and $D_0=\delta_{\mu}\cup D_1''$ if 
$\delta'_{\mu}\cap D_1^-= \emptyset$. 
We may regard $D_0$ as a disk properly embedded in $C_2$. Set 
$\Delta'=D_0\cup D_2\cup \dots \cup D_k$. Then we see that $\Delta'$ is a minimal complete 
meridian system of $C_2$. We can further isotope $D_0$ slightly so that 
$|P\cap \Delta'|<|P\cap \Delta|$. 
This contradicts the minimality of $(|P\cap \Sigma|,|P\cap \Delta|)$. 
\end{proof}

\begin{figure}[htb]\begin{center}
  {\unitlength=1cm
  \begin{picture}(5,3.4)
   \put(0,0){\includegraphics[keepaspectratio]{%
   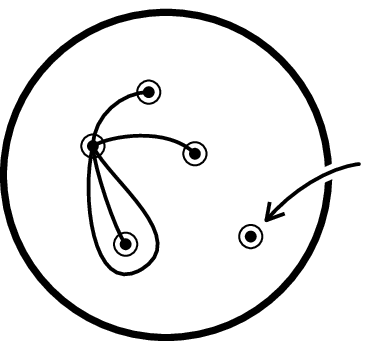}}
   \put(3.9,1.76){an isolated fat-vertex}
  \end{picture}}
  \caption{}
  \label{fig:10}
\end{center}\end{figure}

A fat-vertex of $\Gamma$ is said to be \textit{isolated} if there are no edges of 
$\Gamma$ adjacent to the fat-vertex (\textit{cf}. Figure \ref{fig:10}). 

\begin{lem}\label{A3}
If $\Gamma$ has an isolated fat-vertex, then $(C_1,C_2;S)$ is reducible or $\partial$-reducible. 
\end{lem}

\begin{proof}
Suppose that there is an isolated fat-vertex $D_v$ of $\Gamma$. 
Recall that $D_v$ is a component of $P\cap \Sigma_{\eta}$ which is a 
meridian disk of $C_1$. 
Note that $D_v$ is disjoint from $\Delta$ (\textit{cf}. Figure \ref{fig:10.1}). 

\begin{figure}[htb]\begin{center}
  {\unitlength=1cm
  \begin{picture}(5.5,5)
   \put(0,0){\includegraphics[keepaspectratio]{
   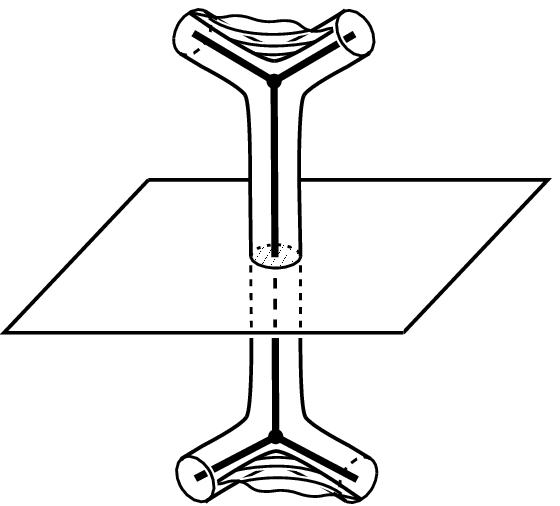}}
   \put(.5,1.9){$P$}
   \put(3.2,2.4){$D_v$}
  \end{picture}}
  \caption{}
  \label{fig:10.1}
\end{center}\end{figure}

Let $C_2'$ be the component obtained by cutting $C_2$ along 
$\Delta$ such that $\partial C_2'$ contains $\partial D_v$. If $\partial D_v$ 
bounds a disk $D'_v$ in $C_2'$, then $D_v$ and $D'_v$ indicates the reducibility of 
$(C_1,C_2;S)$. Otherwise, $C_2'$ is a $($a closed orientable surface$)\times [0,1]$, and 
$\partial D_v$ is a boundary component of a spanning annulus in $C_2'$ (and hence $C_2$). 
Hence we see that $(C_1,C_2;S)$ is $\partial$-reducible. 
\end{proof}

\begin{lem}\label{A4}
Suppose that no fat-vertices of $\Gamma$ are isolated. 
Then each fat-vertex of $\Gamma$ is a base of a loop. 
\end{lem}

\begin{proof}
Suppose that there is a fat-vertex $D_w$ of $\Gamma$ which is not a base of a loop. 
Since no fat-vartices of $\Gamma$ are isolated, there is an edge of $\Gamma$ adjacent to $D_w$. 
Let $\sigma$ be the edge of $\Sigma$ with 
$h^1_{\sigma}\supset D_w$. (Recall that $h^1_{\sigma}$ is a $1$-handle of $\Sigma_{\eta}$ 
corresponding to $\sigma$.) Let $D$ be a component of $\Delta$ with $\partial D
\cap h^1_{\sigma} \ne \emptyset$. Let 
$C_w$ be a union of the arc components of $D \cap P$ which are adjacent to $D_w$. 
Let $\gamma$ be an arc component of $C_w$ which is outermost among the components of $C_w$. 
We call such an arc $\gamma$ an \textit{outermost edge for $D_w$ of $\Gamma$}. 
Let $\delta_{\gamma}\subset D$ be a disk obtained by cutting $D$ along $\gamma$ whose interior 
is disjoint from the edges incident to $D_w$. 
We call such a disk $\delta_{\gamma}$ an \textit{outermost disk for} $(D_w,\gamma)$. 
(Note that $\delta_{\gamma}$ may intersects 
$P$ transversely (\textit{cf}. Figure \ref{fig:11}).) Let $D_{w'}(\ne D_w)$ 
be the fat-vertex of $\Gamma$ attached to $\gamma$. Then we have the following three cases. 

\begin{figure}[htb]\begin{center}
  {\unitlength=1cm
  \begin{picture}(9.3,3.6)
   \put(0,0){\includegraphics[keepaspectratio]{%
   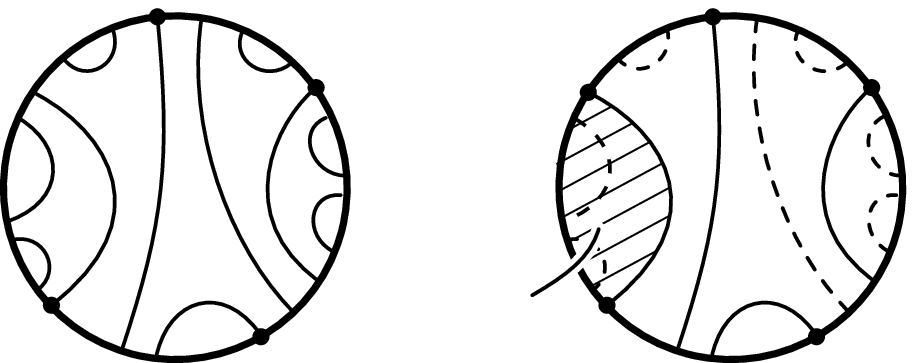}}
   \put(.5,-.3){$D_1$}
   \put(.1,.38){$D_{w}$}
   \put(1.5,3.73){$D_{w}$}
   \put(3.4,2.8){$D_{w}$}
   \put(2.7,-.05){$D_{w}$}
   \put(5.7,.38){$D_{w}$}
   \put(7.14,3.73){$D_{w}$}
   \put(9.04,2.8){$D_{w}$}
   \put(8.34,-.05){$D_{w}$}
   \put(6.8,2.2){$\gamma$}
   \put(5,.5){$\delta_{\gamma}$}
  \end{picture}}
  \caption{}
  \label{fig:11}
\end{center}\end{figure}

\medskip
\noindent\textit{Case 1.}\ \ \ \ 
$(\partial \delta_{\gamma}\setminus \gamma) \subseteq (h^1_{\sigma}\cap D)$. 

\medskip
In this case, we can isotope $\sigma$ along $\delta_{\gamma}$ to reduce $|P\cap \Sigma|$ 
(\textit{cf}. Figure \ref{fig:12}). 

\begin{figure}[htb]\begin{center}
  {\unitlength=1cm
  \begin{picture}(12,2.8)
   \put(-.6,0.15){\includegraphics[keepaspectratio]{%
   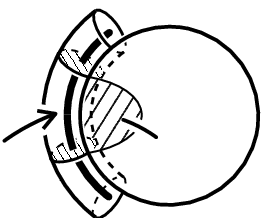}}
   \put(3.36,0){\includegraphics[keepaspectratio]{%
   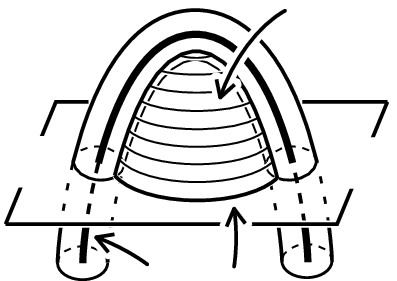}}
   \put(8.2,0){\includegraphics[keepaspectratio]{%
   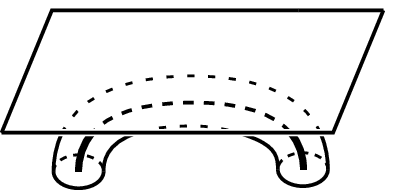}}
   \put(-.9,.7){$\sigma$}
   \put(-.5,.4){$D_w$}
   \put(-.65,1.8){$D_{w'}$}
   \put(.7,1.5){$\gamma$}
   \put(1.1,.7){$\delta_{\gamma}$}
   \put(3,.6){$P$}
   \put(4.8,0){$\sigma$}
   \put(3.5,1.1){$D_w$}
   \put(6.63,1.1){$D_{w'}$}
   \put(5.65,-.1){$\gamma$}
   \put(6.26,2.8){$\delta_{\gamma}$}
   \put(7.48,1.16){$\Longrightarrow$}
  \end{picture}}
  \caption{}
  \label{fig:12}
\end{center}\end{figure}

\medskip
\noindent\textit{Case 2.}\ \ \ \ 
$(\partial \delta_{\gamma}\setminus \gamma) \not\subseteq (h^1_{\sigma}\cap D)$ and 
$D_{w'}\not\subset (h^1_{\sigma}\cap D)$. 

\medskip
Let $p$ be the vertex of $\Sigma$ such that $p\cap \sigma\ne \emptyset$ and 
$h^0_p\cap \delta_{\gamma}\ne \emptyset$. 
Let $\beta$ be the component of $\mathrm{cl}(\sigma\setminus D_w)$ which 
satisfies $\beta\cap p\ne \emptyset$. 
Then we can slide $\beta$ along $\delta_{\gamma}$ so that $\beta$ contains $\gamma$ 
(\textit{cf}. Figure \ref{fig:13}). 
We can further isotope $\beta$ slightly to reduce $|P\cap \Sigma|$, a contradiction. 

\begin{figure}[htb]\begin{center}
  {\unitlength=1cm
  \begin{picture}(12,3.4)
   \put(-.6,0.7){\includegraphics[keepaspectratio]{%
   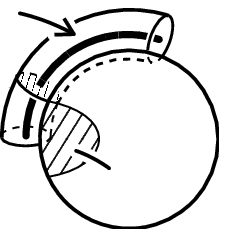}}
   \put(3.1,0){\includegraphics[keepaspectratio]{%
   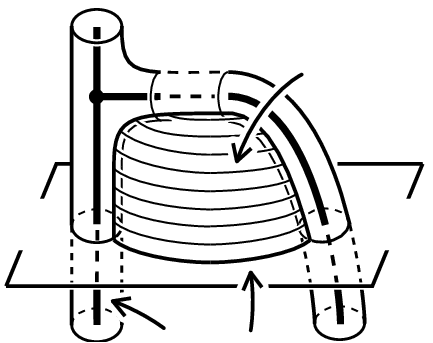}}
   \put(7.94,0){\includegraphics[keepaspectratio]{%
   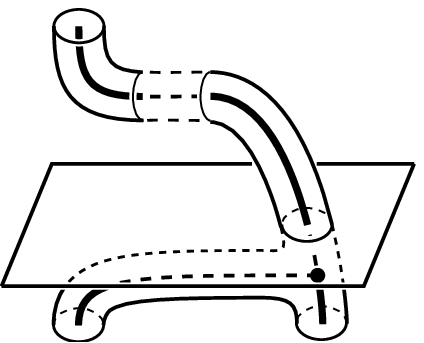}}
   \put(-.7,2.9){$\sigma$}
   \put(.4,1.9){$\gamma$}
   \put(.6,1.1){$\delta_{\gamma}$}
   \put(-1.1,2.2){$D_w$}
%
   \put(2.8,.6){$P$}
   \put(4.8,0){$\sigma$}
   \put(3.5,2.5){$p$}
   \put(3.2,1.16){$D_{w}$}
   \put(6.68,1.16){$D_{w'}$}
   \put(5.55,-.14){$\gamma$}
   \put(6.26,2.8){$\delta_{\gamma}$}
   \put(7.4,1.16){$\Longrightarrow$}
  \end{picture}}
  \caption{}
  \label{fig:13}
\end{center}\end{figure}

\medskip
\noindent\textit{Case 3.}\ \ \ \ 
$(\partial \delta_{\gamma}\setminus \gamma)\not\subseteq (h^1_{\sigma}\cap D)$ and 
$D_{w'}\subset (h^1_{\sigma}\cap D)$. 

\medskip
Let $p$ and $p'$ be the endpoints of $\sigma$. 
Let $\beta$ and $\beta'$ be the components of $\mathrm{cl}(\sigma\setminus (D_w\cup D_{w'}))$ 
which satisfy $p\cap \beta\ne \emptyset$ and $p'\cap \beta'\ne \emptyset$. 
Suppose first that $p\ne p'$. Then we can slide $\beta$ along $\delta_{\gamma}$ 
so that $\beta$ contains $\gamma$ (\textit{cf}. Figure \ref{fig:14}). 
We can further isotope $\beta$ slightly to reduce $|P\cap \Sigma|$, a contradiction. 

\begin{figure}[htb]\begin{center}
  {\unitlength=1cm
  \begin{picture}(12,4.3)
   \put(-.6,0.15){\includegraphics[keepaspectratio]{%
   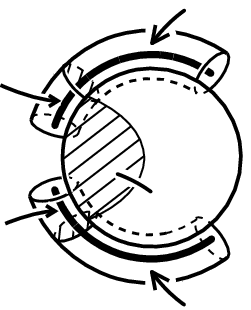}}
   \put(3.46,0){\includegraphics[keepaspectratio]{%
   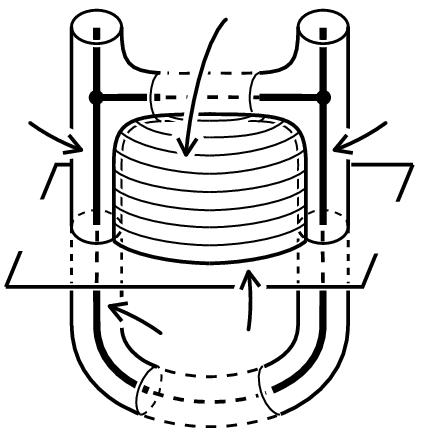}}
   \put(8.26,0){\includegraphics[keepaspectratio]{%
   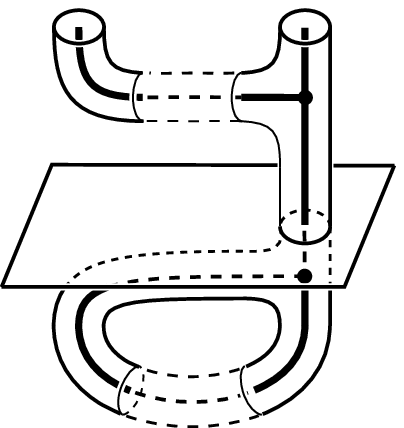}}
   \put(1.32,0.12){$\sigma$}
   \put(1.32,3.15){$\sigma$}
   \put(-.5,2.7){$D_{w}$}
   \put(-.7,.6){$D_{w'}$}
   \put(.8,2.1){$\gamma$}
   \put(1.0,1.2){$\delta_{\gamma}$}
   \put(-0.5,2.0){$p$}
   \put(-0.5,1.4){$p'$}
   \put(-.9,2.3){$\beta$}
   \put(-.9,.9){$\beta'$}
   \put(3.1,1.4){$P$}
   \put(5.1,.8){$\sigma$}
   \put(3.6,1.96){$D_w$}
   \put(7,1.96){$D_{w'}$}
   \put(3.9,3.26){$p$}
   \put(7.04,3.26){$p'$}
   \put(3.5,3.0){$\beta$}
   \put(7.4,3.0){$\beta'$}
   \put(5.85,.8){$\gamma$}
   \put(5.8,4.2){$\delta_{\gamma}$}
   \put(7.7,2){$\Longrightarrow$}
  \end{picture}}
  \caption{}
  \label{fig:14}
\end{center}\end{figure}

Suppose next that $p=p'$. In this case, we perform the following operation which 
is called a \textit{broken edge slide}. 

\begin{figure}[htb]\begin{center}
  {\unitlength=1cm
  \begin{picture}(10.8,4.5)
   \put(0,0){\includegraphics[keepaspectratio,width=4.8cm]{%
   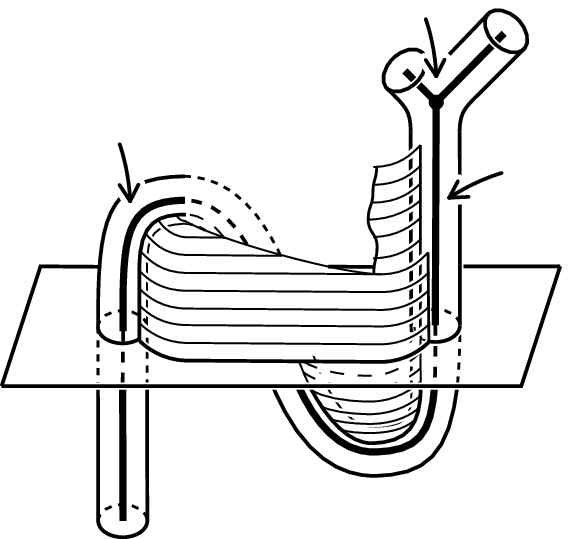}}
   \put(6,0){\includegraphics[keepaspectratio,width=4.8cm]{%
   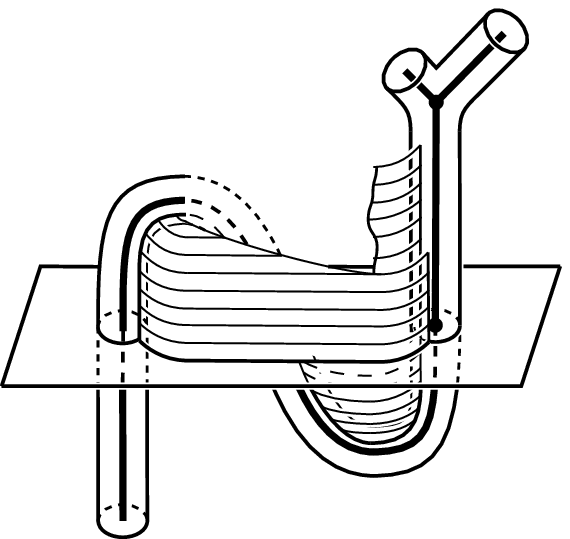}}
   \put(-.3,1.4){$P$}
   \put(.2,1.8){$D_{w}$}
   \put(3.95,1.8){$D_{w'}$}
   \put(.7,3.5){$\beta$}
   \put(4.4,3.2){$\beta'$}
   \put(3.4,4.6){$p$}
   \put(10,1.8){$w'$}
   \put(5.2,2.1){$\Longrightarrow$}
  \end{picture}}
  \caption{}
  \label{fig:15.1}
\end{center}\end{figure}

We first add $w'=D_{w'}\cap \Sigma$ as a vertex of 
$\Sigma$. Then $w'$ cuts $\sigma$ into two edges $\beta'$ and $\mathrm{cl}(\sigma\setminus \beta')$. 
Since $\gamma$ is an outermost edge for $D_w$ of $\Gamma$, we see that 
$\beta'\subset \beta$ (\textit{cf}. Figure \ref{fig:15.1}). Hence we can slide 
$\mathrm{cl}(\beta\setminus \beta')$ along $\delta_{\gamma}$ 
so that $\mathrm{cl}(\beta\setminus \beta')$ contains $\gamma$. 
We now remove the verterx $w'$ of $\Sigma$, that is, we regard a union of $\beta'$ and 
$\mathrm{cl}(\sigma\setminus \beta')$ as an edge of $\Sigma$ again. 
Then we can isotope $\mathrm{cl}(\sigma\setminus \beta')$ slightly to reduce 
$|P\cap \Sigma|$, a contradiction (\textit{cf}. Figure \ref{fig:15.2}). 

\begin{figure}[htb]\begin{center}
  {\unitlength=1cm
  \begin{picture}(10.4,3.5)
   \put(0,0){\includegraphics[keepaspectratio,width=4.8cm]{%
   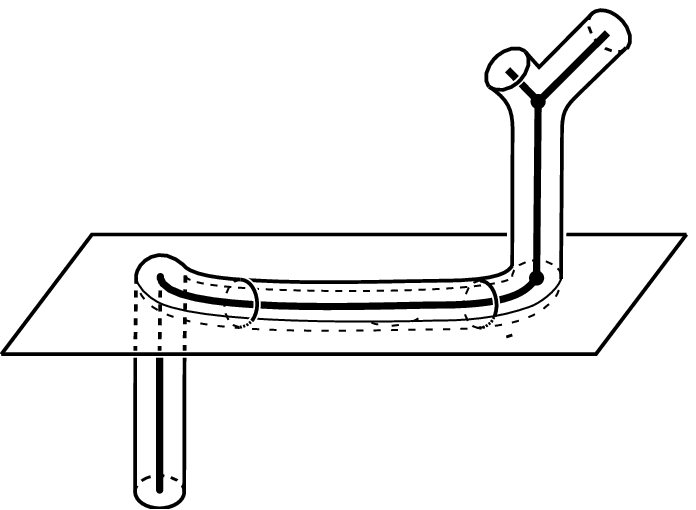}}
   \put(5.6,0){\includegraphics[keepaspectratio,width=4.8cm]{%
   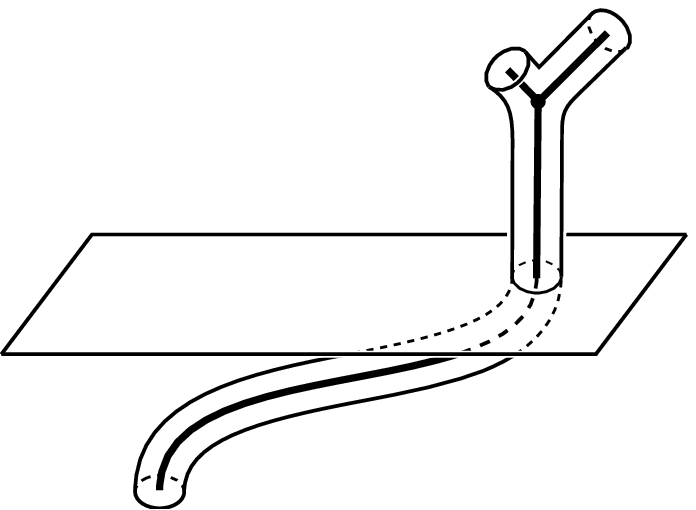}}
   \put(4.0,2.7){$p$}
   \put(3.9,1.4){$w'$}
   \put(4.92,1.5){$\Longrightarrow$}
  \end{picture}}
  \caption{}
  \label{fig:15.2}
\end{center}\end{figure}

\end{proof}

\begin{proof}[Proof of Proposition \ref{A}]
By Lemma \ref{A3}, if there is an isolated fat-vertex of $\Gamma$, then we have the conclusion of 
Proposition \ref{A}. Hence we suppose that no fat-vertices of $\Gamma$ are isolated. 
Then it follows from Lemma \ref{A4} that each fat-vertex of $\Gamma$ is a base of a loop. 
Let $\mu$ be a loop which is innermost in $P$. Then $\mu$ 
cuts a disk $\delta_{\mu}$ from $\mathrm{cl}(P\setminus \Sigma_{\eta})$. 
Since $\mu$ is essential (\textit{cf}. Lemma \ref{A2}), 
we see that $\delta_{\mu}$ contains a fat-vertex of 
$\Gamma$. But since $\mu$ is innermost, such a fat-vertex is not a base of any 
loop. Hence such a fat-vertex is isolated, a contradiction. 
This completes the proof of Proposition \ref{A}.
\end{proof}

\begin{proof}[Proof of (1) in Theorem \ref{S5}]
Suppose that $M$ is reducible. Then by Proposition \ref{A}, we see that 
$(C_1,C_2;S)$ is reducible or $\partial$-reducible, or $C_i$ is reducible for 
$i=1$ or $2$. If $(C_1,C_2;S)$ is reducible or $C_i$ is reducible for $i=1$ or $2$, 
then we are done. So we may assume that $C_1$ and $C_2$ are irreducible and that 
$(C_1,C_2;S)$ is $\partial$-reducible. 
By induction on the genus of the Heegaard surface $S$, we prove that $(C_1,C_2;S)$ is 
reducible. 

Suppose that $\mathit{genus}(S)=0$. Since $C_i$ $(i=1,2)$ are irreducible, we see that 
each of $C_i$ $(i=1,2)$ is a $3$-ball. Hence $M$ is the $3$-sphere and therefore 
$M$ is irreducible, a contradiction. So we may assume that $\mathit{genus}(S)>0$. 
Let $P$ be a $\partial$-reducing disk of $M$ with $|P\cap S|=1$. By changing subs
cripts, 
if necessary, we may assume that $P\cap C_1=D$ is a disk and $P\cap C_2=A$ is a spanning annulus. 

Suppose that $\mathit{genus}(S)=1$. 
Since $C_i$ $(i=1,2)$ are irreducible, we see that $\partial C_i$ contain no $2$-sphere components. 
Since $C_1$ contains an essential disk $D$, we see that $C_1\cong D^2\times S^1$. 
Since $C_2$ contains a spanning annulus $A$, we see that $C_2\cong T^2\times [0,1]$. 
It follows that $M\cong D^2\times S^1$ and hence $M$ is irreducible, a contradiction. 

Suppose that $\mathit{genus}(S)>1$. Let $C_1'$ ($C_2'$ resp.) be the manifold obtained from 
$C_1$ ($C_2$ resp.) by cutting along $D$ ($A$ resp.), and let $A^+$ and $A^-$ be copies of $A$ in 
$\partial C_2'$. Then we see that $C_1'$ consists of either a compression body or a union of two compression bodies 
(\textit{cf}. $(6)$ of Remark \ref{R1}). 
Let $C_2''$ be the manifold obtained from $C_2'$ by attaching $2$-handles along $A^+$ and $A^-$. 
It follows from Remark \ref{R0} that $C_2''$ consists of either a compression body or a union of 
two compression bodies. 

\medskip
Suppose that $C_1'$ consists of a compression body. This implies that $C_2''$ consists of 
a compression body (\textit{cf}. Figure \ref{fig:16}). 
We can naturally obtain a homeomorphism $\partial_+ C_1'\to \partial_+ C_2''$ from 
the homeomorphism $\partial_+ C_1\to \partial_+ C_2$. Set $\partial_+ C_1'=\partial_+ C_2''=S'$. 
Then $(C_1',C_2'';S')$ is a Heegaard splitting of the $3$-manifolds $M'$ obtained by cutting $M$ along $P$. 
Note that $\mathit{genus}(S')=\mathit{genus}(S)-1$. 
Moreover, by using innermost disk arguments, we see that $M'$ is also reducible. 

\begin{figure}[htb]\begin{center}
  {\unitlength=1cm
  \begin{picture}(12,4.3)
   \put(0,2.6){\includegraphics[keepaspectratio]{%
   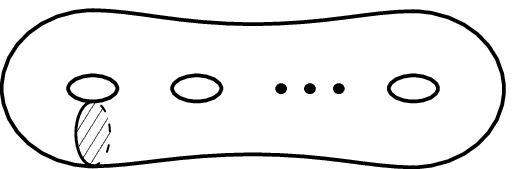}}
   \put(6.8,2.6){\includegraphics[keepaspectratio]{%
   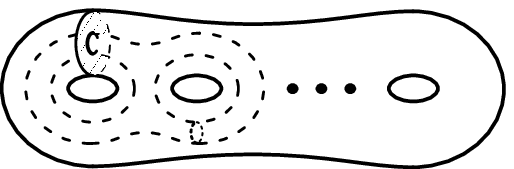}}
   \put(.38,0){\includegraphics[keepaspectratio]{%
   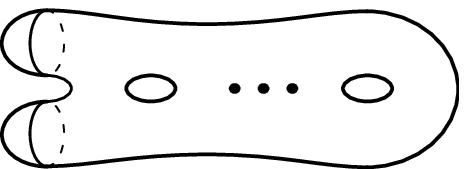}}
   \put(7.18,0){\includegraphics[keepaspectratio]{%
   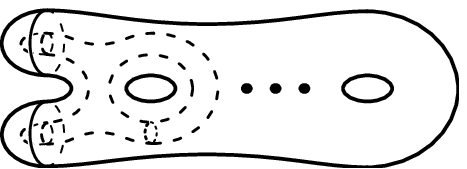}}
   \put(.8,2.3){$D$}
   \put(5.85,3.4){{\Large $\cup$}}
   \put(5.9,3){$S$}
   \put(5,2.6){$C_1$}
   \put(6.54,2.6){$C_2$}
   \put(7.6,4.4){$A$}
   \put(5.9,1.9){{\Large $\Downarrow$}}
   \put(5.9,.9){{\Large $\cup$}}
   \put(5.95,.5){$S'$}
   \put(5,0){$C'_1$}
   \put(6.54,0){$C'_2$}
  \end{picture}}
  \caption{}
  \label{fig:16}
\end{center}\end{figure}

\medskip
\noindent\textit{Claim.}\ \ \ \ 
\begin{enumerate}
\item If $C_1'$ is reducible, then $C_1$ is reducible. 
\item If $C_2''$ is reducible, then one of the following holds. 
\begin{enumerate}
\item $C_2$ is reducible. 
\item The component of $\partial_- C_2$ intersecting $A$ is a torus, say $T$. 
\end{enumerate}
\end{enumerate}

\medskip
\noindent\textit{Proof.}\ \ \ \ Exercise \ref{ES5}. 

\medskip
Recall that we assume that $C_i$ $(i=1,2)$ are irreducible. 
Hence it follows from $(1)$ of the claim that $C_1'$ is irreducible. 
Also it follows from $(2)$ of the claim  that either (I) $C_2''$ is irreducible or (II) 
$C_2''$ is reducible and the condition (b) of $(2)$ in Claim 1 holds. 

Suppose that the condition (I) holds. Then by induction on the genus of a Heegaard surface, 
$(C_1',C_2'';S')$ is reducible, i.e., there are meridian disks $D_1$ and $D_2$ of $C_1'$ and $C_2''$ 
respectively with $\partial D_1=\partial D_2$. Note that this implies that $C_i$ $(i=1,2)$ are non-trivial. 
Let $\alpha^+$ and $\alpha^-$ be the co-cores of the $2$-handles attached to $C_2''$. Then we 
see that $\alpha^+\cup \alpha^-$ is vertical in $C_2''$. 
It follows from Lemma \ref{vertical} that we may assume that $D_2\cap(\alpha^+\cup \alpha^-)=\emptyset$, i.e., 
$D_2$ is disjoint from the $2$-handles. Hence the pair of disks $D_1$ and $D_2$ survives when we restore 
$C_1$ and $C_2$ from $C_1'$ and $C_2''$ respectively. This implies that $(C_1,C_2;S)$ is reducible 
and hence we obtain the conclusion $(1)$ of Theorem \ref{S5}. 

\medskip
Suppose that the condition (II) holds. Then it follows from $(2)$ of Remark \ref{R0} that 
there is a separating disk $E_2$ in $C_2$ such that $E_2$ is disjoint from $A$ and that 
$E_2$ cuts off $T^2\times [0,1]$ from $C_2$ with $T^2\times \{0\}=T$. 
Let $\ell$ be a loop in $S\cap (T^2\times \{1\})$ which intersects $A\cap S(=\partial A\cap S=\partial D)$ 
in a single point. 
Let $E_1$ be a disk properly embedded in $C_1$ which is obtained from two parallel copies 
of $D$ by adding a band along $\ell\setminus ($the product region between the parallel disks$)$. 
Since $\mathit{genus}(S)>1$, we see that $E_1$ is a separating meridian disk of $C_1$. 
Since $\partial E_1$ is isotopic to $\partial E_2$, we see that 
$(C_1,C_2;S)$ is reducible. Hence we obtain the conclusion $(1)$ of Theorem \ref{S5}. 

\medskip
The case that $C_1'$ is a union of two compression bodies is treated analogously, and we leave 
the proof for this case to the reader (Exercise \ref{ES5b}). 
\end{proof}

\begin{exercise}\label{ES5}
\rm
Show the claim in the proof of Theorem \ref{S5}. 
\end{exercise}

\begin{exercise}\label{ES5b}
\rm
Prove that the conclusion $(1)$ of Theorem \ref{S5} holds 
in case that $C_1'$ consists of two compression bodies. 
\end{exercise}

\begin{proof}[Proof of $(2)$ in Theorem \ref{S5}]
Suppose that $M$ is $\partial$-reducible. If $(C_1,C_2;S)$ is $\partial$-reducible, 
then we are done. Let $\widehat{C}_i$ be the compression body obtained by attaching $3$-balls 
to the $2$-sphere boundary components of $C_i$ $(i=1,2)$. 
Set $\widehat{M}=\widehat{C}_1\cup \widehat{C}_2$. Then $\widehat{M}$ is also $\partial$-reducible. 
Then it follows from $(1)$ of Remark \ref{R1} and Proposition \ref{A} that 
$(\widehat{C}_1,\widehat{C}_2;S)$ is reducible or $\partial$-reducible. 
If $(\widehat{C}_1,\widehat{C}_2;S)$ is $\partial$-reducible, then 
we see that $(C_1,C_2;S)$ is also $\partial$-reducible. 
Hence we may assume that $(\widehat{C}_1,\widehat{C}_2;S)$ is reducible. 
By induction on the genus of a Heegaard 
surface, we prove that $(C_1,C_2;S)$ is $\partial$-reducible. Let $P'$ be a reducing 
$2$-sphere of $\widehat{M}$ with $|P'\cap S|=1$. 
For each $i=1$ and $2$, set $D_i=P'\cap \widehat{C}_i$, and let $\widehat{C}_i'$ be the manifold 
obtained by cutting $\widehat{C}_i$ along $D_i$, and let $D_i^+$ and $D_i^-$ be copies 
of $D_i$ in $\partial \widehat{C}_i'$. Then each of $\widehat{C}_i'$ $(i=1,2)$ is either 
$(1)$ a compression body if $D_i$ is non-separating in $\widehat{C}_i$ or $(2)$ 
a union of two compression bodies if $D_i$ is separating in $\widehat{C}_i$. 
Note that we can naturally obtain a homeomorphism 
$\partial_+ \widehat{C}_1'\to \partial_+ \widehat{C}_2'$ from the homeomorphism $\partial_+ 
\widehat{C}_1\to \partial_+ \widehat{C}_2$. Set $\widehat{M}'=\widehat{C}_1'\cup 
\widehat{C}_2'$ and $\partial_+ \widehat{C}_1'=\partial_+ \widehat{C}_2'=S'$. 
Then $(\widehat{C}_1',\widehat{C}_2';S')$ is either $(1)$ a Heegaard splitting 
or $(2)$ a union of two Heegaard splittings (\textit{cf}. Figure \ref{fig:17}). 

\begin{figure}[htb]\begin{center}
  {\unitlength=1cm
  \begin{picture}(12,4)
   \put(.3,2.6){\includegraphics[keepaspectratio]{%
   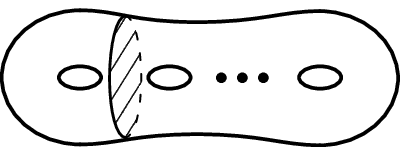}}
   \put(7.5,2.6){\includegraphics[keepaspectratio]{%
   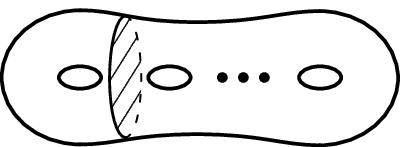}}
   \put(-.5,0){\includegraphics[keepaspectratio]{%
   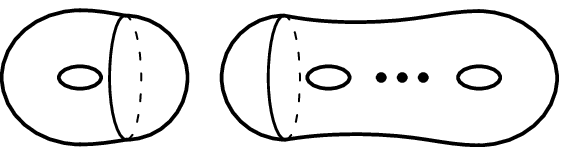}}
   \put(6.8,0){\includegraphics[keepaspectratio]{%
   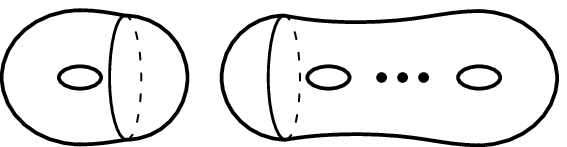}}
   \put(1.4,2.3){$D_1$}
   \put(5.75,3.4){{\Large $\cup$}}
   \put(5.8,3){$S$}
   \put(4.8,2.6){$\widehat{C}_1$}
   \put(6.64,2.6){$\widehat{C}_2$}
   \put(8.6,2.3){$D_2$}
   \put(5.8,1.8){{\Large $\Downarrow$}}
   \put(5.8,.9){{\Large $\cup$}}
   \put(5.85,.5){$S'$}
   \put(5,-.1){$\widehat{C}'_1$}
   \put(6.45,-.1){$\widehat{C}'_2$}
  \end{picture}}
  \caption{}
  \label{fig:17}
\end{center}\end{figure}

By innermost disk arguments, we see that there is a $\partial$-reducing disk of $\widehat{M}$ disjoint 
from $P'$. This implies that a component of $\widehat{M}'$ is $\partial$-reducible and 
hence one of the Heegaard splittings of $(\widehat{C}_1',\widehat{C}_2';S')$ is $\partial$-reducible. 
By induction on the genus of a Heegaard surface, we see that $(\widehat{C}_1,\widehat{C}_2;S)$ 
is $\partial$-reducible. Therefore $(C_1,C_2;S)$ is also $\partial$-reducible and 
hence we have $(2)$ of Theorem \ref{S5}. 
\end{proof}

\subsection{Waldhausen's theorem}\label{Waldhausen}

We devote this subsection to a simplified proof of the following theorem originally due to Waldhausen \cite
{Waldhausen}. To prove the theorem, we exploit Gabai's idea of  ``thin position'' 
(\textit{cf}. \cite{Gabai}), Johannson's technique (\textit{cf}. \cite{Johannson}) and Otal's idea 
(\textit{cf}. \cite{Otal}) of viewing the Heegaard splittings as a graph in the three dimensional space. 

\begin{thm}[Waldhausen]\label{S6}
Any Heegaard splitting of $S^3$ is standard, i.e., is obtained from the trivial 
Heegaard splitting by stabilization. 
\end{thm}

\bigskip\noindent
\textbf{Thin position of graphs in the \mbox{\boldmath $3$}-sphere.}
Let $\Gamma\subset S^3$ be a finite graph in which all vertices are of valence three. 
Let $h:S^3 \to [-1,1]$ be a height 
function such that $h^{-1}(t)=P(t)\cong S^2$ for $t\in (-1,1)$, $h^{-1}(-1)=($the south pole of 
$S^3)$, and $h^{-1}(1)=($the north pole of $S^3)$. Let $\mathcal{V}$ denote the set of vertices of $\Gamma$. 

\begin{defi}
\rm
The graph $\Gamma$ is in \textit{Morse position} with respect to $h$ if the following conditions are satisfied. 
\begin{enumerate}
\item $h|_{\Gamma\setminus \mathcal{V}}$ has finitely many non-degenerate critical points. 
\item The height of critical points of $h|_{\Gamma\setminus \mathcal{V}}$ and the vertices $\mathcal{V}$ 
are mutually different. 
\end{enumerate}
\end{defi}

A set of the \textit{critical heights} for $\Gamma$ is the set of height at which there is either a critical 
point of $h|_{\Gamma\setminus \mathcal{V}}$ or a component of $\mathcal{V}$. 
We can deform $\Gamma$ by an isotopy so that a regular neighborhood of each vertex $v$ of $\Gamma$ 
is either of \textit{Type-$y$} (i.e., two edges incident to $v$ is above $v$ and the remaining edge is below $v$) 
or of \textit{Type-$\lambda$} (i.e., two edges incident to $v$ is below $v$ and the remaining edge is above $v$). 
Such a graph is said to be in \textit{normal form}. We call a vertex $v$ a \textit{$y$-vertex} 
(a \textit{$\lambda$-vertex} resp.) if $\eta(v;\Gamma)$ 
is of Type-$y$ (Type-$\lambda$ resp.). 

Suppose that $\Gamma$ is in Morse position and in normal form. Note that $\eta(\Gamma;S^3)$ can be 
regarded as the union of $0$-handles corresponding to the regular neighborhood of the vertices and 
$1$-handles corresponding to the regular neighborhood of the edges. A simple loop $\alpha$ in 
$\partial \eta(\Gamma;S^3)$ is in \textit{normal form} if the following conditions are satisfied. 

\medskip\noindent
(a) For each $1$-handle $(\cong D^2\times [0,1])$, each component of $\alpha\cap (\partial D^2\times [0,1])$ 
is an essential arc in the annulus $\partial D^2\times [0,1]$. 

\noindent
(b) For each $0$-handle $(\cong B^3)$, each component of $\alpha\cap \partial B^3$ 
is an arc which is essential in the $2$-sphere with three holes 
$\mathrm{cl}(\partial B^3\setminus ($the $1$-handles incident to $B^3))$. 

\medskip
Let $D$ be a disk properly embedded in $\mathrm{cl}(S^3\setminus \eta(\Gamma;S^3))$. 
We say that $D$ is in \textit{normal form} if the following conditions are satisfied. 

\begin{figure}[htb]\begin{center}
  {\unitlength=1cm
  \begin{picture}(9.8,6.5)
   \put(0.2,2.5){\includegraphics[keepaspectratio]{%
             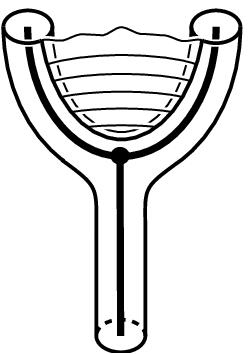}}
   \put(3.7,3.5){\includegraphics[keepaspectratio]{%
             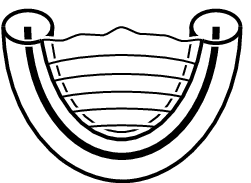}}
   \put(7.2,2){\includegraphics[keepaspectratio]{%
             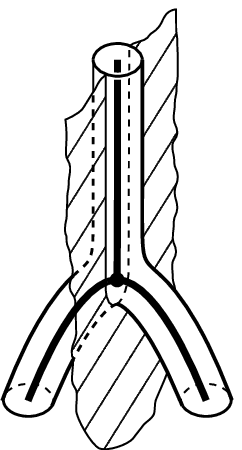}}
   \put(0,0){\includegraphics[keepaspectratio,width=2.8cm]{%
             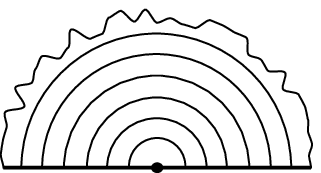}}
   \put(3.5,0){\includegraphics[keepaspectratio,width=2.8cm]{%
             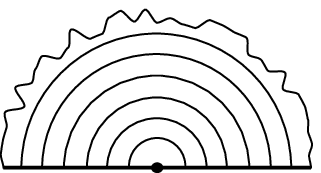}}
   \put(7,0){\includegraphics[keepaspectratio,width=2.8cm]{%
             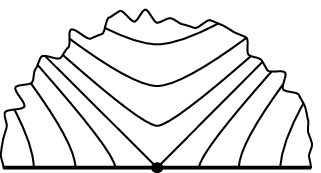}}
  \end{picture}}
  \caption{}
  \label{fig:A03}
\end{center}\end{figure}

\begin{enumerate}
\item $\partial D$ is in normal form. 
\item Each critical point of $h|_{\mathrm{int}(D)}$ is non-degenerate. 
\item No critical points of $h|_{\mathrm{int}(D)}$ occur at critical heights of $\Gamma$. 
\item No two critical points of $h|_{\mathrm{int}(D)}$ occur at the same height. 
\item $h|_{\partial D}$ is a Morse function on $\partial D$ satisfying the following 
(\textit{cf}. Figure \ref{fig:A03}). 
\begin{enumerate}
\item Each minimum of $h|_{\partial D}$ occurs either  at a $y$-vertex in ``half-center'' singularity or 
at a minimum of $\Gamma$ in ``half-center'' singularity. 
\item Each maximum of $h|_{\partial D}$ occurs either at a $\lambda$-vertex in ``half-center'' singularity or 
at a maximum of $\Gamma$ in ``half-center'' singularity. 
\end{enumerate}
\end{enumerate}

By Morse theory (\textit{cf}. \cite{Milnor}), it is known that $D$ can be put in normal form. 

\begin{figure}[htb]\begin{center}
  {\unitlength=1cm
  \begin{picture}(4.0,6.0)
   \put(.3,0){\includegraphics[keepaspectratio,height=6cm]{%
   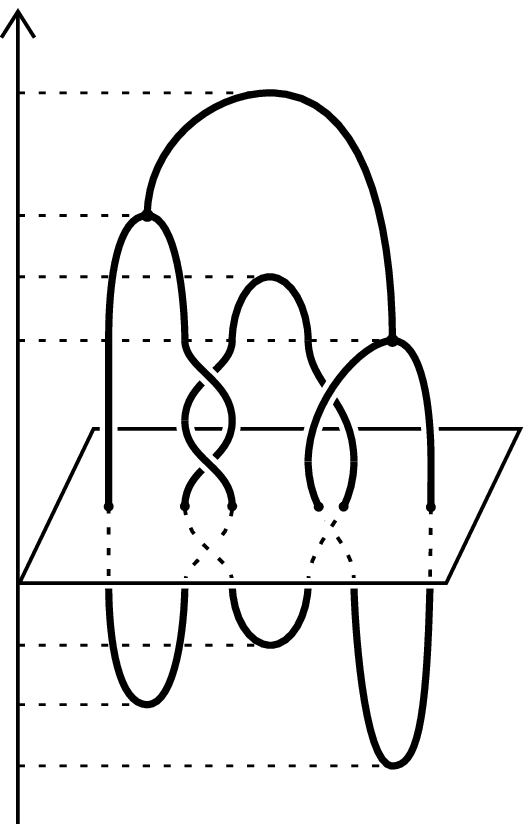}}
   \put(0,5.3){$c_6$}
   \put(0,4.4){$c_5$}
   \put(0,3.95){$c_4$}
   \put(0,3.46){$c_3$}
   \put(0,1.24){$c_2$}
   \put(0,.8){$c_1$}
   \put(0,.36){$c_0$}
   \put(3.86,2.0){$P(t)$}
  \end{picture}}
  \caption{}
  \label{fig:17.1}
\end{center}\end{figure}

Recall that $h:S^3 \to [-1,1]$ is a height 
function such that $h^{-1}(t)=P(t)\cong S^2$ for $t\in (-1,1)$, $h^{-1}(-1)=($the south pole of 
$S^3)$, and $h^{-1}(1)=($the north pole of $S^3)$. We isotope $\Gamma$ to be 
in Morse position and in normal form. 
For $t\in (-1,1)$, set $w_{\Gamma}(t)=|P(t)\cap \Gamma|$. Note that $w_{\Gamma}(t)$ is constant 
on each component of $(-1,1)\setminus ($the critical heights of $\Gamma)$. 
Set $W_{\Gamma}=\max \{w_{\Gamma}(t)|t\in (-1,1)\}$ (\textit{cf}. Figure \ref{fig:17.1}). 

Let $n_{\Gamma}$ be the number of the components of $(-1,1)\setminus ($the critical heights of $\Gamma)$ 
on which the value $W_{\Gamma}$ is attained. 

\begin{defi}
\rm
A graph $\Gamma\subset S^3$ is said to be in \textit{thin position} if $(W_{\Gamma},n_{\Gamma})$ is minimal 
with respect to lexicographic order among all graphs which are obtained from $\Gamma$ by 
ambient isotopies and edge slides and are in Morse position and in normal form. 
\end{defi}

\bigskip\noindent
\textbf{A proof of Theorem \ref{S6}.}
Let $(C_1,C_2;S)$ be a genus $g>0$ Heegaard splitting of $S^3$. 
Let $\Sigma$ be a trivalent spine of $C_1$. Note that $\eta(\partial_- C_1\cup \Sigma;M)$ is obtained from 
regular neighborhoods of $\partial_- C_1$ and the vertices of $\Sigma$ by attaching $1$-handles 
corresponding to the edges of $\Sigma$. Set $\Sigma_{\eta}=\eta(\Sigma;M)$. 
As in Section \ref{Haken}, the notation $h^0_v$, called a \textit{vertex} of $\Sigma_{\eta}$, 
means a regular neighborhood of a vertex $v$ of $\Sigma$. 
Also, the notation $h^1_{\sigma}$, called an \textit{edge} of $\Sigma_{\eta}$, 
means a $1$-handle corresponding to an edge $\sigma$ of $\Sigma$. 
Let $\Delta_1$ ($\Delta_2$ resp.) be a complete meridian system of $C_1$ 
($C_2$ resp.).  

\begin{prop}\label{B}
There is an edge of $\Sigma_{\eta}$ which is disjoint 
from $\Delta_2$, or $\Sigma$ is modified by edge slides so that the modified graph contains 
an unknotted cycle (i.e., the modified graph contains a graph 
$\alpha$ so that $\alpha$ bounds a disk in $S^3$). 
\end{prop}

\begin{proof}[Proof of Theorem \ref{S6} via Proposition \ref{B}]
We prove Theorem \ref{S6} by induction on the genus of a Heegaard surface. If $\mathit{genus}(S)=0$, 
then $(C_1,C_2;S)$ is standard (\textit{cf}. Definition \ref{trivial}). So we may assume that 
$\mathit{genus}(S)>0$ for a Heegaard splitting $(C_1,C_2;S)$. 

Suppose first that $\Sigma$ has an unknotted cycle $\alpha$. 
Then $\eta(\alpha;C_1)$ is a standard solid torus in $S^3$, that is, the exterior of $\eta(\alpha;C_1)$ 
is a solid torus. Since 
$C_1^-=\mathrm{cl}(C_1\setminus \eta(\alpha;C_1))$ is a compression body, we see that $(C_1^-,C_2;S)$ is a 
Heegaard splitting of the solid torus $\mathrm{cl}(S^3\setminus \eta(\alpha;C_1))$. Since a solid 
torus is $\partial$-reducible, $(C_1^-,C_2;S)$ is $\partial$-reducible by 
Theorem \ref{S5}, that is, there is a $\partial$-reducing disk $D_{\alpha}$ for $(C_1^-,C_2;S)$ 
with $|D_{\alpha}\cap S|=1$. 
Since $\eta(\alpha;C_1)$ is a standard solid torus in $S^3$, $D_{\alpha}$ intersects a meridian disk 
$D'_{\alpha}$ of $\eta(\alpha;C_1)$ transversely in a single point. 
Set $D_2=D_{\alpha}\cap C_2$. Then by extending $D'_{\alpha}$, we obtain a meridian disk $D_1$ of $C_1$ such that 
$\partial D_1$ intersects $\partial D_2$ transversely in a single point, i.e., $D_1$ and $D_2$ give 
stabilization of $(C_1,C_2;S)$. Hence we obtain a Heegaard 
splitting $(C_1',C_2';S')$ with $\mathit{genus}(S')<\mathit{genus}(S)$ (\textit{cf}. Figure \ref{fig:18}). 
By induction on the genus of a Heegaard surface, we can see that $(C_1,C_2;S)$ is standard. 

\begin{figure}[htb]\begin{center}
  {\unitlength=1cm
  \begin{picture}(7.5,5)
   \put(0,0){\includegraphics[height=5cm,keepaspectratio]{%
   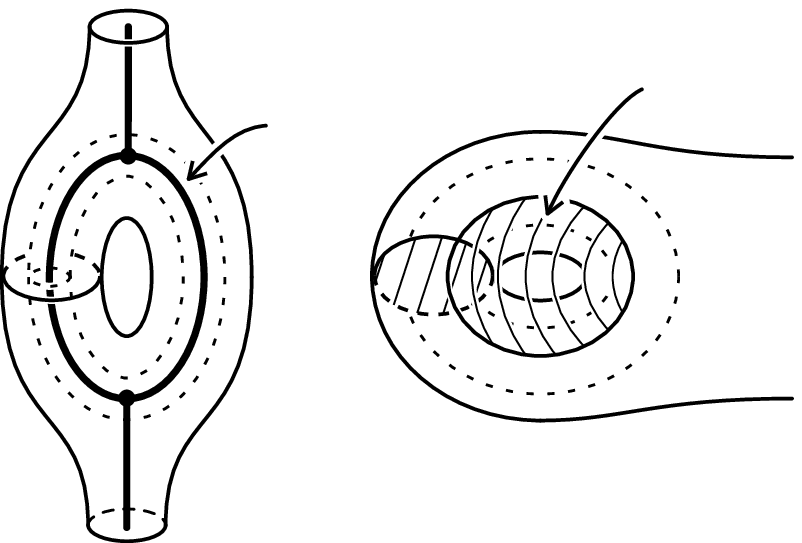}}
   \put(0,4.1){$C'_1$}
   \put(2.6,3.8){$\alpha$}
   \put(3,2.2){$D_1$}
   \put(6.1,4.2){$D_2$}
  \end{picture}}
  \caption{}
  \label{fig:18}
\end{center}\end{figure}

Suppose next that there is an edge $\sigma$ of $\Sigma$ with $h^1_{\sigma}\cap \Delta_2=\emptyset$. 
Let $D_{\sigma}$ be a meridian disk of $C_1$ which is co-core of the $1$-handle $h^1_{\sigma}$. 
Note that $D_{\sigma}\cap \Delta_2=\emptyset$. 
Cutting $C_2$ along $\Delta_2$, we obtain a union of $3$-balls and hence we see that $\partial D_{\sigma}$ 
bounds a disk, say $D'_{\sigma}$, properly embedded in one of the $3$-balls. Note that $D'_{\sigma}$ 
corresponds to a meridian disk of $C_2$. Hence we see that $(C_1,C_2;S)$ is reducible. 
It follows from a generalized Sch\"onflies theorem that every $2$-sphere in $S^3$ separates it into 
two $3$-balls (\textit{cf}. Section 2.F.5 of \cite{Rolfsen}). 
Hence by cutting $S^3$ along the reducing $2$-sphere and capping off $3$-balls, we obtain 
two Heegaard splittings of $S^3$ such that the genus of each Heegaard surface is less than that of $S$. 
Then we see that $(C_1,C_2;S)$ is standard by induction on the genus of a Heegaard surface. 
\end{proof}

In the remainder, we prove Proposition \ref{B}. Let $h:S^3 \to [-1,1]$ be a height 
function such that $h^{-1}(t)=P(t)\cong S^2$ for $t\in (-1,1)$, $h^{-1}(-1)=($the south pole of 
$S^3)$, and $h^{-1}(1)=($the north pole of $S^3)$. We may assume that $\Sigma$ is in thin position. 
We also assume that each component of $\Delta_2$ is in normal form, $\Delta_1$ intersects $\Delta_2$ 
transversely and $|\Delta_1\cap \Delta_2|$ is minimal. 

For the proof of Proposition \ref{B}, it is enough to show the following; if there are no edges of 
$\Sigma_{\eta}$ which are disjoint from $\Delta_2$, then $\Sigma$ is modified by edge slides so that 
the modified graph contains an unknotted cycle. 
Hence we suppose that there are no edges of $\Sigma_{\eta}$ which are disjoint from $\Delta_2$. 
Set $\Lambda(t)=P(t)\cap (\Sigma_{\eta}\cup \Delta_2)$. We note that $P(t)$, $\Sigma$ and $\Delta_2$ intersect 
transversely at a regular height $t$. In the following, we mainly consider such a regular height 
$t$ with $\Lambda(t)\ne \emptyset$ unless otherwise denoted. 
We also note that we may assume that $\Lambda(t)$ does not contain a loop component 
by an argument similar to the proof of Lemma \ref{A1}. Hence 
$\Lambda(t)$ is regarded as a graph in $P(t)$ which consists of 
\textit{fat-vertices} $P(t)\cap \Sigma_{\eta}$ and edges $P(t)\cap \Delta_2$. 

\begin{lem}\label{B1}
If there is a fat-vertex of $\Lambda(t)$ with valence less than two, then 
$\Sigma$ is modified by edge slides so that the modified graph contains 
an unknotted cycle. 
\end{lem}

\begin{proof}
Suppose that there is a fat-vertex $D_v$ of $\Lambda(t)$ with valence less than two. 
Let $\sigma$ be the edge of $\Sigma$ with $h^1_{\sigma}\supset D_v$ and $p$ one of the endpoints of $\sigma$. 
Since we assume that there are no edges of $\Sigma_{\eta}$ which are disjoint from $\Delta_2$, 
we see that any fat-vertex of $\Lambda(t)$ is of valence greater than zero. 
Hence $D_v$ is of valence one. Then there is the disk component $D$ of $\Delta_2$ with 
$h^1_{\sigma}\cap D \ne \emptyset$. Since $\partial D$ intersects the fat-vertex $D_v$ in a 
single point and hence $\partial D$ intersects $h^1_{\sigma}$ in a single arc, 
we can perform an edge slide on $\sigma$ along 
$\mathrm{cl}(\partial D\setminus h^1_{\sigma})$ to obtain a new graph $\Sigma'$ from $\Sigma$ 
(\textit{cf}. Figure \ref{fig:19}). 
Clearly, $\Sigma'$ contains an unknotted cycle (bounding a disk corresponding to $D_2$). 

\begin{figure}[htb]\begin{center}
  {\unitlength=1cm
  \begin{picture}(5,5)
   \put(0,0){\includegraphics[keepaspectratio]{%
   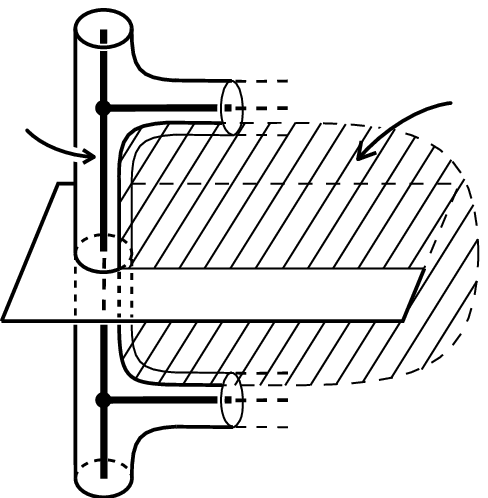}}
   \put(.05,3.83){$\sigma$}
   \put(-.8,1.8){$P(t)$}
   \put(4.6,4){$D$}
  \end{picture}}
  \caption{}
  \label{fig:19}
\end{center}\end{figure}

\end{proof}

An edge of a graph $\Lambda(t)$ is said to be \textit{simple} if the edge joins 
distinct two fat-vertices of $\Lambda(t)$. 
Recall that an edge of a graph $\Lambda(t)$ is called a \textit{loop} if the edge is not simple. 

\begin{lem}\label{B2}
Suppose that there are no fat-vertices of valence less than two. 
Then there exists a fat-vertex $D_w$ of $\Lambda(t)$ such that any 
outermost edge for $D_w$ of $\Lambda(t)$ is simple. 
\end{lem}

\begin{proof}
If $\Lambda(t)$ does 
not contain a loop, then we are done. So we may assume that $\Lambda(t)$ contains a 
loop, say $\mu$. Let $D_v$ be a fat-vertex of $\Lambda(t)$ which is a bese of $\mu$. 
By an argument similar to the proof of Lemma \ref{A2}, we can see 
that $\mu$ cuts $\mathrm{cl}(P(t)\setminus D_v)$ into two disks, and each of the two disks contains a 
fat-vertex of $\Lambda(t)$. Let $\mu_0$ be a loop of $\Lambda(t)$ which is innermost in $P(t)$. 
Let $D_w$ be a fat-vertex contained in the interior of the innermost disk bounded by $\mu_0$. 
Note that $D_w$ is not isolated and that every edge contained in the interior of the innermost disk 
is simple. Hence any outermost edge for $D_w$ of $\Lambda(t)$ is simple. 
\end{proof}

Let $D_w$ be a fat-vertex of $\Lambda(t)$ with a simple edge $\gamma(\subset \Lambda(t))$. 
We may assume that $\gamma$ is a simple outermost edge for $D_w$ of $\Lambda(t)$ and 
$\gamma$ is contained in a disk component $D$ of $\Delta_2$. It follows from Lemma \ref{B2} that we can always 
find such a fat-vertex $D_w$ and an edge $\gamma$ if each fat-vertex of $\Lambda(t)$ 
is of valence greater than one. Let $\delta_{\gamma}$ be the outermost disk for $(D_w,\gamma)$. 
We say an outermost edge $\gamma$ is \textit{upper} (\textit{lower} resp.) if 
$\eta(\gamma;\delta_{\gamma})$ is above (below resp.) $\gamma$ with respect to 
the height function $h$. Let $t_0$ be a regular height with $w_{\Sigma}(t_0)=W_{\Sigma}$. 

\begin{lem}\label{B3}
Let $D_w$ be a fat-vertex of $\Lambda(t_0)$ with a simple outermost edge for $D_w$ of $\Lambda(t)$. 
Then we have one of the following. 
\begin{enumerate}
\item All the simple outermost edges for $D_w$ of $\Lambda(t)$ are either upper or lower. 
\item $\Sigma$ is modified by edge slides so that the modified graph contains an unknotted cycle. 
\end{enumerate}
\end{lem}

\begin{proof}
Suppose that $\Lambda(t)$ contains simple outermost edges for $D_w$, say $\gamma$ and $\gamma'$, 
such that $\gamma$ is upper and $\gamma'$ is lower. For the proof of Lemma \ref{B3}, 
it is enough to show that $\Sigma$ is 
modified by edge slides so that there is an unknotted cycle. 
Let $\delta_{\gamma}$ and $\delta_{\gamma'}$ be the outermost disk for $(D_w,\gamma)$ and 
$(D_{w'},\gamma')$ respectively. 
Let $\sigma$ be the edge of $\Sigma$ with $h^1_{\sigma} \supset D_w$. 
Let $\bar{\gamma}$ ($\bar{\gamma}'$ resp.) be a union of the components obtained by cutting $\sigma$ 
by the two fat-vertices of $\Lambda(t_0)$ incident to $\gamma$ ($\gamma'$ resp.) 
such that a $1$-handle correponding to each component intersects 
$\partial \delta_{\gamma}\setminus \gamma$ ($\partial \delta_{\gamma'}\setminus \gamma'$ resp.). 
We note that $\bar{\gamma}$ ($\bar{\gamma}'$ resp.) satisfies one of the following conditions. 

\begin{enumerate}
\item $\bar{\gamma}$ ($\bar{\gamma}'$ resp.) consists of an arc such that 
$\bar{\gamma}$ ($\bar{\gamma}'$ resp.) and $\sigma$ share a single endpoint. 
\item $\bar{\gamma}$ ($\bar{\gamma}'$ resp.) consists of an arc with $\bar{\gamma}\subset \mathrm{int}(\sigma)$ 
($\bar{\gamma}'\subset \mathrm{int}(\sigma)$ resp.). 
\item $\bar{\gamma}$ ($\bar{\gamma}'$ resp.) consists of two subarcs of $\sigma$ such that each component of 
$\bar{\gamma}$ ($\bar{\gamma}'$ resp.) and 
$\sigma$ share  a single endpoint. 
\end{enumerate}

In each of the conditions above, corresponding figures are illustrated in Figure \ref{fig:21.1}. 
(We remark that there is a case in the condition $(3)$ which is similar to the latter half of 
Case $3$ in the proof of Lemma \ref{A4}. Recall that we only have to use broken edge slides in that case. 
Hence we omit details in such a case here.)

\begin{figure}[htb]\begin{center}
  {\unitlength=1cm
  \begin{picture}(9.0,12.0)(.5,0)
   \put(1.0,8.6){\includegraphics[keepaspectratio]{
             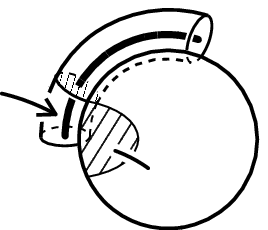}}
   \put(5.0,8){\includegraphics[keepaspectratio,width=3.8cm]{%
             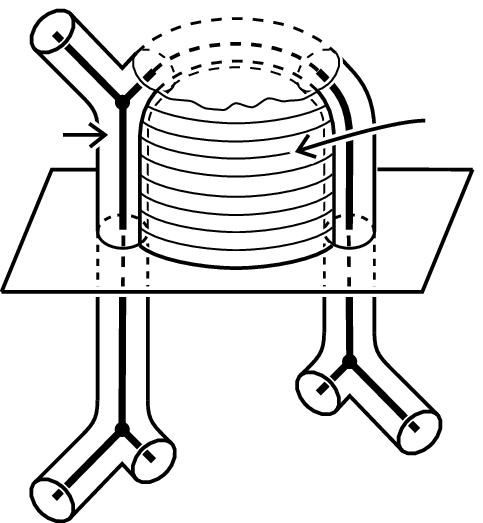}}
   \put(1.0,4.6){\includegraphics[keepaspectratio]{
             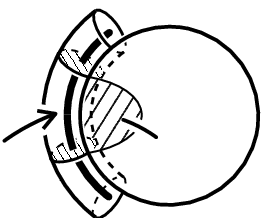}}
   \put(5.0,4){\includegraphics[keepaspectratio,width=3.8cm]{%
             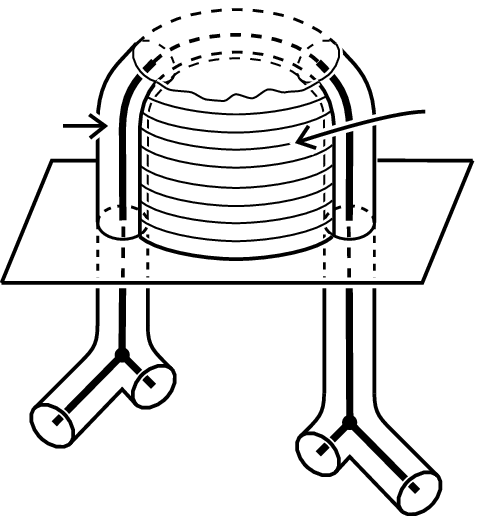}}
   \put(1.0,.6){\includegraphics[keepaspectratio]{
             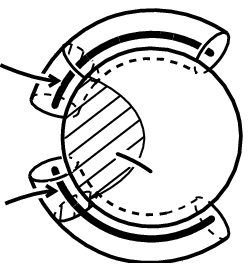}}
   \put(5.0,0){\includegraphics[keepaspectratio,width=3.8cm]{%
             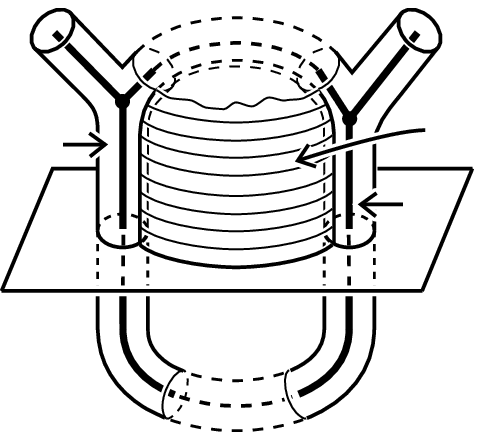}}
   \put(0,9.4){{\large\bfseries I}}
   \put(0,5.4){{\large\bfseries II}}
   \put(0,1.4){{\large\bfseries III}}
   \put(2.55,9.0){$\delta_\gamma$}
   \put(2.6,5.2){$\delta_\gamma$}
   \put(2.6,1.4){$\delta_\gamma$}
   \put(.7,10){$\bar{\gamma}$}
   \put(.8,5.2){$\bar{\gamma}$}
   \put(.7,1.0){$\bar{\gamma}$}
   \put(.7,2.1){$\bar{\gamma}$}
   \put(8.5,11.1){$\delta_\gamma$}
   \put(8.5,7.1){$\delta_\gamma$}
   \put(8.5,2.3){$\delta_\gamma$}
   \put(5.2,11){$\bar{\gamma}$}
   \put(5.2,7){$\bar{\gamma}$}
   \put(5.2,2.2){$\bar{\gamma}$}
   \put(8.3,1.7){$\bar{\gamma}$}
   \put(5.2,10.1){$D_w$}
   \put(5.2,6.1){$D_w$}
   \put(5.2,1.4){$D_w$}
  \end{picture}}
  \caption{}
  \label{fig:21.1}
\end{center}\end{figure}

\medskip
\noindent\textit{Case} (1)-(1).\ \ \ \ 
Both $\bar{\gamma}$ and $\bar{\gamma}'$ satisfy the condition $(1)$. 

\medskip
If the endpoints of $\gamma$ are the same as those of $\gamma'$, then we can slide 
$\bar{\gamma}$ ($\bar{\gamma}'$ resp.) to $\gamma$ ($\gamma'$ resp.) along the disk $\delta_{\gamma}$ 
($\delta_{\gamma'}$ resp.) and hence we obtain an unknotted cycle (\textit{cf}. Figure \ref{fig:23}). 
Otherwise, we can perform a Whitehead move on $\Sigma$ to reduce $(W_{\Sigma},n_{\Sigma})$, a contradiction 
(\textit{cf}. Figure \ref{fig:24}). 

\begin{figure}[htb]\begin{center}
  {\unitlength=1cm
  \begin{picture}(9.8,4.5)
   \put(0,0){\includegraphics[keepaspectratio,width=4.0cm]{%
   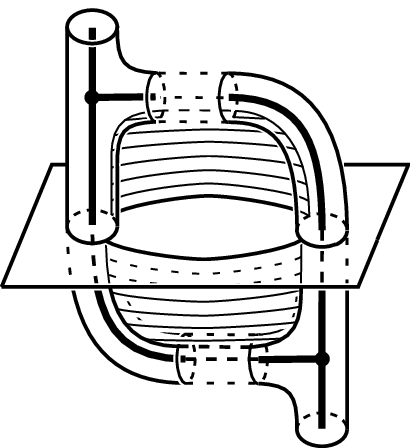}}
   \put(5.6,0){\includegraphics[keepaspectratio,width=4.0cm]{%
   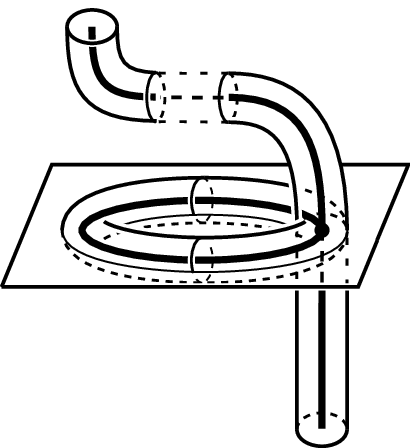}}
   \put(-.9,1.7){$P(t_0)$}
   \put(4.57,2.1){$\Longrightarrow$}
  \end{picture}}
  \caption{}
  \label{fig:23}
\end{center}\end{figure}

\begin{figure}[htb]\begin{center}
  {\unitlength=1cm
  \begin{picture}(12,2.6)
   \put(0,0){\includegraphics[keepaspectratio,width=3.0cm]{
             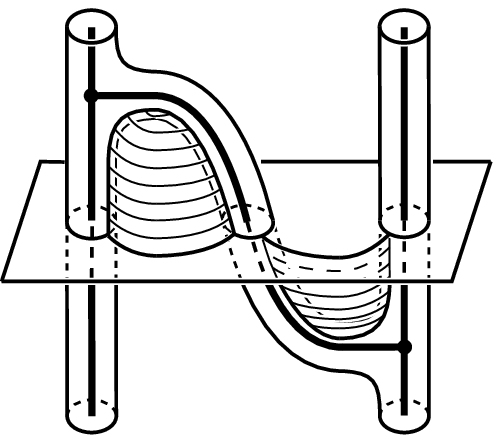}}
   \put(4.5,0){\includegraphics[keepaspectratio,width=3.0cm]{
               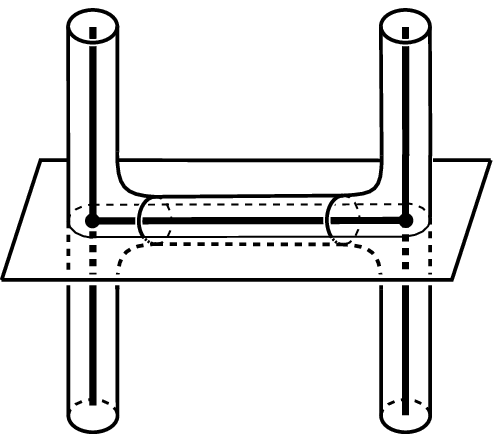}}
   \put(9,0){\includegraphics[keepaspectratio,width=3.0cm]{
             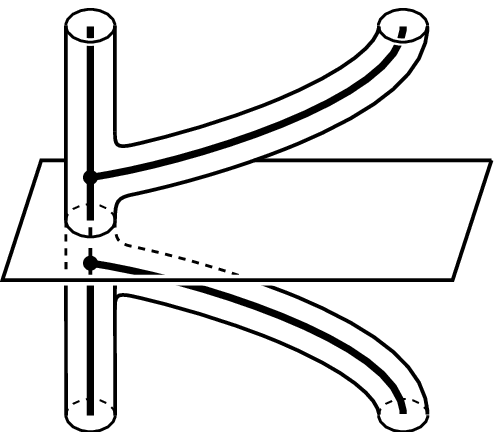}}   
   \put(3.5,1.2){$\Longrightarrow$}
   \put(8.0,1.2){$\Longrightarrow$}
  \end{picture}}
  \caption{}
  \label{fig:24}
\end{center}\end{figure}

\medskip
\noindent\textit{Case} (1)-(2).\ \ \ \ 
Either $\bar{\gamma}$ or $\bar{\gamma}'$, say $\bar{\gamma}$, satisfies the condition $(1)$ and 
$\bar{\gamma}'$ satisfies the condition $(2)$. 

\medskip
Then we can slide $\bar{\gamma}$ ($\bar{\gamma}'$ resp.) to $\gamma$ ($\gamma'$ resp.) along the disk 
$\delta_{\gamma}$ ($\delta_{\gamma'}$ resp.). Then $\Sigma$ is further isotoped to 
reduce $(W_{\Sigma},n_{\Sigma})$, a contradiction (\textit{cf}. Figure \ref{fig:23.1}). 

\begin{figure}[htb]\begin{center}
  {\unitlength=1cm
  \begin{picture}(12.0,4.5)
   \put(0,1.0){\includegraphics[keepaspectratio,width=5.5cm]{%
               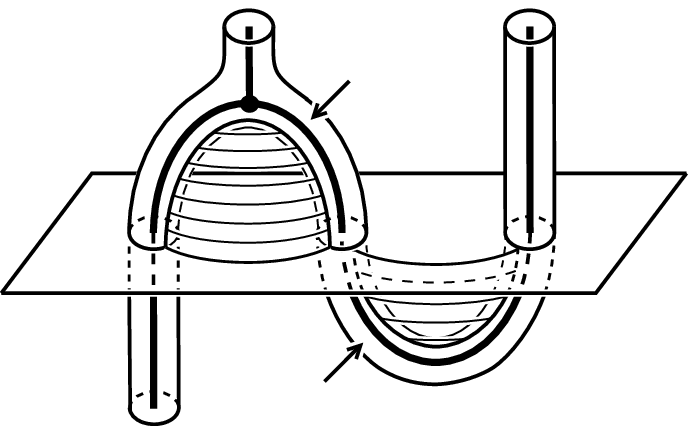}}
   \put(6.6,1.0){\includegraphics[keepaspectratio,width=5.5cm]{%
               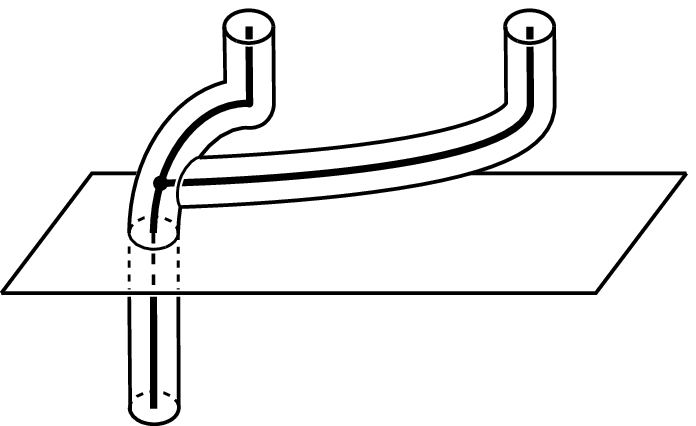}}
   \put(4.2,0){\includegraphics[keepaspectratio,height=.6cm]{%
               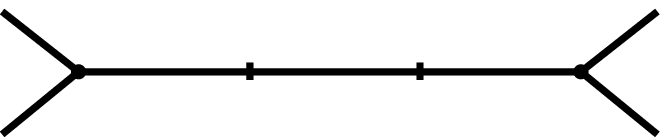}}
   \put(5.8,2.4){$\Longrightarrow$}
   \put(2.9,3.8){$\bar{\gamma}$}
   \put(2.2,1.0){$\bar{\gamma}'$}
   \put(3.0,2.6){$D_w$}
   \put(4.85,0.85){$\bar{\gamma}$}
   \put(4.6,0.5){$\overbrace{\hspace{0.5cm}}$}
   \put(5.65,0.85){$\bar{\gamma}'$}
   \put(5.4,0.5){$\overbrace{\hspace{0.5cm}}$}
   \put(5.2,-0.1){$w$}
  \end{picture}}
  \caption{}
  \label{fig:23.1}
\end{center}\end{figure}

\medskip
\noindent\textit{Case} (1)-(3).\ \ \ \ 
Either $\bar{\gamma}$ or $\bar{\gamma}'$, say $\bar{\gamma}$, satisfies the condition $(1)$ and 
$\bar{\gamma}'$ satisfies the condition $(3)$. 

\medskip
Let $\bar{\gamma}'_1$ and $\bar{\gamma}'_2$ be the components of $\bar{\gamma}'$ with 
$h^1_{\bar{\gamma}'_1}\supset D_w$. Note that $\bar{\gamma}\supset \bar{\gamma}'_2$ and hence 
$\mathrm{int}(\bar{\gamma})\supset \partial \bar{\gamma}'_2$. This implies that 
$\mathrm{int}(\bar{\gamma})\cap P(t_0)\ne \emptyset$. Hence we can slide $\bar{\gamma}$ to $\gamma$ 
along the disk $\delta_{\gamma}$. We can further isotope $\Sigma$ slightly to reduce 
$(W_{\Sigma},n_{\Sigma})$, a contradiction 
(\textit{cf}. Figure \ref{fig:23.2}). 

\begin{figure}[htb]\begin{center}
  {\unitlength=1cm
  \begin{picture}(12.0,5.4)
   \put(0,1.4){\includegraphics[keepaspectratio,width=5.5cm]{%
               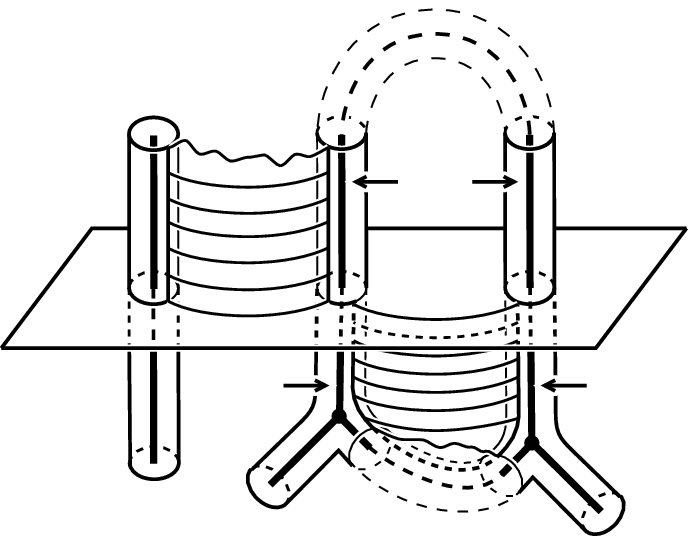}}
   \put(6.6,1.1){\includegraphics[keepaspectratio,width=5.5cm]{%
               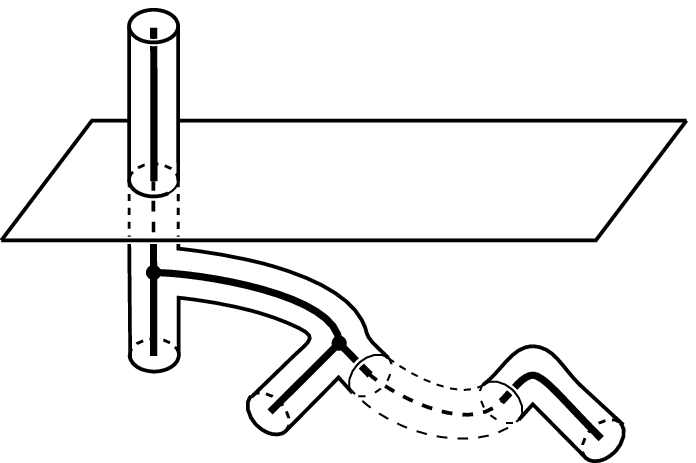}}
   \put(4.5,0.4){\includegraphics[keepaspectratio,height=.6cm]{%
               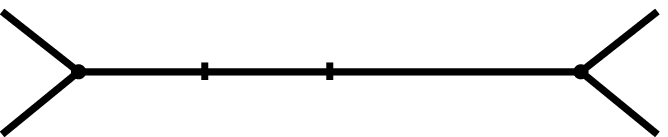}}
   \put(5.8,3.3){$\Longrightarrow$}
   \put(3.34,4.1){$\bar{\gamma}$}
   \put(1.9,2.5){$\bar{\gamma}'_1$}
   \put(4.8,2.5){$\bar{\gamma}'_2$}
   \put(3.0,3.3){$D_w$}
   \put(5.3,1.14){$\bar{\gamma}$}
   \put(4.9,0.8){$\overbrace{\hspace{1cm}}$}
   \put(4.84,0.6){$\underbrace{\hspace{0cm}}$}
   \put(6.1,0.6){$\underbrace{\hspace{1cm}}$}
   \put(6,0.85){$w$}
   \put(5.2,0.1){$\bar{\gamma}'_2$}
   \put(6.2,0.1){$\bar{\gamma}'_1$}
  \end{picture}}
  \caption{}
  \label{fig:23.2}
\end{center}\end{figure}

\medskip
\noindent\textit{Case} (2)-(2).\ \ \ \ 
Both $\bar{\gamma}$ and $\bar{\gamma}'$ satisfy the condition $(2)$. 

\medskip
Then we can slide $\bar{\gamma}$ ($\bar{\gamma}'$ resp.) to $\gamma$ ($\gamma'$ resp.) along the disk 
$\delta_{\gamma}$ ($\delta_{\gamma'}$ resp.). Moreover, we can isotope $\sigma$ slightly to reduce 
$(W_{\Sigma},n_{\Sigma})$, a contradiction (\textit{cf}. Figure \ref{fig:22}). 

\begin{figure}[htb]\begin{center}
  {\unitlength=1cm
  \begin{picture}(12.0,4.6)
   \put(0,1.4){\includegraphics[keepaspectratio]{
   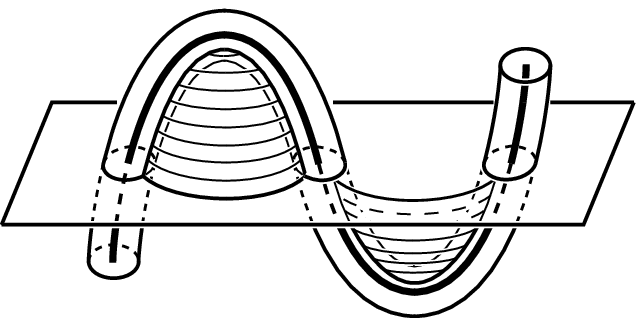}}
   \put(7.6,1.78){\includegraphics[keepaspectratio]{%
   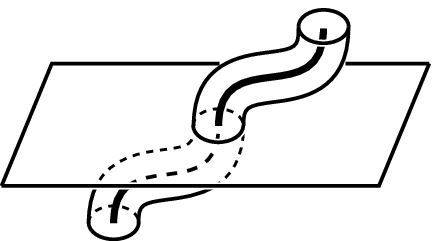}}
   \put(4.6,0.4){\includegraphics[keepaspectratio,height=.6cm]{%
   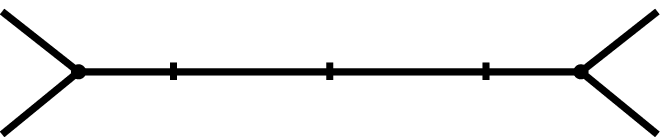}}
   \put(6.8,2.9){$\Longrightarrow$}
   \put(-.9,2.4){$P(t_0)$}
   \put(3.6,3.0){$D_w$}
   \put(1.6,2.5){$\gamma$}
   \put(4.4,2.8){$\gamma'$}
   \put(6.0,0.85){$w$}
   \put(5.35,0.6){$\underbrace{\hspace{0cm}}$}
   \put(5.7,0){$\bar{\gamma}$}
   \put(6.2,0.6){$\underbrace{\hspace{0cm}}$}
   \put(6.6,0){$\bar{\gamma}'$}
  \end{picture}}
  \caption{}
  \label{fig:22}
\end{center}\end{figure}

\medskip
\noindent\textit{Case} (2)-(3).\ \ \ \ 
Either $\bar{\gamma}$ or $\bar{\gamma}'$, say $\bar{\gamma}$, satisfies the condition $(2)$ and 
$\bar{\gamma}'$ satisfies the condition $(3)$. 

\medskip
Note that $\bar{\gamma}'$ consists of two arcs, say $\bar{\gamma}'_1$ and $\bar{\gamma}'_2$, with 
$\bar{\gamma}_1'\cap \bar{\gamma}=D_w$. Then we have the following cases. 

(i) $\bar{\gamma}_2'$ is disjoint from $\bar{\gamma}$. 
In this case, we can slide $\bar{\gamma}$ ($\bar{\gamma}'$ resp.) to $\gamma$ ($\gamma'$ resp.) along 
the disk $\delta_{\gamma}$ ($\delta_{\gamma'}$ resp.). Moreover, we can isotope $\sigma$ slightly to reduce 
$(W_{\Sigma},n_{\Sigma})$, a contradiction (\textit{cf}. Figure \ref{fig:23.3}). 

\begin{figure}[htb]\begin{center}
  {\unitlength=1cm
  \begin{picture}(12.0,9.0)
   \put(0,5){\includegraphics[keepaspectratio,width=5.5cm]{%
               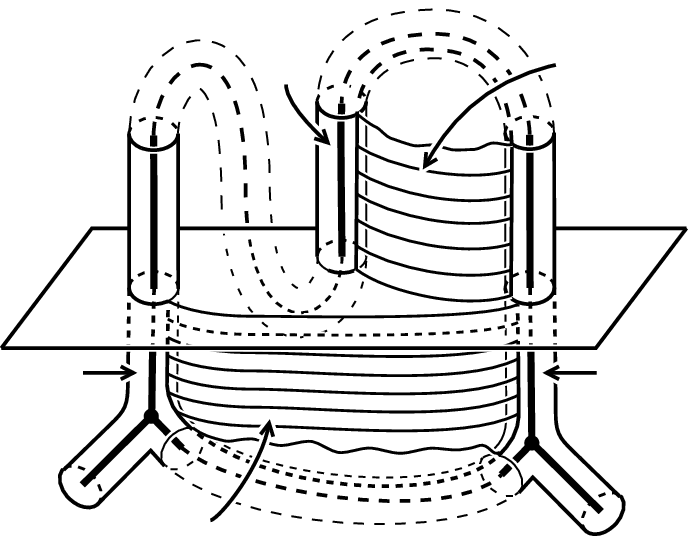}}
   \put(6.6,5){\includegraphics[keepaspectratio,width=5.5cm]{%
               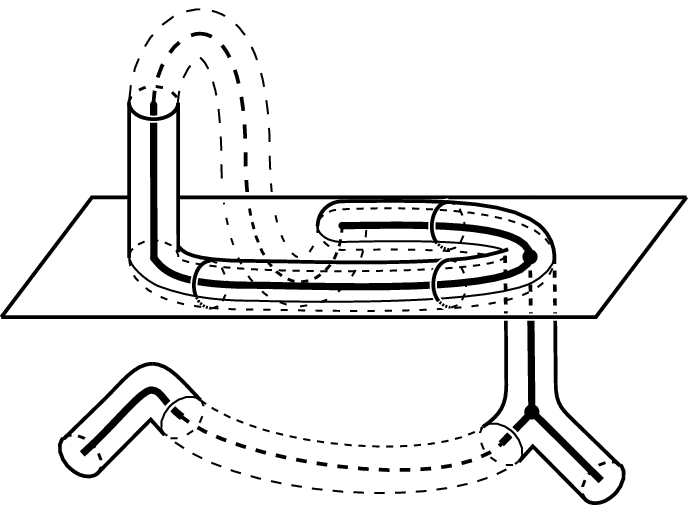}}
   \put(6.6,0.0){\includegraphics[keepaspectratio,width=5.5cm]{%
               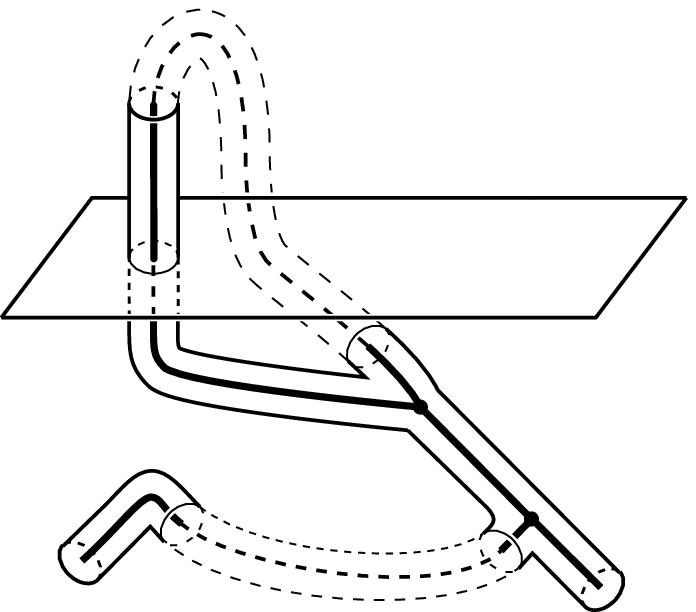}}
   \put(1.5,2.0){\includegraphics[keepaspectratio,height=.6cm]{%
               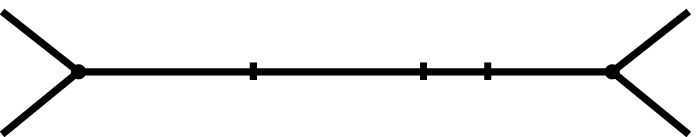}}
   \put(5.8,6.8){$\Longrightarrow$}
   \put(9.3,4.4){$\Downarrow$}
   \put(4.5,8.7){$\delta_\gamma$}
   \put(1.4,4.9){$\delta_{\gamma}'$}
   \put(2.1,8.7){$\bar{\gamma}$}
   \put(.3,6.1){$\bar{\gamma}'_2$}
   \put(4.8,6.1){$\bar{\gamma}'_1$}
   \put(4.43,6.9){$D_w$}
   \put(2.4,2.5){$w$}
   \put(2.7,2.4){$\overbrace{\hspace{0.8cm}}$}
   \put(3.0,2.74){$\bar{\gamma}$}
   \put(1.9,2.2){$\underbrace{\hspace{0.8cm}}$}
   \put(3.8,2.2){$\underbrace{\hspace{0.0cm}}$}
   \put(2.1,1.6){$\bar{\gamma}'_1$}
   \put(3.9,1.6){$\bar{\gamma}'_2$}
  \end{picture}}
  \caption{}
  \label{fig:23.3}
\end{center}\end{figure}

(ii) $\bar{\gamma}_2'\cap \bar{\gamma}$ consists of a point, i.e., $\bar{\gamma}_2'$ and $\bar{\gamma}$ share 
one endpoint. Note that $\bar{\gamma}\cup \bar{\gamma}'=\sigma$ and 
$\bar{\gamma}\cap \bar{\gamma}'=\partial \bar{\gamma}$. In this case, we can slide 
$\bar{\gamma}$ ($\bar{\gamma}'$ resp.) to $\gamma$ ($\gamma'$ resp.) along the disk $\delta_{\gamma}$ 
($\delta_{\gamma'}$ resp.) and hence we obtain an unknotted cycle (\textit{cf}. Figure \ref{fig:23.4}). 

\begin{figure}[htb]\begin{center}
  {\unitlength=1cm
  \begin{picture}(12.0,5.4)
   \put(0,1.4){\includegraphics[keepaspectratio,width=5.5cm]{%
               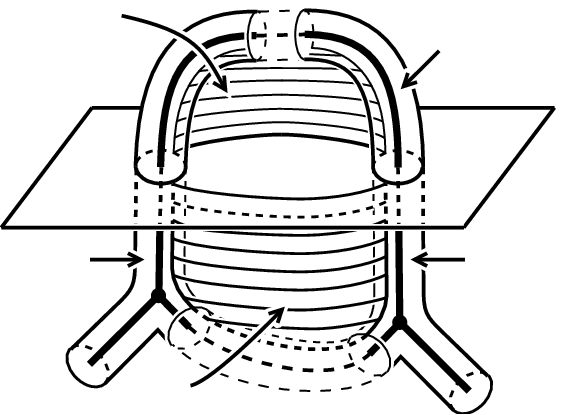}}
   \put(6.6,1.4){\includegraphics[keepaspectratio,width=5.5cm]{%
               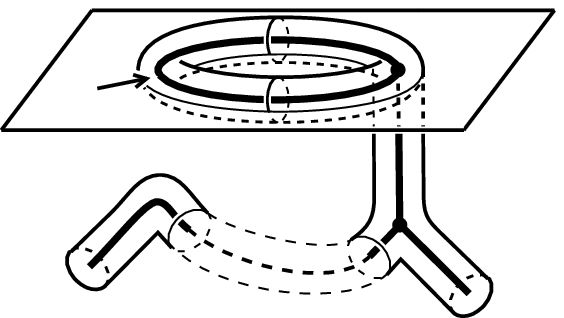}}
   \put(4.2,0.4){\includegraphics[keepaspectratio,height=.6cm]{%
               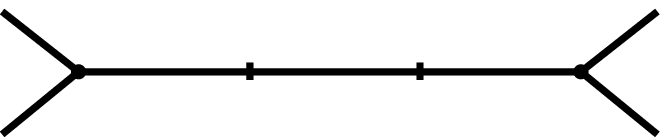}}
   \put(5.7,3.7){$\Longrightarrow$}
   \put(.8,5.3){$\delta_\gamma$}
   \put(1.5,1.5){$\delta_{\gamma'}$}
   \put(4.4,4.9){$\bar{\gamma}$}
   \put(.5,2.8){$\bar{\gamma}'_2$}
   \put(4.6,2.8){$\bar{\gamma}'_1$}
   \put(4.3,3.8){$D_w$}
   \put(7.3,3.5){$\sigma$}
   \put(5.1,0.85){$w$}
   \put(5.65,1.2){$\bar{\gamma}$}
   \put(5.4,0.85){$\overbrace{\hspace{.7cm}}$}
   \put(4.6,0.6){$\underbrace{\hspace{.7cm}}$}
   \put(6.2,0.6){$\underbrace{\hspace{.7cm}}$}
   \put(4.8,0.0){$\bar{\gamma}'_1$}
   \put(6.4,0.0){$\bar{\gamma}'_2$}
  \end{picture}}
  \caption{}
  \label{fig:23.4}
\end{center}\end{figure}

(iii) $\bar{\gamma}_2'\cap \bar{\gamma}$ consists of an arc. In this case, 
we can slide $\bar{\gamma}$ to $\gamma$ along the disk 
$\delta_{\gamma}$. Since $\bar{\gamma}_2'\cap \bar{\gamma}$ consists of an arc, $\bar{\gamma}$ contains 
at least three critical points. Hence we can further isotope $\sigma$ slightly to reduce 
$(W_{\Sigma},n_{\Sigma})$, a contradiction (\textit{cf}. Figure \ref{fig:23.5}). 

\begin{figure}[htb]\begin{center}
  {\unitlength=1cm
  \begin{picture}(12.0,11.5)
   \put(0,5){\includegraphics[keepaspectratio,width=5.5cm]{%
               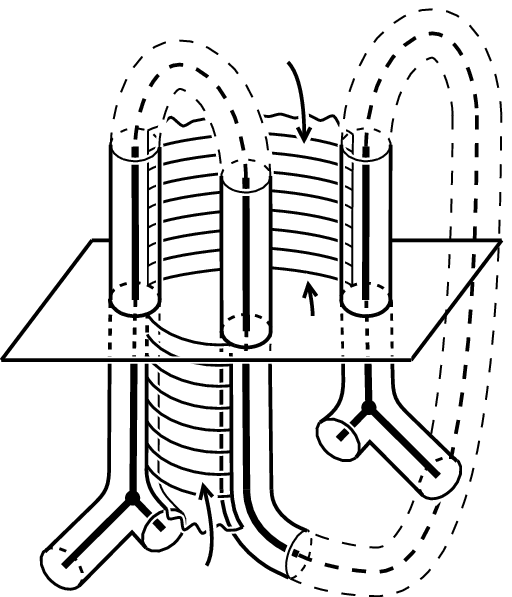}}
   \put(6.6,5){\includegraphics[keepaspectratio,width=5.5cm]{%
               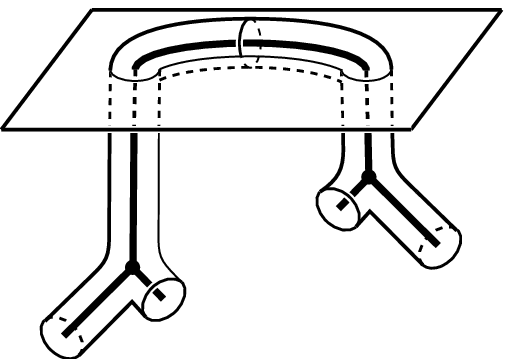}}
   \put(6.6,0.0){\includegraphics[keepaspectratio,width=5.5cm]{%
               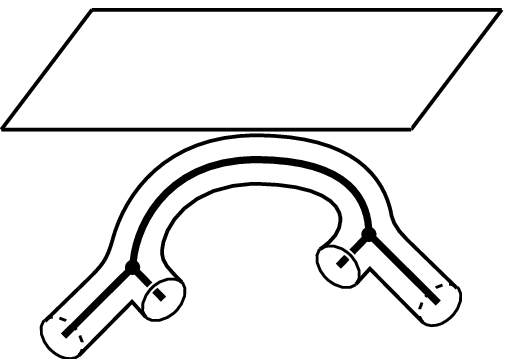}}
   \put(1.5,2.0){\includegraphics[keepaspectratio,height=.6cm]{%
               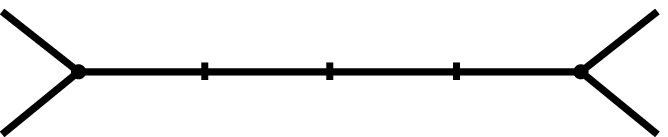}}
   \put(5.8,8.0){$\Longrightarrow$}
   \put(9.3,4.4){$\Downarrow$}
   \put(2.7,10.8){$\delta_\gamma$}
   \put(1.7,5.1){$\delta_{\gamma'}$}
   \put(2.0,7.9){$\gamma'$}
   \put(3.3,7.8){$\gamma$}
   \put(.64,8.1){$D_w$}
   \put(2.2,2.5){$w$}
   \put(2.46,2.4){$\overbrace{\hspace*{1.2cm}}$}
   \put(2.9,2.74){$\bar{\gamma}$}
   \put(1.8,2.2){$\underbrace{\hspace*{0cm}}$}
   \put(3.1,2.2){$\underbrace{\hspace*{1.1cm}}$}
   \put(2.0,1.6){$\bar{\gamma}'_1$}
   \put(3.4,1.6){$\bar{\gamma}'_2$}
  \end{picture}}
  \caption{}
  \label{fig:23.5}
\end{center}\end{figure}

\medskip
\noindent\textit{Case} (3)-(3).\ \ \ \ 
Both $\bar{\gamma}$ and $\bar{\gamma}'$ satisfy the condition $(3)$. 

\medskip
Let $\bar{\gamma}_1$ and $\bar{\gamma}_2$ ($\bar{\gamma}'_1$ and $\bar{\gamma}'_2$ resp.) 
be the components of $\bar{\gamma}$ ($\bar{\gamma}'$ resp.) with $h^1_{\bar{\gamma}_1}\supset D_w$ 
($h^1_{\bar{\gamma}'_1}\supset D_w$ resp.). 
Note that $\bar{\gamma}_1\supset \bar{\gamma}'_2$ and $\bar{\gamma}'_1\supset \bar{\gamma}_2$. 
In this case, we can slide $\bar{\gamma}'_1$ to $\gamma'$ along the disk 
$\delta_{\gamma'}$. Since $\bar{\gamma}'_1\supset \bar{\gamma}_2$, $\bar{\gamma}'_1$ contains 
at least one critical point. Hence we can further isotope $\sigma$ slightly to reduce $(W_{\Sigma},n_{\Sigma})$, 
a contradiction (\textit{cf}. Figure \ref{fig:23.6}). 

\begin{figure}[htb]\begin{center}
  {\unitlength=1cm
  \begin{picture}(12.0,11.5)
   \put(0,5){\includegraphics[keepaspectratio,width=5.5cm]{%
               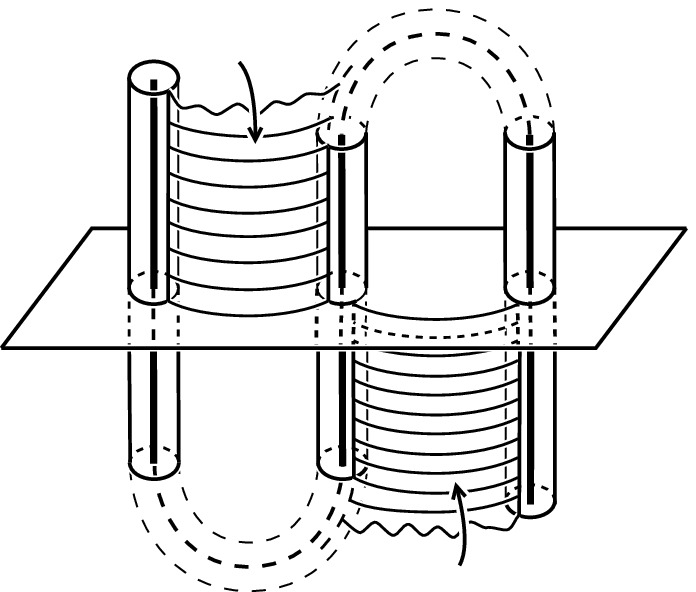}}
   \put(6.6,5){\includegraphics[keepaspectratio,width=5.5cm]{%
               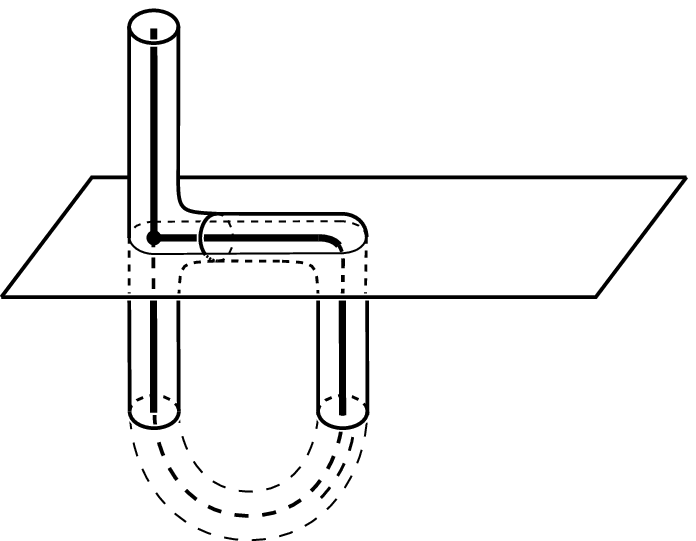}}
   \put(6.6,0.0){\includegraphics[keepaspectratio,width=5.5cm]{%
               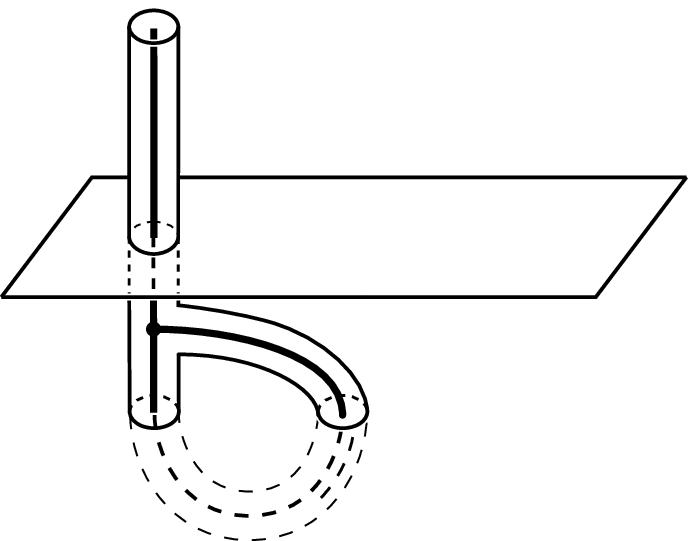}}
   \put(1.5,2.3){\includegraphics[keepaspectratio,height=.6cm]{%
               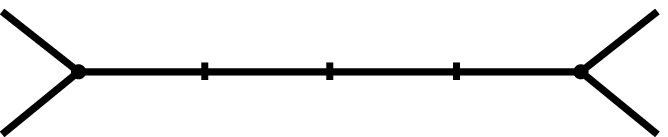}}
   \put(5.8,7.3){$\Longrightarrow$}
   \put(9.3,4.4){$\Downarrow$}
   \put(1.5,9.3){$\delta_\gamma$}
   \put(3.1,5.0){$\delta_{\gamma'}$}
   \put(3.0,7.4){$D_w$}
   \put(2.9,2.25){$w$}
   \put(2.3,3.2){$\bar{\gamma}_1$}
   \put(3.7,3.2){$\bar{\gamma}_2$}
   \put(1.9,2.8){$\overbrace{\hspace*{1.1cm}}$}
   \put(3.6,2.8){$\overbrace{\hspace*{0cm}}$}
   \put(1.8,2.2){$\underbrace{\hspace*{0cm}}$}
   \put(3.1,2.2){$\underbrace{\hspace*{1.1cm}}$}
   \put(2.0,1.6){$\bar{\gamma}'_1$}
   \put(3.5,1.6){$\bar{\gamma}'_2$}
  \end{picture}}
  \caption{}
  \label{fig:23.6}
\end{center}\end{figure}

\end{proof}

\medskip
Suppose that $D_w$ is a fat-vertex of $\Lambda(t_0)$ such that there are no loops based on $D_w$. 
It follows from Lemma \ref{B3} that all the simple outermost edges for $D_w$ of $\Lambda(t_0)$ 
are either upper or lower. 

\begin{lem}\label{B4}
Suppose that all of the simple outermost edges for $D_w$ of $\Lambda(t_0)$ are upper (lower resp.). 
Then one of the following holds. 
\begin{enumerate}
\item For each fat-vertex $D_{w'}$ of $\Lambda(t_0)$, every simple outermost edges for 
$D_{w'}$ of $\Lambda(t_0)$ is upper (lower resp.). 
\item $\Sigma$ is modified by edge slides so that the modified graph contains an unknotted cycle. 
\end{enumerate}
\end{lem}

\begin{proof}
Since the arguments are symmetric, we may suppose that all the simple outermost edges for $D_w$ 
of $\Lambda(t_0)$ are upper. 
Let $\gamma$ be a simple outermost edge for $D_w$ of $\Lambda(t_0)$. 
Note that $\gamma$ is upper. 
Suppose that there is a fat-vertex $D_{w'}$ such that $\Lambda(t_0)$ contains 
a lower simple outermost edge $\gamma'$ for $D_w$. 
Let $\delta_{\gamma}$ ($\delta_{\gamma'}$ resp.) be the outermost disk for $(D_w,\gamma)$ 
($(D_{w'},\gamma')$ resp.). 
Let $\sigma$ ($\sigma'$ resp.) be the edge of $\Sigma$ with $h^1_{\sigma} \supset D_w$ 
($h^1_{\sigma'} \supset D_{w'}$ resp.). 
Let $\bar{\gamma}$ ($\bar{\gamma}'$ resp.) be a union of the components obtained by cutting $\sigma$ 
by the two fat-vertices of $\Lambda(t_0)$ incident to $\gamma$ ($\gamma'$ resp.) 
such that a $1$-handle correponding to each component intersects 
$\partial \delta_{\gamma}\setminus \gamma$ 
($\partial \delta_{\gamma'}\setminus \gamma'$ resp.). 
Then $\bar{\gamma}$ ($\bar{\gamma}'$ resp.) satisfies one of the conditions $(1)$, $(2)$ and $(3)$ 
in the proof of Lemma \ref{B3}. The proof of Lemma \ref{B4} is divided into the following cases. 

\medskip
\noindent\textit{Case A}.\ \ \ \ 
$\bar{\gamma}\cap \bar{\gamma}'=\emptyset$. 

Then we have the following six cases. In each case, we can slide 
(a component of) $\bar{\gamma}$ ($\bar{\gamma}'$ resp.) to $\gamma$ ($\gamma'$ resp.) along the disk 
$\delta_{\gamma}$ ($\delta_{\gamma'}$ resp.). Moreover, we can isotope $\sigma$ and $\sigma'$ slightly 
to reduce $(W_{\Sigma},n_{\Sigma})$ is reduced, a contradiction. 

\medskip
\noindent\textit{Case A}-(1)-(1).\ \ \ \ 
Both $\bar{\gamma}$ and $\bar{\gamma}'$ satisfy the condition $(1)$. 

\medskip
See Figure \ref{fig:CaseA-(1)-(1)}. 

\begin{figure}[htb]\begin{center}
  {\unitlength=1cm
  \begin{picture}(12,4)
   \put(0,0){\includegraphics[keepaspectratio,width=5cm]{%
   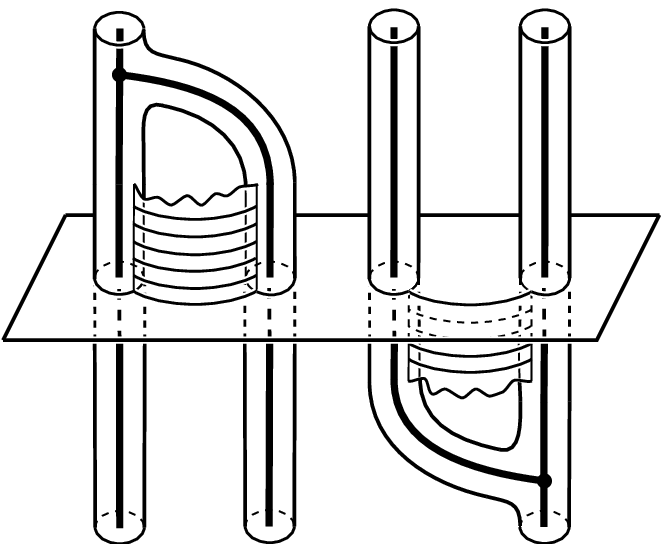}}
   \put(7.0,0){\includegraphics[keepaspectratio,width=5cm]{%
   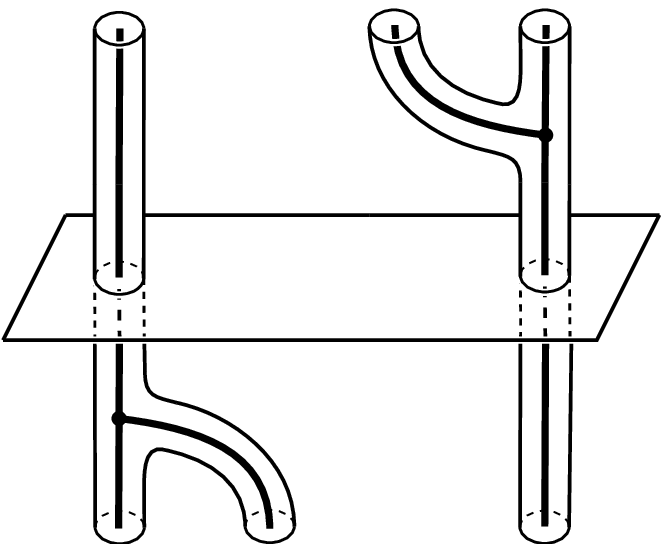}}
   \put(5.7,1.9){$\Longrightarrow$}
   \put(1.4,1.6){$\gamma$}
   \put(3.4,1.9){$\gamma'$}
  \end{picture}}
  \caption{}
  \label{fig:CaseA-(1)-(1)}
\end{center}\end{figure}

\medskip
\noindent\textit{Case A}-(1)-(2).\ \ \ \ 
Either $\bar{\gamma}$ or $\bar{\gamma}'$, say $\bar{\gamma}$, satisfies the condition $(1)$ and 
$\bar{\gamma}'$ satisfies the condition $(2)$. 

\medskip
See Figure \ref{fig:CaseA-(1)-(2)}. 

\begin{figure}[htb]\begin{center}
  {\unitlength=1cm
  \begin{picture}(12,4)
   \put(0,0){\includegraphics[keepaspectratio,width=5cm]{%
   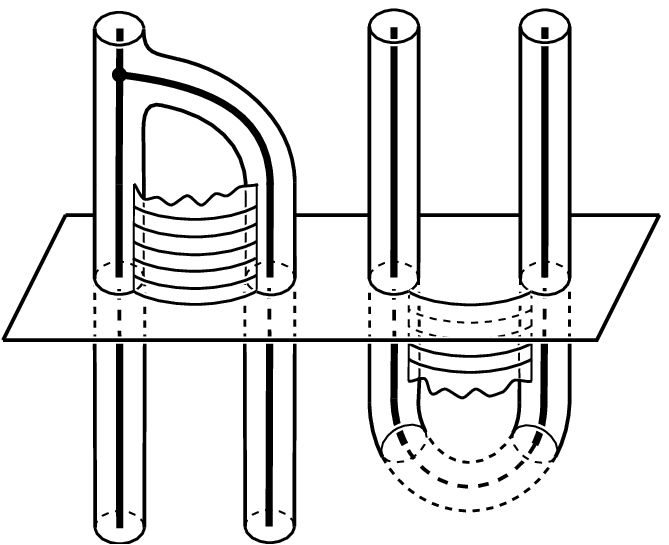}}
   \put(7.0,0){\includegraphics[keepaspectratio,width=5cm]{%
   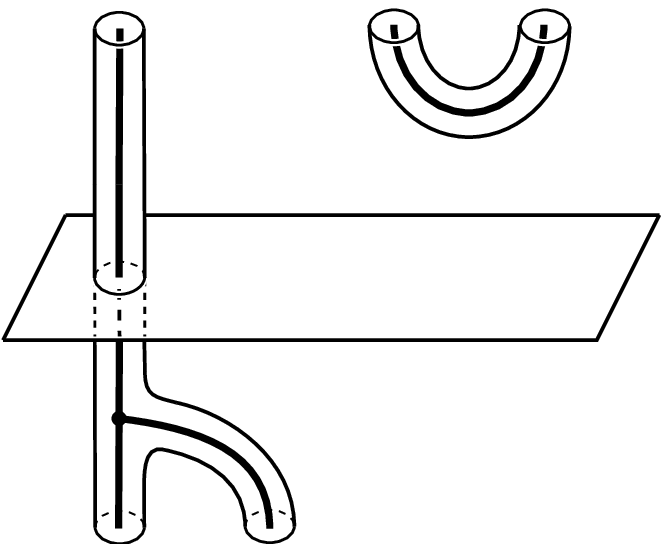}}
   \put(5.7,1.9){$\Longrightarrow$}
   \put(1.4,1.6){$\gamma$}
   \put(3.4,1.9){$\gamma'$}
  \end{picture}}
  \caption{}
  \label{fig:CaseA-(1)-(2)}
\end{center}\end{figure}

\medskip
\noindent\textit{Case A}-(1)-(3).\ \ \ \ 
Either $\bar{\gamma}$ or $\bar{\gamma}'$, say $\bar{\gamma}$, satisfies the condition $(1)$ and 
$\bar{\gamma}'$ satisfies the condition $(3)$. 

\medskip
See Figure \ref{fig:CaseA-(1)-(3)}. 

\begin{figure}[htb]\begin{center}
  {\unitlength=1cm
  \begin{picture}(12,4)
   \put(0,0){\includegraphics[keepaspectratio,width=5cm]{%
   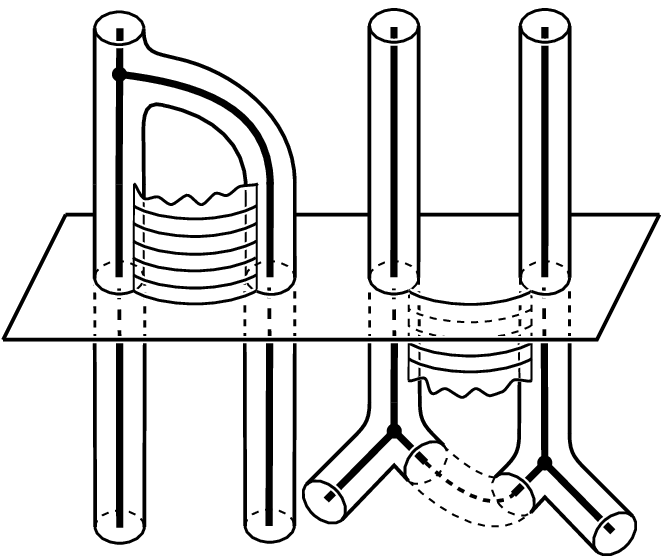}}
   \put(7.0,0){\includegraphics[keepaspectratio,width=5cm]{%
   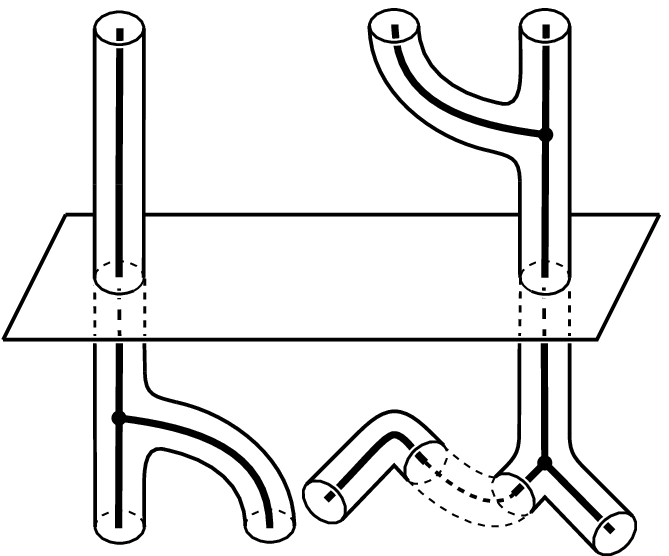}}
   \put(5.7,1.9){$\Longrightarrow$}
   \put(1.4,1.7){$\gamma$}
   \put(3.4,2.0){$\gamma'$}
  \end{picture}}
  \caption{}
  \label{fig:CaseA-(1)-(3)}
\end{center}\end{figure}

\medskip
\noindent\textit{Case A}-(2)-(2).\ \ \ \ 
Both $\bar{\gamma}$ and $\bar{\gamma}'$ satisfy the condition $(2)$. 

\medskip
See Figure \ref{fig:CaseA-(2)-(2)}. 

\begin{figure}[htb]\begin{center}
  {\unitlength=1cm
  \begin{picture}(12,4)
   \put(0,0){\includegraphics[keepaspectratio,width=5cm]{%
   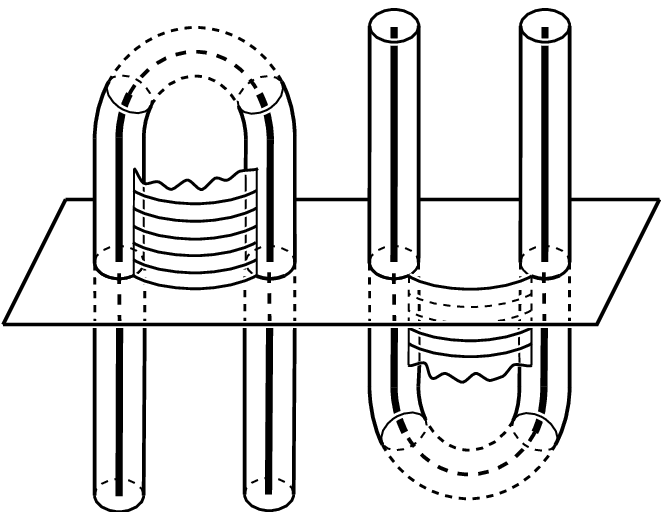}}
   \put(7.0,0){\includegraphics[keepaspectratio,width=5cm]{%
   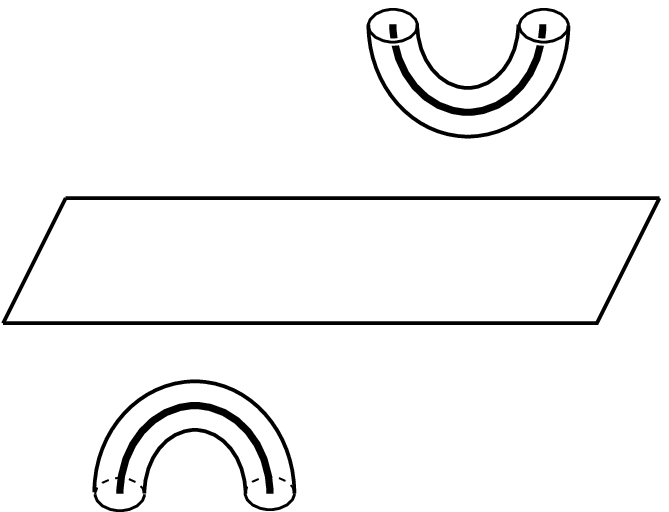}}
   \put(5.7,1.9){$\Longrightarrow$}
   \put(1.4,1.5){$\gamma$}
   \put(3.4,1.8){$\gamma'$}
  \end{picture}}
  \caption{}
  \label{fig:CaseA-(2)-(2)}
\end{center}\end{figure}

\medskip
\noindent\textit{Case A}-(2)-(3).\ \ \ \ 
Either $\bar{\gamma}$ or $\bar{\gamma}'$, say $\bar{\gamma}$, satisfies the condition $(2)$ and 
$\bar{\gamma}'$ satisfies the condition $(3)$. 

\medskip
See Figure \ref{fig:CaseA-(2)-(3)}. 

\begin{figure}[htb]\begin{center}
  {\unitlength=1cm
  \begin{picture}(12,4)
   \put(0,0){\includegraphics[keepaspectratio,width=5cm]{%
   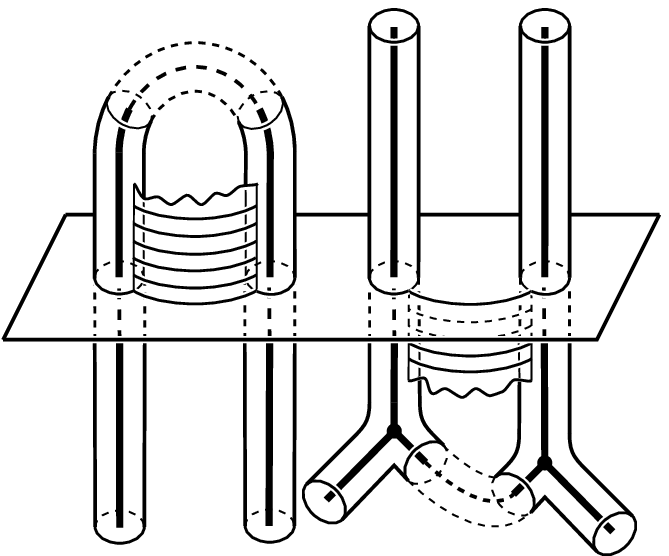}}
   \put(7.0,0){\includegraphics[keepaspectratio,width=5cm]{%
   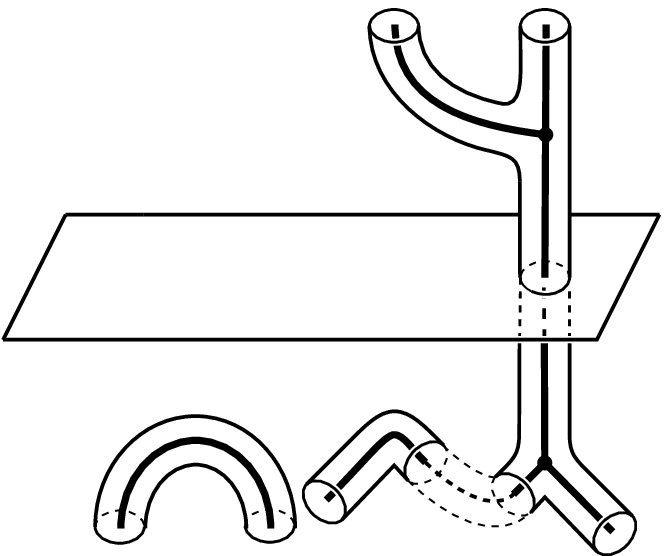}}
   \put(5.7,1.9){$\Longrightarrow$}
   \put(1.4,1.7){$\gamma$}
   \put(3.4,2.0){$\gamma'$}
  \end{picture}}
  \caption{}
  \label{fig:CaseA-(2)-(3)}
\end{center}\end{figure}

\medskip
\noindent\textit{Case A}-(3)-(3).\ \ \ \ 
Both $\bar{\gamma}$ and $\bar{\gamma}'$ satisfy the condition $(3)$. 

\medskip

See Figure \ref{fig:CaseA-(3)-(3)}. 

\begin{figure}[htb]\begin{center}
  {\unitlength=1cm
  \begin{picture}(12,4)
   \put(0,0){\includegraphics[keepaspectratio,width=5cm]{%
   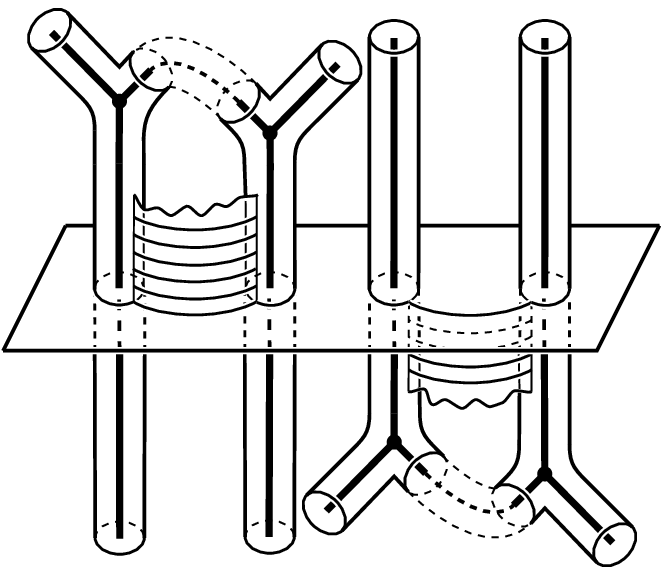}}
   \put(7.0,0){\includegraphics[keepaspectratio,width=5cm]{%
   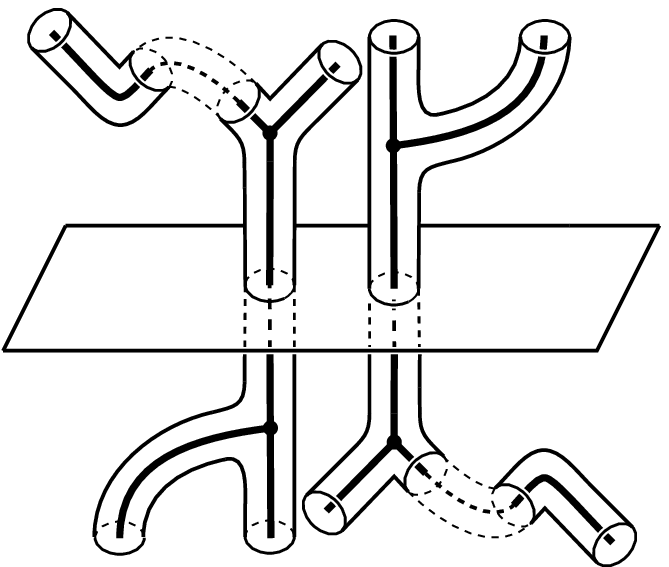}}
   \put(5.7,1.9){$\Longrightarrow$}
   \put(1.4,1.7){$\gamma$}
   \put(3.4,2.0){$\gamma'$}
  \end{picture}}
  \caption{}
  \label{fig:CaseA-(3)-(3)}
\end{center}\end{figure}

\medskip
\noindent\textit{Case B}.\ \ \ \ 
$\bar{\gamma}\cap \bar{\gamma}'\ne \emptyset$. 

\medskip
\noindent\textit{Case B}-(1)-(1).\ \ \ \ 
Both $\bar{\gamma}$ and $\bar{\gamma}'$ satisfy the condition $(1)$. 

\medskip
We first suppose that $\mathrm{int} (\bar{\gamma})\cap \mathrm{int} (\bar{\gamma}')=\emptyset$. Then we can slide 
$\bar{\gamma}$ ($\bar{\gamma}'$ resp.) to $\gamma$ ($\gamma'$ resp.) along the disk $\delta_{\gamma}$ 
($\delta_{\gamma'}$ resp.). If $\partial \gamma=\partial \gamma'(=\{ w,w'\} )$, then $\bar{\gamma}\cup \bar{\gamma}'$ 
composes an unknotted cycle and hence Lemma \ref{B4} holds (\textit{cf}. Figure \ref{fig:23a}). 
Otherwise, we can perform a Whitehead move on $\Sigma$ and hence we can reduce $(W_{\Sigma},n_{\Sigma})$, 
a contradiction (\textit{cf}. Figure \ref{fig:23.1a}). 

\begin{figure}[htb]\begin{center}
  {\unitlength=1cm
  \begin{picture}(9.6,5.6)
   \put(0,1.4){\includegraphics[keepaspectratio,width=4.0cm]{%
   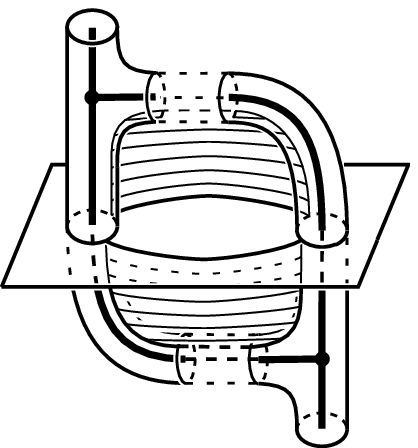}}
   \put(5.6,1.4){\includegraphics[keepaspectratio,width=4.0cm]{%
   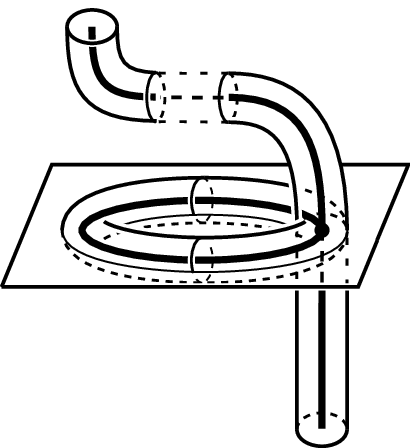}}
   \put(3.5,0.4){\includegraphics[keepaspectratio,height=.6cm]{%
   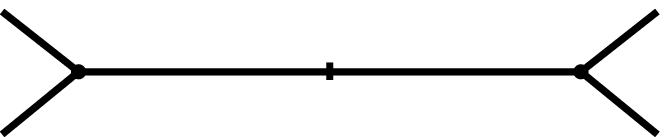}}
   \put(-.9,3.1){$P(t_0)$}
   \put(4.57,3.5){$\Longrightarrow$}
   \put(4.9,.35){$w$}
   \put(4.3,1.3){$\bar{\gamma}$}
   \put(3.85,.9){$\overbrace{\hspace*{1.1cm}}$}
   \put(5.5,1.3){$\bar{\gamma}'$}
   \put(5.06,.9){$\overbrace{\hspace*{1.1cm}}$}
  \end{picture}}
  \caption{}
  \label{fig:23a}
\end{center}\end{figure}

\begin{figure}[htb]\begin{center}
  {\unitlength=1cm
  \begin{picture}(10.5,8.0)
   \put(0,5){\includegraphics[keepaspectratio,width=4.0cm]{%
               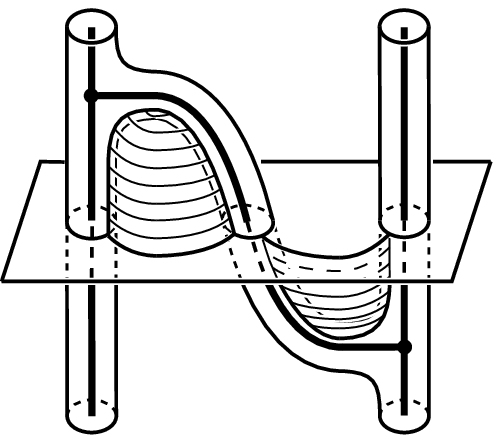}}
   \put(6.6,5){\includegraphics[keepaspectratio,width=4.0cm]{%
               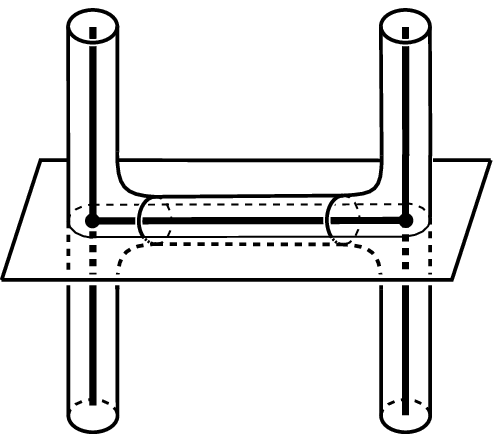}}
   \put(6.6,0.0){\includegraphics[keepaspectratio,width=4.0cm]{%
               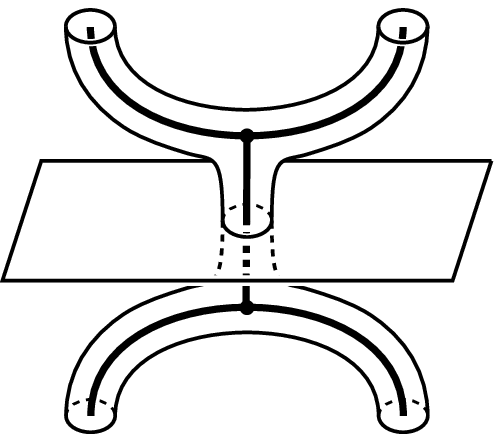}}
   \put(1.5,2.3){\includegraphics[keepaspectratio,height=.6cm]{%
               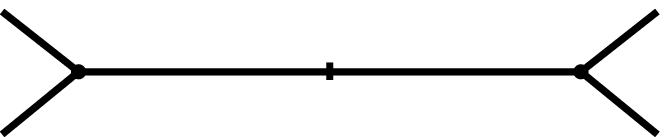}}
   \put(5.0,6.6){$\Longrightarrow$}
   \put(8.5,4.0){$\Downarrow$}
%
   \put(2.3,6.7){$D_w$}
   \put(2.9,2.25){$w$}
   \put(2.3,3.2){$\bar{\gamma}$}
   \put(1.9,2.8){$\overbrace{\hspace*{1.1cm}}$}
   \put(3.5,3.2){$\bar{\gamma}'$}
   \put(3.1,2.8){$\overbrace{\hspace*{1.1cm}}$}
  \end{picture}}
  \caption{}
  \label{fig:23.1a}
\end{center}\end{figure}

\medskip
We next suppose that $\mathrm{int} (\bar{\gamma})\cap \mathrm{int} (\bar{\gamma}')\ne \emptyset$. 
Then there are two possibilities: $(1)$ $\bar{\gamma}\subset \bar{\gamma}'$ or $\bar{\gamma}'\subset \bar{\gamma}$, 
say the 
latter holds and $(2)$ $\bar{\gamma}\not\subset \bar{\gamma}'$ and $\bar{\gamma}'\not\subset \bar{\gamma}$. 
In each case, 
we can slide $\bar{\gamma}$ to $\gamma$ along the disk $\delta_{\gamma}$. Moreover, we can isotope 
$\sigma$ slightly to reduce $(W_{\Sigma},n_{\Sigma})$, a contradiction 
(\textit{cf}. Figures \ref{fig:CaseB-(1)-(1)-(a)} and \ref{fig:CaseB-(1)-(1)-(b)}). 

\begin{figure}[htb]\begin{center}
  {\unitlength=1cm
  \begin{picture}(12,5.8)
   \put(-.4,1.4){\includegraphics[keepaspectratio,width=6cm]{%
   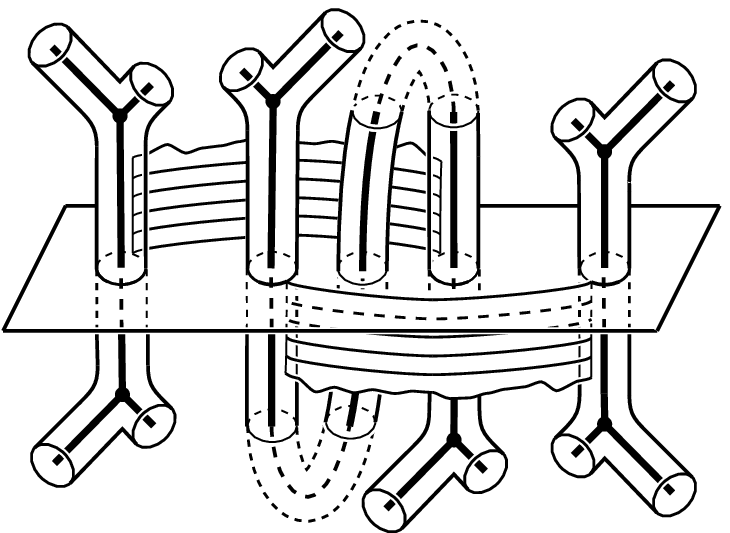}}
   \put(6.4,1.4){\includegraphics[keepaspectratio,width=6cm]{%
   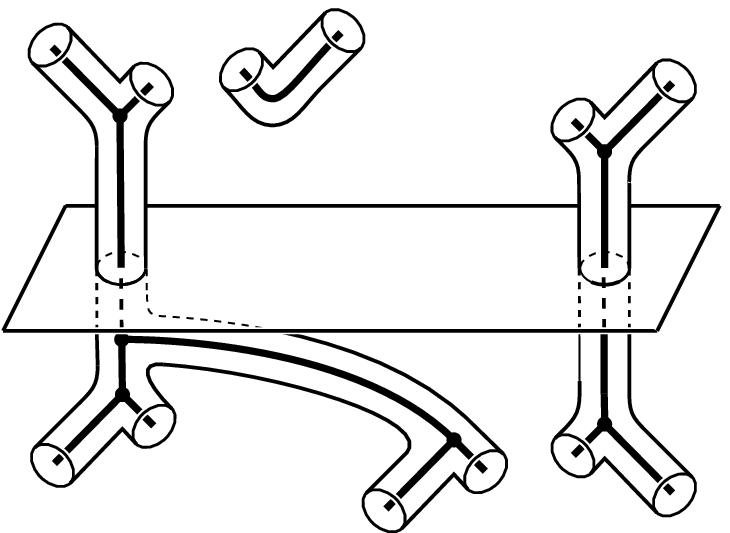}}
   \put(4.5,0.4){\includegraphics[keepaspectratio,height=.6cm]{%
   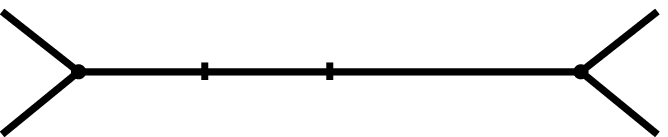}}
   \put(5.7,3.4){$\Longrightarrow$}
   \put(1.1,3.5){$\gamma$}
   \put(3.8,3.5){$\gamma'$}
   \put(6.1,0.85){$w$}
   \put(5.3,1.2){$\bar{\gamma}$}
   \put(4.9,0.85){$\overbrace{\hspace{1cm}}$}
   \put(5.0,0.3){$w'$}
   \put(5.5,0.6){$\underbrace{\hspace{1.6cm}}$}
   \put(6.15,0){$\bar{\gamma}'$}
  \end{picture}}
  \caption{}
  \label{fig:CaseB-(1)-(1)-(a)}
\end{center}\end{figure}

\begin{figure}[htb]\begin{center}
  {\unitlength=1cm
  \begin{picture}(12,5.8)
   \put(-.4,1.4){\includegraphics[keepaspectratio,width=6cm]{%
   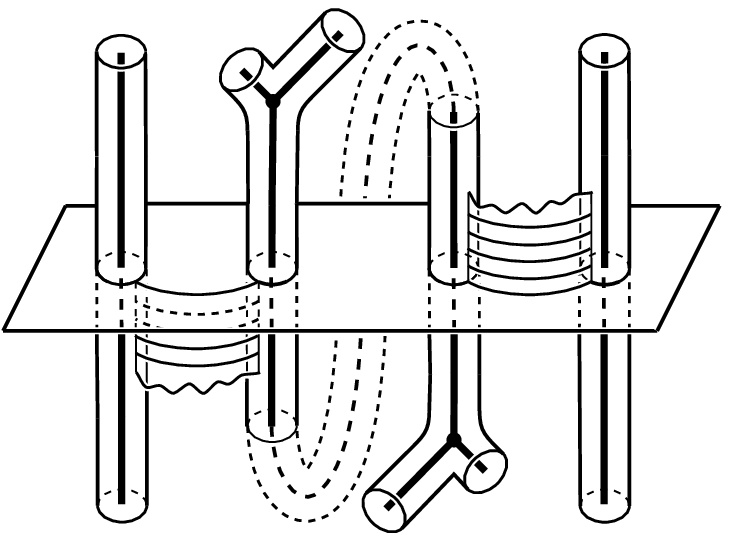}}
   \put(6.4,1.4){\includegraphics[keepaspectratio,width=6cm]{%
   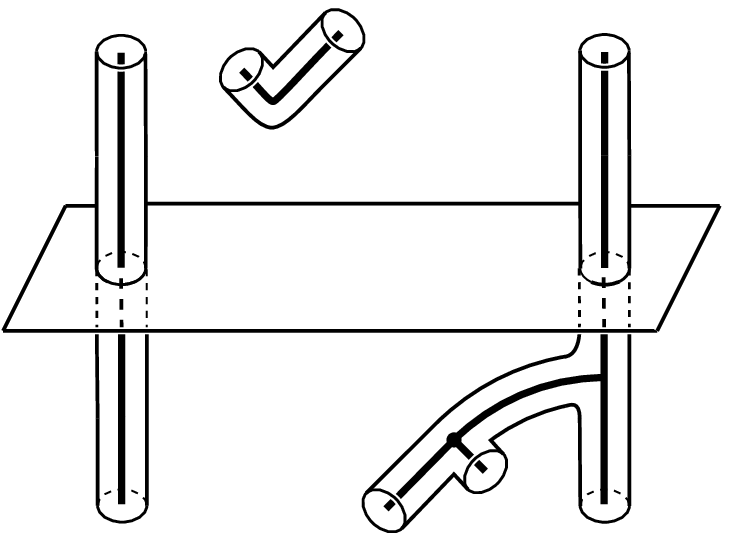}}
   \put(4.5,0.4){\includegraphics[keepaspectratio,height=.6cm]{%
   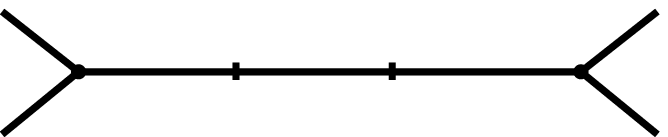}}
   \put(5.7,3.4){$\Longrightarrow$}
   \put(1.1,3.5){$\gamma'$}
   \put(3.9,3.2){$\gamma$}
   \put(5.7,0.85){$w'$}
   \put(5.1,1.2){$\bar{\gamma}$}
   \put(4.9,0.85){$\overbrace{\hspace{.5cm}}$}
   \put(6.3,0.3){$w$}
   \put(4.9,0.6){$\underbrace{\hspace{1.4cm}}$}
   \put(5.5,0){$\bar{\gamma}'$}
  \end{picture}}
  \caption{}
  \label{fig:CaseB-(1)-(1)-(b)}
\end{center}\end{figure}

\medskip
\noindent\textit{Case B}-(1)-(2).\ \ \ \ 
Either $\bar{\gamma}$ or $\bar{\gamma}'$, say $\bar{\gamma}$, satisfies the condition $(1)$ and 
$\bar{\gamma}'$ satisfies the condition $(2)$. 

\medskip
We first suppose that $\mathrm{int} (\bar{\gamma})\cap \mathrm{int} (\bar{\gamma}')=\emptyset$. 
Then we can slide $\bar{\gamma}$ ($\bar{\gamma}'$ resp.) to $\gamma$ ($\gamma'$ resp.) along the disk 
$\delta_{\gamma}$ 
($\delta_{\gamma'}$ resp.). Moreover, we can isotope $\sigma$ slightly to reduce 
$(W_{\Sigma},n_{\Sigma})$, a contradiction 
(\textit{cf}. Figure \ref{fig:CaseB-(1)-(2)-(b)}). 

\begin{figure}[htb]\begin{center}
  {\unitlength=1cm
  \begin{picture}(12.0,4.5)
   \put(0,1.0){\includegraphics[keepaspectratio,width=5.5cm]{%
               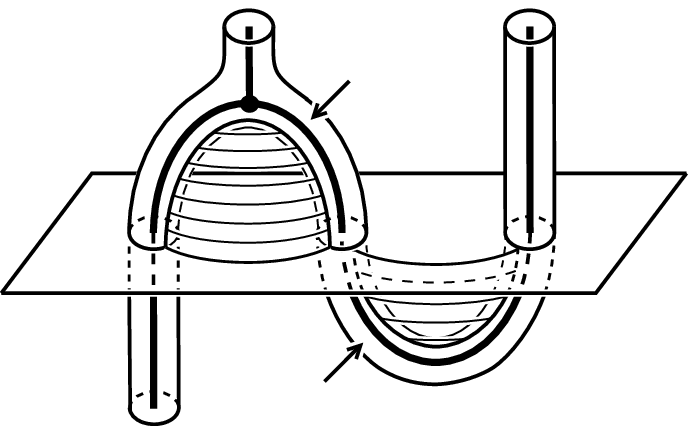}}
   \put(6.6,1.0){\includegraphics[keepaspectratio,width=5.5cm]{%
               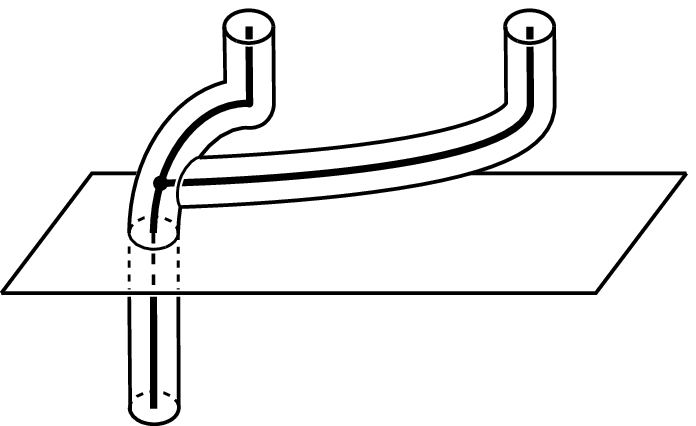}}
   \put(4.2,0){\includegraphics[keepaspectratio,height=.6cm]{%
               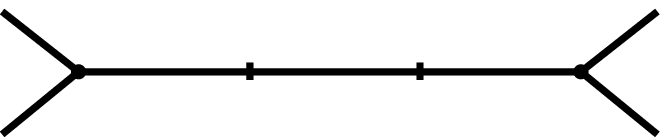}}
   \put(5.8,2.4){$\Longrightarrow$}
   \put(2.9,3.8){$\bar{\gamma}$}
   \put(2.2,1.0){$\bar{\gamma}'$}
   \put(3.0,2.6){$D_w$}
   \put(4.85,0.85){$\bar{\gamma}$}
   \put(4.6,0.5){$\overbrace{\hspace{0.5cm}}$}
   \put(5.65,0.85){$\bar{\gamma}'$}
   \put(5.4,0.5){$\overbrace{\hspace{0.5cm}}$}
   \put(5.2,-0.1){$w$}
  \end{picture}}
  \caption{}
  \label{fig:CaseB-(1)-(2)-(b)}
\end{center}\end{figure}

\medskip
We next suppose that $\mathrm{int} (\bar{\gamma})\cap \mathrm{int} (\bar{\gamma}')\ne \emptyset$. 
Note that it is impossible that $\bar{\gamma}\subset \bar{\gamma}'$. Hence there are 
two possibilities: $\bar{\gamma}'\subset \bar{\gamma}$ and $\bar{\gamma}'\not\subset \bar{\gamma}$. In each case, 
we can slide $\bar{\gamma}$ to $\gamma$ along the disk $\delta_{\gamma}$. Moreover, we can isotope 
$\sigma$ slightly to reduce $(W_{\Sigma},n_{\Sigma})$, a contradiction 
(\textit{cf}. Figures \ref{fig:CaseB-(1)-(2)-(a)} and \ref{fig:CaseB-(1)-(2)-(c)}). 

\begin{figure}[htb]\begin{center}
  {\unitlength=1cm
  \begin{picture}(12,5.8)
   \put(0,1.4){\includegraphics[keepaspectratio,width=5cm]{%
   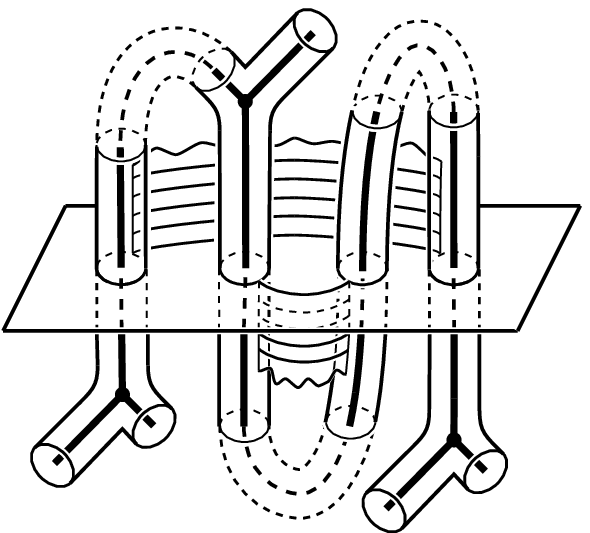}}
   \put(7,1.4){\includegraphics[keepaspectratio,width=5cm]{%
   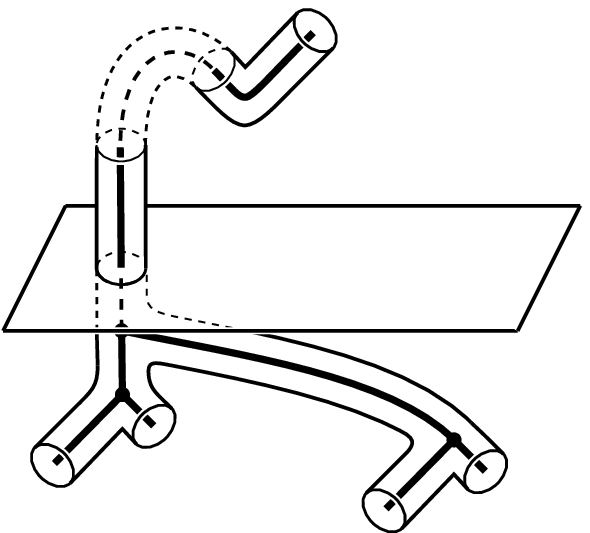}}
   \put(4.5,0.4){\includegraphics[keepaspectratio,height=.6cm]{%
   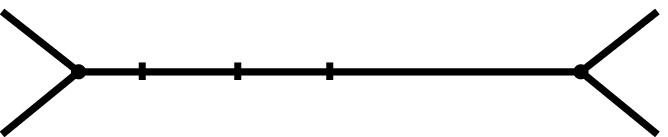}}
   \put(5.7,3.4){$\Longrightarrow$}
   \put(1.4,3.6){$\gamma$}
   \put(2.5,3.54){$\gamma'$}
   \put(6.1,0.85){$w$}
   \put(5.4,1.2){$\bar{\gamma}$}
   \put(4.9,0.85){$\overbrace{\hspace{1.1cm}}$}
%
   \put(5.1,0.6){$\underbrace{\hspace{.0cm}}$}
   \put(5.3,0){$\bar{\gamma}'$}
  \end{picture}}
  \caption{}
  \label{fig:CaseB-(1)-(2)-(a)}
\end{center}\end{figure}

\begin{figure}[htb]\begin{center}
  {\unitlength=1cm
  \begin{picture}(12,5.8)
   \put(-.4,1.4){\includegraphics[keepaspectratio,width=6cm]{%
   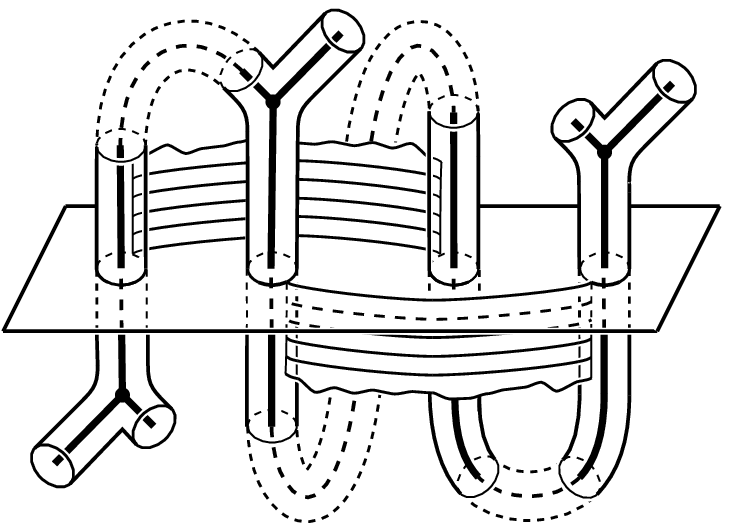}}
   \put(6.4,1.4){\includegraphics[keepaspectratio,width=6cm]{%
   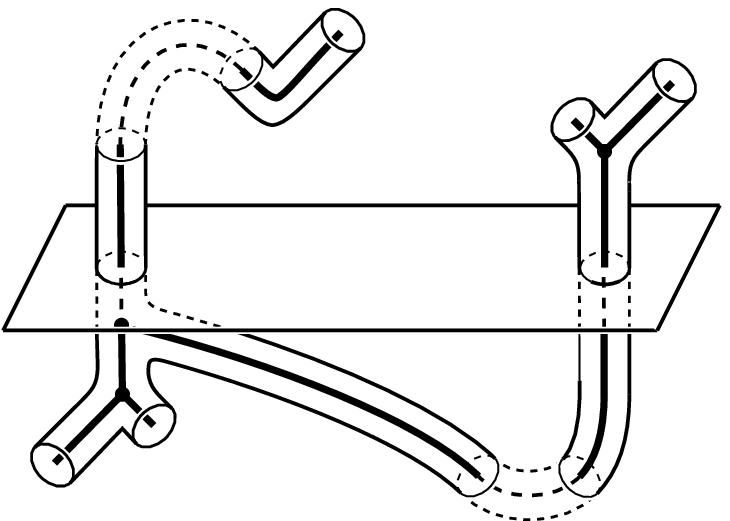}}
   \put(4.5,0.4){\includegraphics[keepaspectratio,height=.6cm]{%
   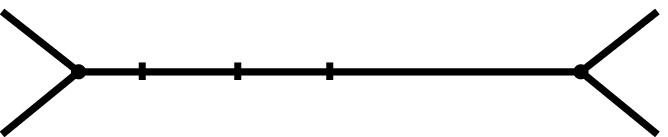}}
   \put(5.7,3.4){$\Longrightarrow$}
   \put(1.1,3.5){$\gamma$}
   \put(3.8,3.4){$\gamma'$}
   \put(5.6,0.85){$w$}
   \put(5.1,1.2){$\bar{\gamma}$}
   \put(4.9,0.85){$\overbrace{\hspace{0cm}}$}
%
   \put(5.2,0.6){$\underbrace{\hspace{.8cm}}$}
   \put(5.5,0){$\bar{\gamma}'$}
  \end{picture}}
  \caption{}
  \label{fig:CaseB-(1)-(2)-(c)}
\end{center}\end{figure}

\medskip
\noindent\textit{Case B}-(1)-(3).\ \ \ \ 
Either $\bar{\gamma}$ or $\bar{\gamma}'$, say $\bar{\gamma}$, satisfies the condition $(1)$ and 
$\bar{\gamma}'$ satisfies the condition $(3)$. 

\medskip
Let $\bar{\gamma}'_1$ and $\bar{\gamma}'_2$ be the components of $\bar{\gamma}'$ with 
$h^1_{\bar{\gamma}'_1}\supset D_{w'}$. 

We first suppose that $\bar{\gamma}\subset \bar{\gamma}_1'$. Then we can slide $\bar{\gamma}'_1$ into $\gamma'$ 
along the disk $\delta_{\gamma'}$. Moreover, we can isotope 
$\sigma$ slightly to reduce $(W_{\Sigma},n_{\Sigma})$, a contradiction 
(\textit{cf}. Figure \ref{fig:CaseB-(1)-(3)-(a)}). 

\begin{figure}[htb]\begin{center}
  {\unitlength=1cm
  \begin{picture}(12,5.8)
   \put(0.0,1.4){\includegraphics[keepaspectratio,width=5cm]{%
   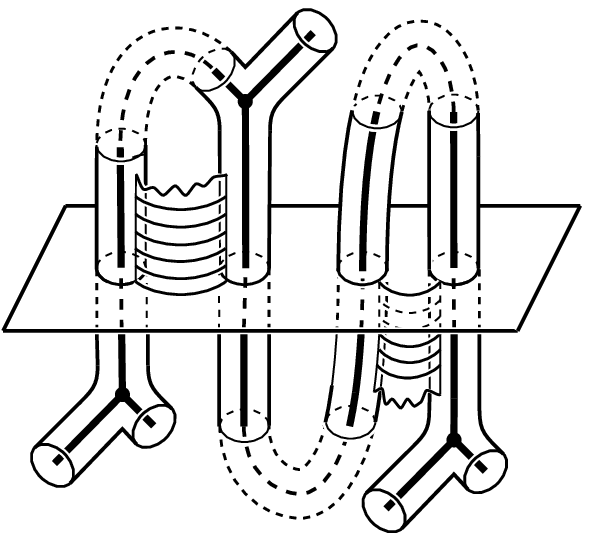}}
   \put(7.0,1.4){\includegraphics[keepaspectratio,width=5cm]{%
   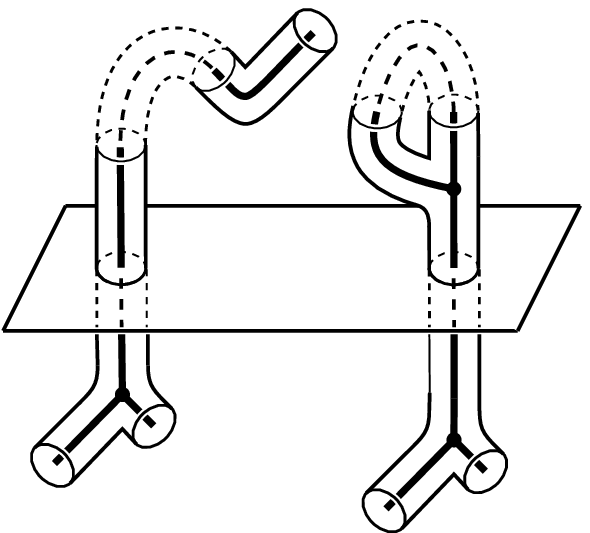}}
   \put(4.5,0.4){\includegraphics[keepaspectratio,height=.6cm]{%
   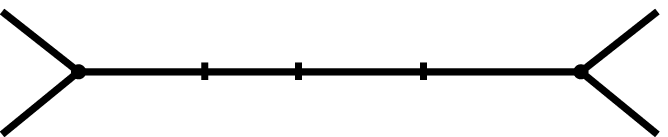}}
   \put(5.7,3.4){$\Longrightarrow$}
   \put(1.4,3.2){$\gamma$}
   \put(3.36,3.56){$\gamma'$}
   \put(5.5,0.85){$w$}
   \put(5.0,1.2){$\bar{\gamma}$}
   \put(4.8,0.85){$\overbrace{\hspace{0cm}}$}
%
   \put(4.9,0.6){$\underbrace{\hspace{1.0cm}}$}
   \put(5.2,0){$\bar{\gamma}'_1$}
   \put(6.5,0.6){$\underbrace{\hspace{0.6cm}}$}
   \put(6.7,0){$\bar{\gamma}'_2$}
  \end{picture}}
  \caption{}
  \label{fig:CaseB-(1)-(3)-(a)}
\end{center}\end{figure}

\medskip
We next suppose that $\bar{\gamma}_1'\subset \bar{\gamma}$. Then there are 
two possibilities: $\bar{\gamma}'_2\cap \bar{\gamma}=\emptyset$ and $\bar{\gamma}_2'\cap\bar{\gamma}\ne \emptyset$. 
In each case, 
we can slide $\bar{\gamma}$ to $\gamma$ along the disk $\delta_{\gamma}$. Moreover, we can isotope 
$\sigma$ slightly to reduce $(W_{\Sigma},n_{\Sigma})$, a contradiction 
(\textit{cf}. Figures \ref{fig:CaseB-(1)-(3)-(b)} and \ref{fig:CaseB-(1)-(3)-(c)}). 

\begin{figure}[htb]\begin{center}
  {\unitlength=1cm
  \begin{picture}(12,5.8)
   \put(-.4,1.4){\includegraphics[keepaspectratio,width=6cm]{%
   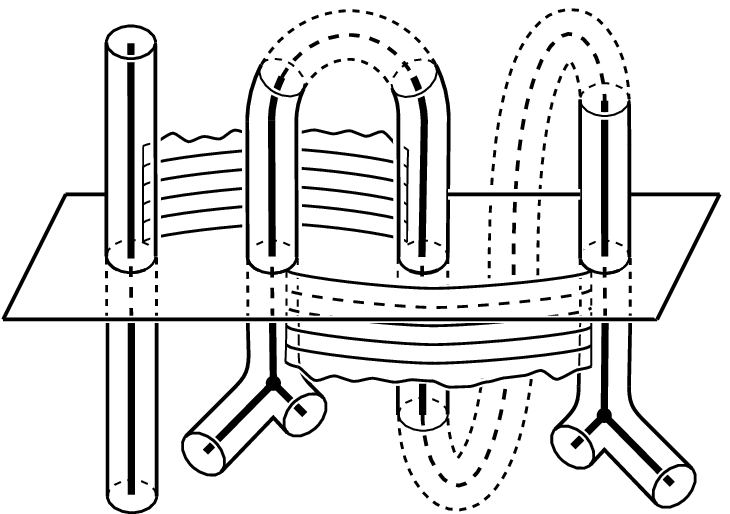}}
   \put(6.4,1.4){\includegraphics[keepaspectratio,width=6cm]{%
   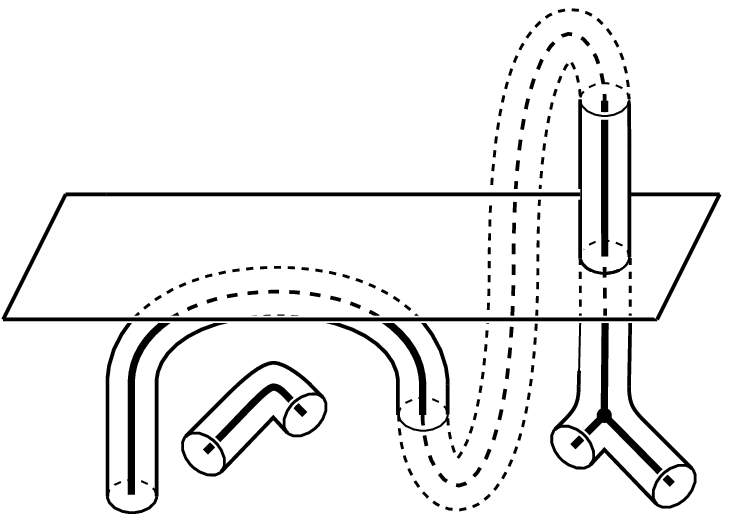}}
   \put(4.5,0.4){\includegraphics[keepaspectratio,height=.6cm]{%
   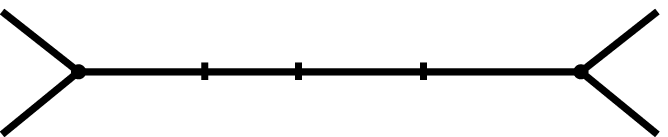}}
   \put(5.7,3.4){$\Longrightarrow$}
   \put(1.1,3.5){$\gamma$}
   \put(3.35,3.4){$\gamma'$}
   \put(5.8,0.85){$w$}
   \put(5.2,1.2){$\bar{\gamma}'_1$}
   \put(4.8,0.85){$\overbrace{\hspace{0cm}}$}
   \put(6.5,0.85){$\overbrace{\hspace{0.6cm}}$}
   \put(6.7,1.2){$\bar{\gamma}'_2$}
%
   \put(4.9,0.6){$\underbrace{\hspace{1.0cm}}$}
   \put(5.3,0){$\bar{\gamma}$}
  \end{picture}}
  \caption{}
  \label{fig:CaseB-(1)-(3)-(b)}
\end{center}\end{figure}

\begin{figure}[htb]\begin{center}
  {\unitlength=1cm
  \begin{picture}(12,5.8)
   \put(-.4,1.4){\includegraphics[keepaspectratio,width=6cm]{%
   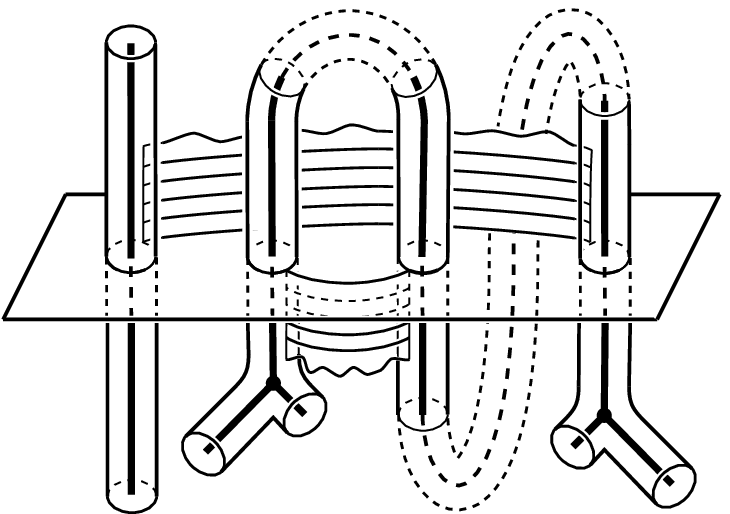}}
   \put(6.4,1.4){\includegraphics[keepaspectratio,width=6cm]{%
   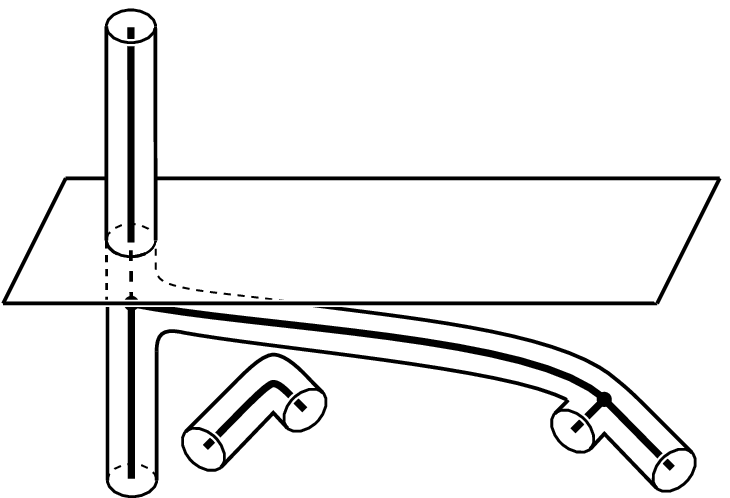}}
   \put(4.5,0.4){\includegraphics[keepaspectratio,height=.6cm]{%
   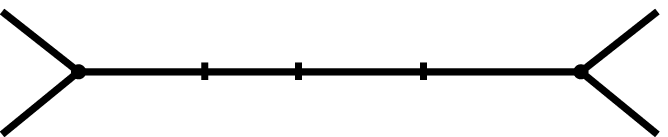}}
   \put(5.7,3.4){$\Longrightarrow$}
   \put(1.1,3.5){$\gamma$}
   \put(2.3,3.4){$\gamma'$}
%
   \put(5.6,1.2){$\bar{\gamma}$}
   \put(4.9,0.85){$\overbrace{\hspace{1.6cm}}$}
%
   \put(4.9,0.6){$\underbrace{\hspace{.6cm}}$}
   \put(5.1,0){$\bar{\gamma}'_1$}
   \put(5.9,0.6){$\underbrace{\hspace{1.2cm}}$}
   \put(6.4,0){$\bar{\gamma}'_2$}
  \end{picture}}
  \caption{}
  \label{fig:CaseB-(1)-(3)-(c)}
\end{center}\end{figure}

\medskip
\noindent\textit{Case B}-(2)-(2).\ \ \ \ 
Both $\bar{\gamma}$ and $\bar{\gamma}'$ satisfy the condition $(2)$. 

\medskip
We first suppose that $\mathrm{int} (\bar{\gamma})\cap \mathrm{int} (\bar{\gamma}')=\emptyset$. 
Then we can slide $\bar{\gamma}$ ($\bar{\gamma}'$ resp.) to $\gamma$ ($\gamma'$ resp.) along the disk 
$\delta_{\gamma}$ 
($\delta_{\gamma'}$ resp.). Moreover, we can isotope 
$\sigma$ slightly to reduce $(W_{\Sigma},n_{\Sigma})$, a contradiction 
(\textit{cf}. Figure \ref{fig:22}). 

\begin{figure}[htb]\begin{center}
  {\unitlength=1cm
  \begin{picture}(12.0,4.6)
   \put(0,1.4){\includegraphics[keepaspectratio]{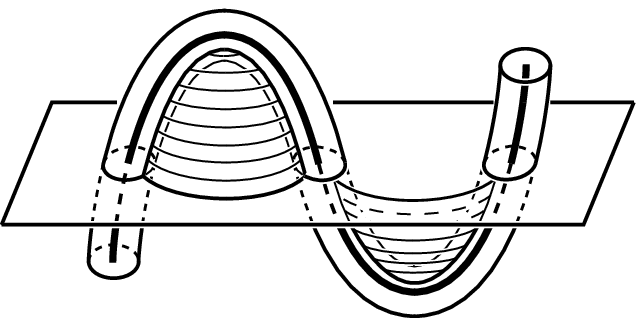}}
   \put(7.6,1.78){\includegraphics[keepaspectratio]{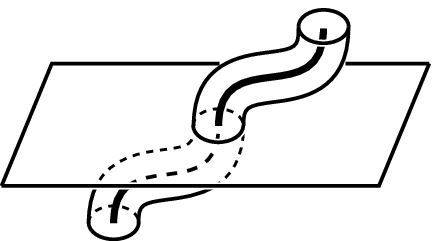}}
   \put(4.6,0.4){\includegraphics[keepaspectratio,height=.6cm]{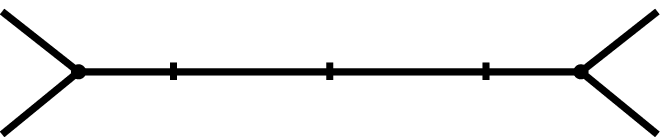}}
   \put(6.8,2.9){$\Longrightarrow$}
   \put(-.9,2.4){$P(t_0)$}
   \put(3.6,3.0){$D_w$}
   \put(1.6,2.5){$\gamma$}
   \put(4.4,2.8){$\gamma'$}
   \put(6.0,0.85){$w$}
   \put(5.35,0.6){$\underbrace{\hspace{0cm}}$}
   \put(5.7,0){$\bar{\gamma}$}
   \put(6.2,0.6){$\underbrace{\hspace{0cm}}$}
   \put(6.6,0){$\bar{\gamma}'$}
  \end{picture}}
  \caption{}
  \label{fig:22}
\end{center}\end{figure}

We next suppose that $\mathrm{int} (\bar{\gamma})\cap \mathrm{int} (\bar{\gamma}')\ne \emptyset$. 
Then there are two possibilities: $(1)$ $\bar{\gamma}\subset \bar{\gamma}'$ or $\bar{\gamma}'\subset \bar{\gamma}$, 
say the latter holds and $(2)$ $\bar{\gamma}\not\subset \bar{\gamma}'$ and $\bar{\gamma}'\not\subset \bar{\gamma}$. 
In each case, 
we can slide $\bar{\gamma}$ to $\gamma$ along the disk $\delta_{\gamma}$. Moreover, we can isotope 
$\sigma$ slightly to reduce $(W_{\Sigma},n_{\Sigma})$, a contradiction 
(\textit{cf}. Figures \ref{fig:CaseB-(2)-(2)-(a)} and \ref{fig:CaseB-(2)-(2)-(b)}).

\begin{figure}[htb]\begin{center}
  {\unitlength=1cm
  \begin{picture}(12,5.8)
   \put(-.4,1.4){\includegraphics[keepaspectratio,width=6cm]{%
   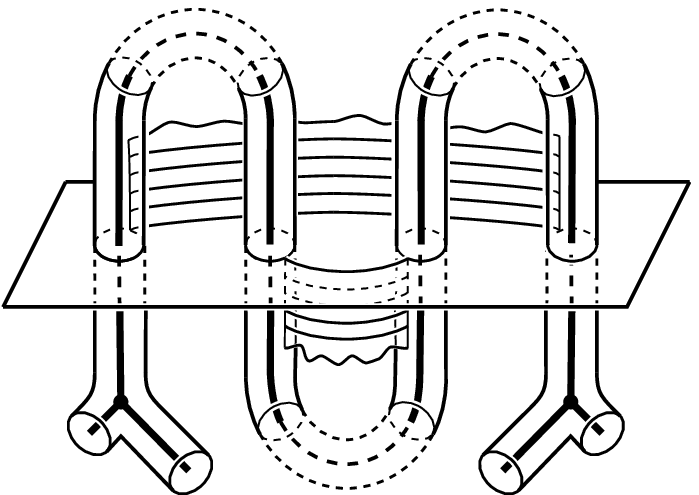}}
   \put(6.4,1.4){\includegraphics[keepaspectratio,width=6cm]{%
   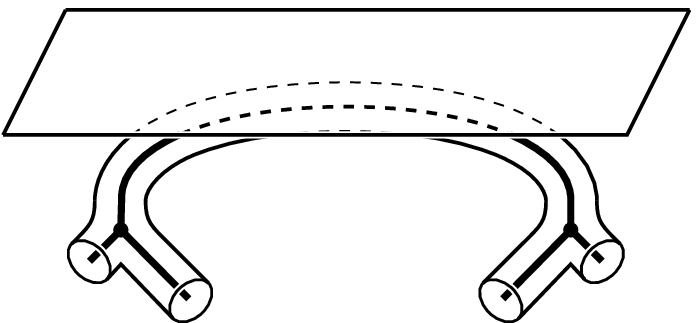}}
   \put(4.5,0.4){\includegraphics[keepaspectratio,height=.6cm]{%
   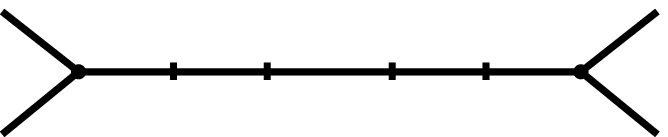}}
   \put(5.7,3.4){$\Longrightarrow$}
   \put(1.1,3.5){$\gamma$}
   \put(2.4,3.4){$\gamma'$}
%
   \put(5.9,1.2){$\bar{\gamma}$}
   \put(5.3,0.85){$\overbrace{\hspace{1.4cm}}$}
%
   \put(5.7,0.6){$\underbrace{\hspace{.6cm}}$}
   \put(5.9,0){$\bar{\gamma}'$}
  \end{picture}}
  \caption{}
  \label{fig:CaseB-(2)-(2)-(a)}
\end{center}\end{figure}

\begin{figure}[htb]\begin{center}
  {\unitlength=1cm
  \begin{picture}(12,5.8)
   \put(-.4,1.4){\includegraphics[keepaspectratio,width=6cm]{%
   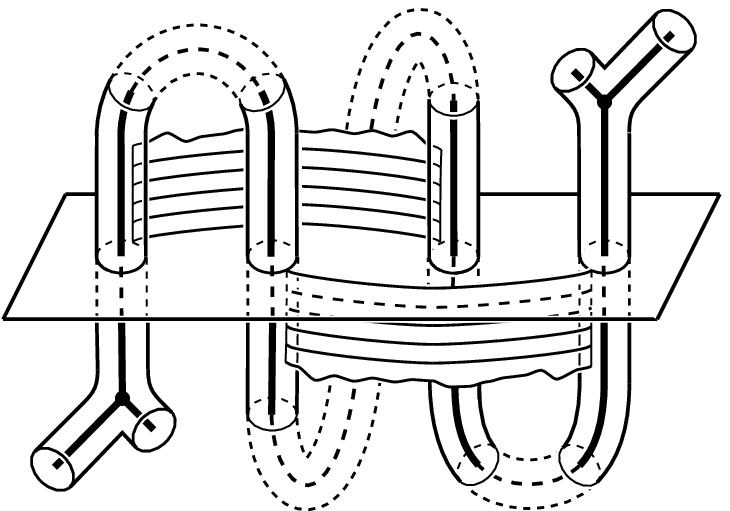}}
   \put(6.4,1.4){\includegraphics[keepaspectratio,width=6cm]{%
   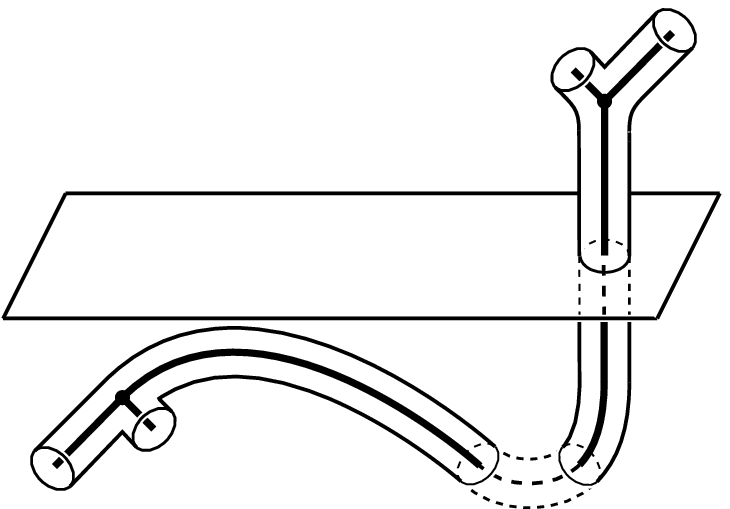}}
   \put(4.5,0.4){\includegraphics[keepaspectratio,height=.6cm]{%
   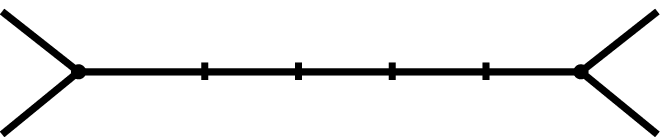}}
   \put(5.7,3.4){$\Longrightarrow$}
   \put(1.1,3.5){$\gamma$}
   \put(3.8,3.4){$\gamma'$}
%
   \put(5.9,1.2){$\bar{\gamma}$}
   \put(5.5,0.85){$\overbrace{\hspace{0.8cm}}$}
%
   \put(5.9,0.6){$\underbrace{\hspace{.9cm}}$}
   \put(6.2,0){$\bar{\gamma}'$}
  \end{picture}}
  \caption{}
  \label{fig:CaseB-(2)-(2)-(b)}
\end{center}\end{figure}

\medskip
\noindent\textit{Case B}-(2)-(3).\ \ \ \ 
Either $\bar{\gamma}$ or $\bar{\gamma}'$, say $\bar{\gamma}$, satisfies the condition $(2)$ and 
$\bar{\gamma}'$ satisfies the condition $(3)$. 

\medskip
Let $\bar{\gamma}'_1$ and $\bar{\gamma}'_2$ be the components of $\bar{\gamma}'$ with 
$\partial \bar{\gamma}'_1\supset D_{w'}$. 

We first suppose that $\mathrm{int} (\bar{\gamma})\cap \mathrm{int} (\bar{\gamma}')=\emptyset$. 
Since $\bar{\gamma}\cap \bar{\gamma'}\ne \emptyset$, we may suppose that 
$\bar{\gamma}'_1\cap \bar{\gamma}(=\partial \bar{\gamma}'_1\cap \partial \bar{\gamma})$ consists of 
a single point. Then we can slide $\bar{\gamma}$ ($\bar{\gamma}'_1$ resp.) to $\gamma$ ($\gamma'$ resp.) 
along the disk $\delta_{\gamma}$ ($\delta_{\gamma'}$ resp.). 
If $\bar{\gamma}'_2\cap \bar{\gamma}\ne \emptyset$, then $\bar{\gamma}'_2\cap \bar{\gamma}=
\partial \bar{\gamma}'_2\cap \bar{\gamma}$ consists of a single point. 
Hence $\bar{\gamma}'_1\cup \bar{\gamma}$ 
composes an unknotted cycle and hence Lemma \ref{B4} holds (\textit{cf}. Figure \ref{fig:23.4}). 
Otherwise,  we can further isotope $\Sigma$ to reduce $(W_{\Sigma},n_{\Sigma})$, a contradiction 
(\textit{cf}. Figure \ref{fig:CaseB-(2)-(3)-(a)}). 

\begin{figure}[htb]\begin{center}
  {\unitlength=1cm
  \begin{picture}(12.0,5.4)
   \put(0,1.4){\includegraphics[keepaspectratio,width=5.5cm]{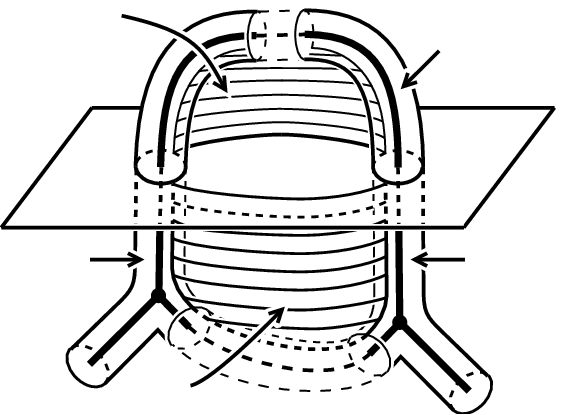}}
   \put(6.6,1.4){\includegraphics[keepaspectratio,width=5.5cm]{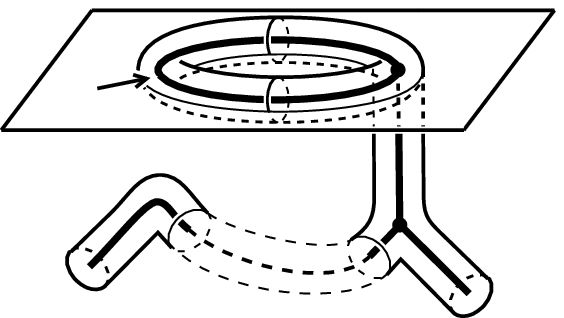}}
   \put(4.2,0.4){\includegraphics[keepaspectratio,height=.6cm]{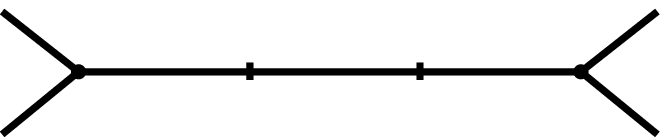}}
   \put(5.7,3.7){$\Longrightarrow$}
   \put(.8,5.3){$\delta_\gamma$}
   \put(1.5,1.5){$\delta_{\gamma'}$}
   \put(4.4,4.9){$\bar{\gamma}$}
   \put(.5,2.8){$\bar{\gamma}'_2$}
   \put(4.6,2.8){$\bar{\gamma}'_1$}
   \put(4.3,3.8){$D_w$}
   \put(7.3,3.5){$\sigma$}
   \put(5.1,0.85){$w$}
   \put(5.65,1.2){$\bar{\gamma}$}
   \put(5.4,0.85){$\overbrace{\hspace{.7cm}}$}
   \put(4.6,0.6){$\underbrace{\hspace{.7cm}}$}
   \put(6.2,0.6){$\underbrace{\hspace{.7cm}}$}
   \put(4.8,0.0){$\bar{\gamma}'_1$}
   \put(6.4,0.0){$\bar{\gamma}'_2$}
  \end{picture}}
  \caption{}
  \label{fig:23.4}
\end{center}\end{figure}

\begin{figure}[htb]\begin{center}
  {\unitlength=1cm
  \begin{picture}(12,5.8)
   \put(0.0,1.4){\includegraphics[keepaspectratio,width=5cm]{%
   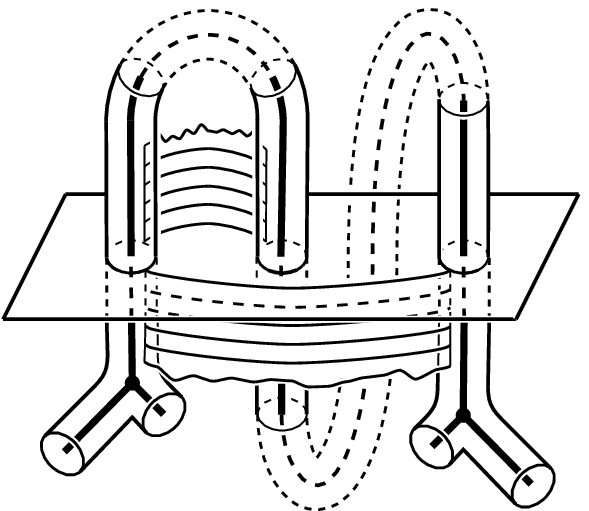}}
   \put(7.0,1.4){\includegraphics[keepaspectratio,width=5cm]{%
   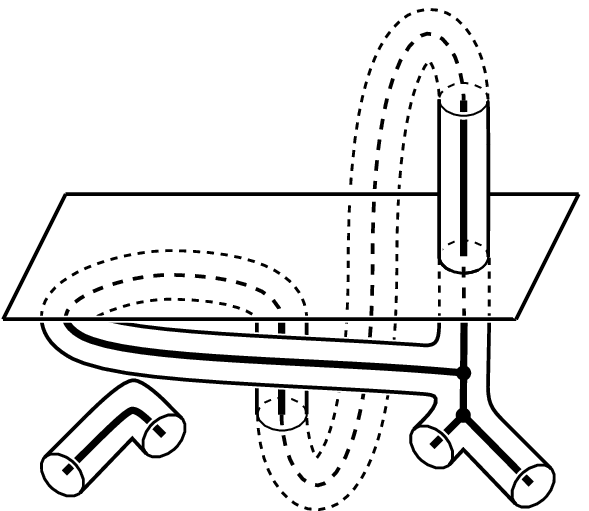}}
   \put(4.5,0.4){\includegraphics[keepaspectratio,height=.6cm]{%
   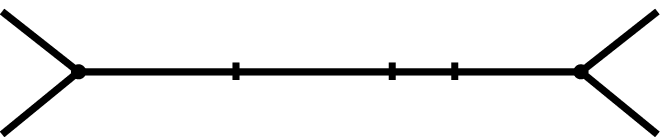}}
   \put(5.7,3.4){$\Longrightarrow$}
   \put(1.6,3.6){$\gamma$}
   \put(2.66,3.43){$\gamma'$}
%
   \put(5.9,1.2){$\bar{\gamma}$}
   \put(5.55,0.85){$\overbrace{\hspace{0.8cm}}$}
%
   \put(4.9,0.6){$\underbrace{\hspace{0cm}}$}
   \put(5.1,0){$\bar{\gamma}'_1$}
   \put(6.6,0.6){$\underbrace{\hspace{0cm}}$}
   \put(6.8,0){$\bar{\gamma}'_2$}
  \end{picture}}
  \caption{}
  \label{fig:CaseB-(2)-(3)-(a)}
\end{center}\end{figure}

We next suppose that $\mathrm{int} \bar{\gamma}\cap \mathrm{int} \bar{\gamma}'\ne \emptyset$. 
We may assume that $\mathrm{int} \bar{\gamma}\cap \mathrm{int} \bar{\gamma}'_1\ne \emptyset$. 

Then there are two possibilities: 
$\mathrm{int} \bar{\gamma}\cap \mathrm{int} \bar{\gamma}'_2= \emptyset$ and 
$\mathrm{int} \bar{\gamma}\cap \mathrm{int} \bar{\gamma}'_2\ne \emptyset$. 
In each case, we can slide $\bar{\gamma}$ to $\gamma$ along the disk 
$\delta_{\gamma}$. Moreover, we can isotope $\sigma$ and $\sigma'$ 
slightly to reduce $(W_{\Sigma},n_{\Sigma})$, a contradiction 
(\textit{cf}. Figure \ref{fig:CaseB-(2)-(3)-(b)} and Figure \ref{fig:CaseB-(2)-(3)-(c)}). 

\begin{figure}[htb]\begin{center}
  {\unitlength=1cm
  \begin{picture}(12,5.8)
   \put(-.4,1.4){\includegraphics[keepaspectratio,width=6cm]{%
   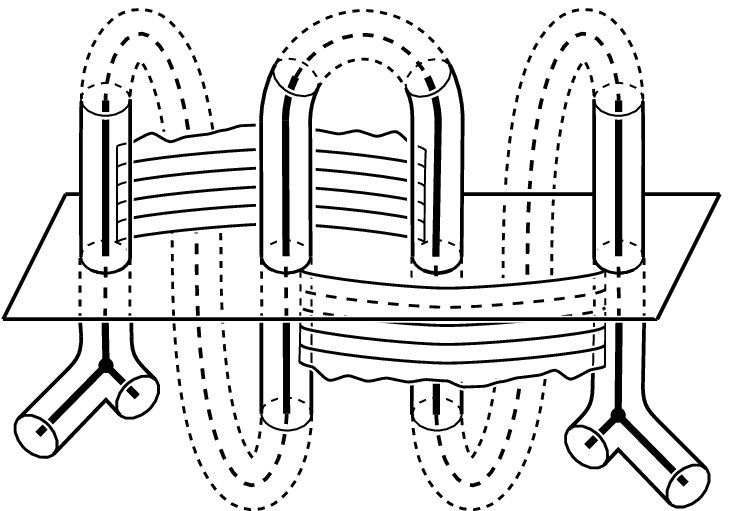}}
   \put(6.4,1.4){\includegraphics[keepaspectratio,width=6cm]{%
   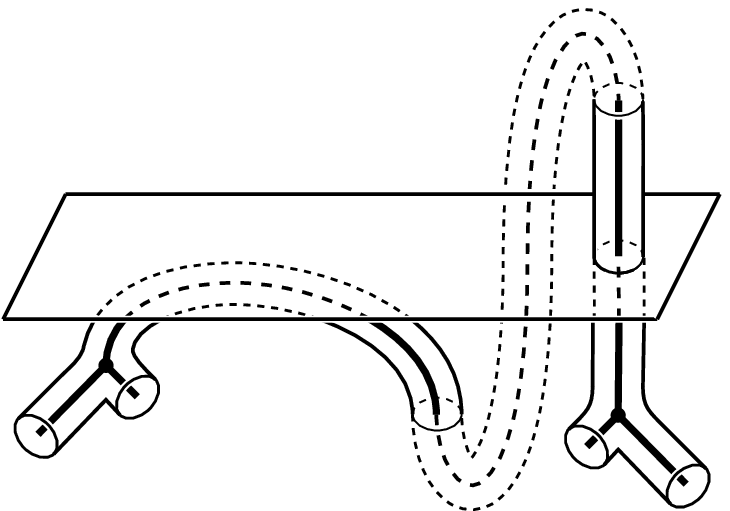}}
   \put(4.5,0.4){\includegraphics[keepaspectratio,height=.6cm]{%
   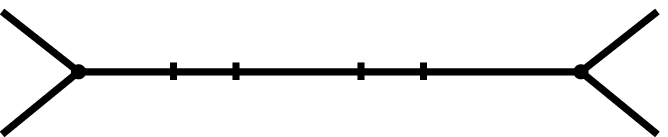}}
   \put(5.7,3.4){$\Longrightarrow$}
   \put(1.5,3.5){$\gamma$}
   \put(3.45,3.4){$\gamma'$}
%
   \put(5.7,1.2){$\bar{\gamma}$}
   \put(5.4,0.85){$\overbrace{\hspace{0.8cm}}$}
%
   \put(4.9,0.6){$\underbrace{\hspace{0cm}}$}
   \put(5.1,0){$\bar{\gamma}'_1$}
   \put(6.5,0.6){$\underbrace{\hspace{0cm}}$}
   \put(6.7,0){$\bar{\gamma}'_2$}
  \end{picture}}
  \caption{}
  \label{fig:CaseB-(2)-(3)-(b)}
\end{center}\end{figure}

\begin{figure}[htb]\begin{center}
  {\unitlength=1cm
  \begin{picture}(12,5.8)
   \put(-.4,1.4){\includegraphics[keepaspectratio,width=6cm]{%
   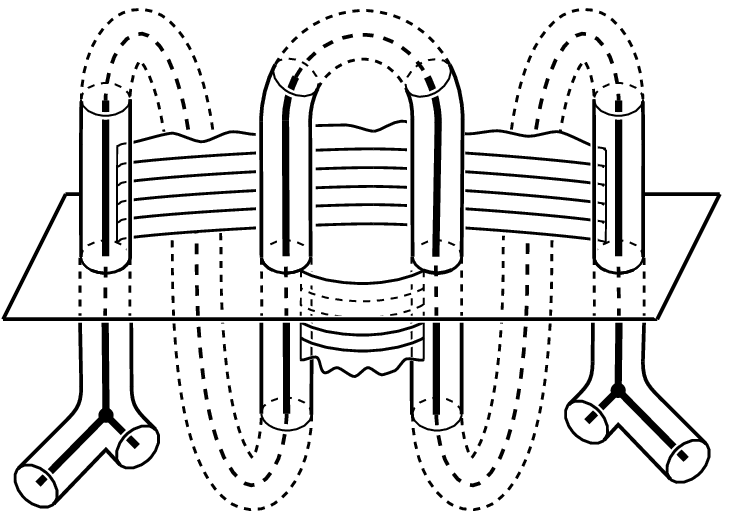}}
   \put(6.4,1.4){\includegraphics[keepaspectratio,width=6cm]{%
   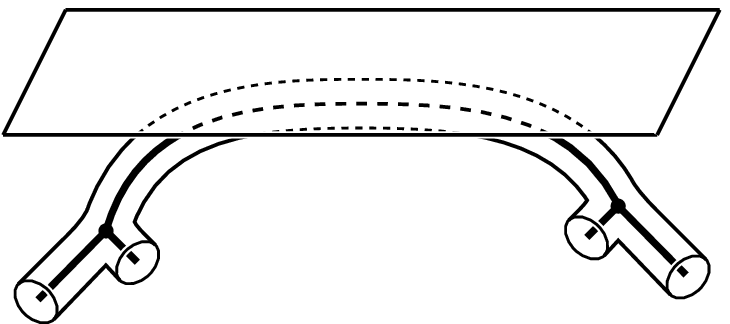}}
   \put(4.5,0.4){\includegraphics[keepaspectratio,height=.6cm]{%
   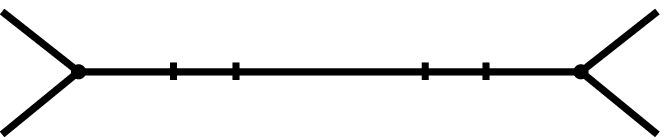}}
   \put(5.7,3.4){$\Longrightarrow$}
   \put(1.5,3.5){$\gamma$}
   \put(2.46,3.4){$\gamma'$}
%
   \put(5.9,1.2){$\bar{\gamma}$}
   \put(5.3,0.85){$\overbrace{\hspace{1.4cm}}$}
%
   \put(4.9,0.6){$\underbrace{\hspace{0cm}}$}
   \put(5.1,0){$\bar{\gamma}'_1$}
   \put(6.5,0.6){$\underbrace{\hspace{0cm}}$}
   \put(6.7,0){$\bar{\gamma}'_2$}
  \end{picture}}
  \caption{}
  \label{fig:CaseB-(2)-(3)-(c)}
\end{center}\end{figure}

\medskip
\noindent\textit{Case B}-(3)-(3).\ \ \ \ 
Both $\bar{\gamma}$ and $\bar{\gamma}'$ satisfy the condition $(3)$. 

\medskip
Let $\bar{\gamma}_1$ and $\bar{\gamma}_2$ ($\bar{\gamma}'_1$ and $\bar{\gamma}'_2$ resp.) 
be the components of $\bar{\gamma}$ ($\bar{\gamma}'$ resp.) with $h^1_{\bar{\gamma}_1}\supset D_w$ 
($h^1_{\bar{\gamma}'_1}\supset D_{w'}$ resp.). 
Without loss of generality, we may suppose that $\bar{\gamma}_1\subset \bar{\gamma}'_1$. 
Then there are teo possibilities: $(1)$ $\bar{\gamma}_2\subset \bar{\gamma}'_2$ 
and $(2)$ $\bar{\gamma}_2\supset \bar{\gamma}'_2$. 
In each case, we can slide $\bar{\gamma}'_1$ into 
$\gamma'$ along the disk $\delta_{\gamma'}$. Moreover, we can isotope $\Sigma$ to 
reduce $(W_{\Sigma},n_{\Sigma})$ is reduced, a contradiction 
(\textit{cf}. Figure \ref{fig:CaseB-(3)-(3)-(a)} and Figure \ref{fig:CaseB-(3)-(3)-(b)}). 

\begin{figure}[htb]\begin{center}
  {\unitlength=1cm
  \begin{picture}(12,5.8)
   \put(-.4,1.4){\includegraphics[keepaspectratio,width=6cm]{%
   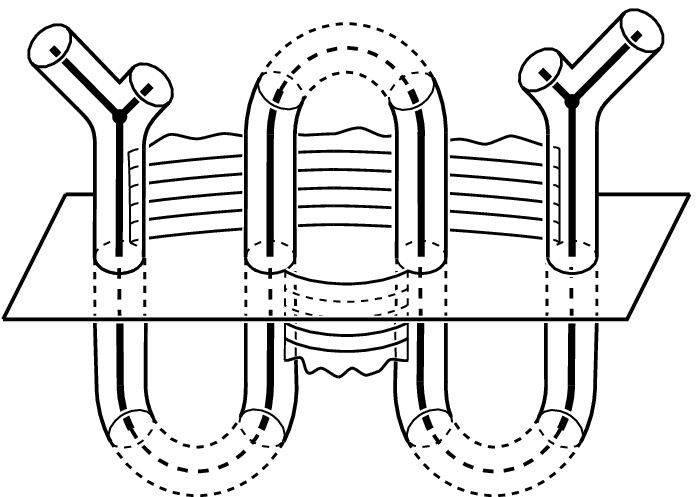}}
   \put(6.4,1.4){\includegraphics[keepaspectratio,width=6cm]{%
   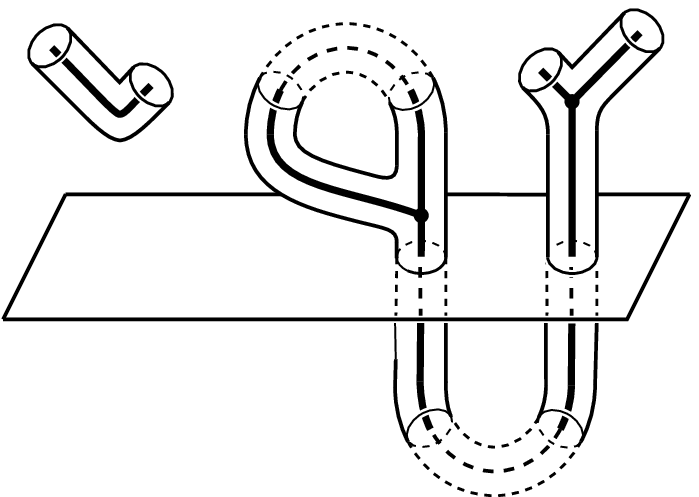}}
   \put(4.5,0.4){\includegraphics[keepaspectratio,height=.6cm]{%
   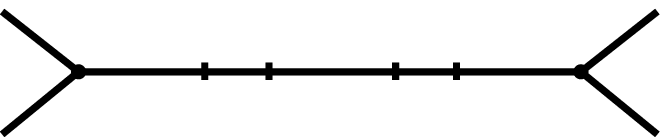}}
   \put(5.7,3.4){$\Longrightarrow$}
   \put(1.1,3.4){$\gamma$}
   \put(2.45,3.35){$\gamma'$}
%
   \put(5.0,1.2){$\bar{\gamma}_1$}
   \put(4.8,0.85){$\overbrace{\hspace{0cm}}$}
   \put(6.8,1.2){$\bar{\gamma}_2$}
   \put(6.6,0.85){$\overbrace{\hspace{0cm}}$}
%
   \put(4.9,0.6){$\underbrace{\hspace{0.8cm}}$}
   \put(5.1,0){$\bar{\gamma}'_1$}
   \put(6.4,0.6){$\underbrace{\hspace{0.8cm}}$}
   \put(6.7,0){$\bar{\gamma}'_2$}
  \end{picture}}
  \caption{}
  \label{fig:CaseB-(3)-(3)-(a)}
\end{center}\end{figure}

\begin{figure}[htb]\begin{center}
  {\unitlength=1cm
  \begin{picture}(12,5.8)
   \put(-.4,1.4){\includegraphics[keepaspectratio,width=6cm]{%
   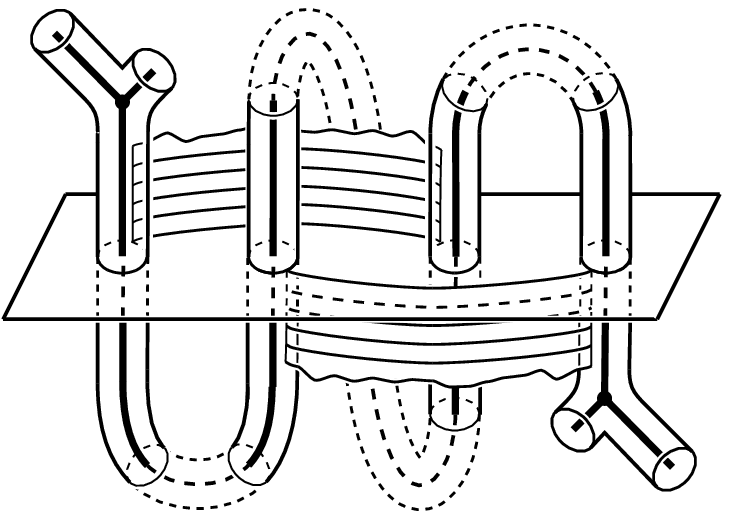}}
   \put(6.4,1.4){\includegraphics[keepaspectratio,width=6cm]{%
   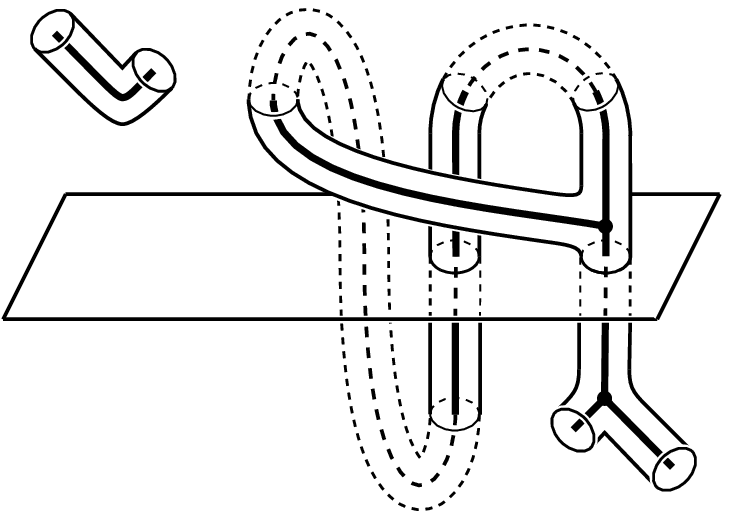}}
   \put(4.5,0.4){\includegraphics[keepaspectratio,height=.6cm]{%
   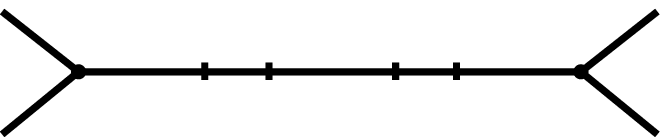}}
   \put(5.7,3.4){$\Longrightarrow$}
   \put(1.1,3.5){$\gamma$}
   \put(3.8,3.4){$\gamma'$}
%
   \put(5.0,1.2){$\bar{\gamma}_1$}
   \put(4.8,0.85){$\overbrace{\hspace{0cm}}$}
   \put(6.7,1.2){$\bar{\gamma}_2$}
   \put(6.4,0.85){$\overbrace{\hspace{0.8cm}}$}
%
   \put(4.9,0.6){$\underbrace{\hspace{0.8cm}}$}
   \put(5.1,0){$\bar{\gamma}'_1$}
   \put(6.6,0.6){$\underbrace{\hspace{0cm}}$}
   \put(6.8,0){$\bar{\gamma}'_2$}
  \end{picture}}
  \caption{}
  \label{fig:CaseB-(3)-(3)-(b)}
\end{center}\end{figure}

\end{proof}

\begin{figure}[htb]\begin{center}
  {\unitlength=1cm
  \begin{picture}(3.6,4.4)
   \put(0.2,0){\includegraphics[keepaspectratio]{%
   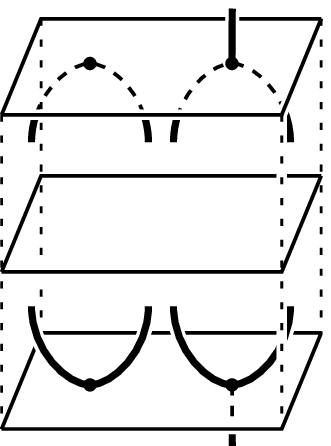}}
   \put(-.3,0){$t_0^-$}
   \put(-.3,1.7){$t_0$}
   \put(-.3,3.2){$t_0^+$}
  \end{picture}}
  \caption{}
  \label{fig:21}
\end{center}\end{figure}

Let $t_0^+$ ($t_0^-$ resp.) be the first critical 
height above $t_0$ (below $t_0$ resp.). Since $|P(t_0)\cap \Sigma|=
W_{\Sigma}=\max \{w_{\Sigma}(t)|t\in (-1,1)\}$, we see that the critical point of the height 
$t_0^+$ ($t_0^-$ resp.) is a maximum or a $\lambda$-vertex (a minimum or a $y$-vertex resp.) 
(see Figure \ref{fig:21}). 

\begin{lem}\label{B4_5}
The critical height $t_0^-$ is a $y$-vertex (not a minimum), or $\Sigma$ is modified 
by edge slides so that the modified graph contains an unknotted cycle. 
\end{lem}

\begin{proof}
Suppose that the critical point of the height $t_0^-$ is a minimum. Let $t_0^{-+}$ be a regular height 
just above $t_0^-$. Then $\Lambda (t_0^{-+})$ contains a fat-vertex with a lower simple outermost edge for 
the fat-vertex of $\Lambda (t_0^{-+})$. Hence it follows from Lemma \ref{B4} 
that every simple outermost edge for each fat-vertex of $\Lambda (t_0^{-+})$ is lower. Similarly, 
every simple outermost edge for each fat-vertex of $\Lambda (t_0)$ is upper. We now vary $t$ for $t_0^{-+}$ to $t_0$. 
Note that for each regular height $t$, all the simple outermost edges for each fat-vertex of $\Lambda (t)$ 
are either upper or lower (Lemma \ref{B4}); such a regular height $t$ is said to be upper or lower 
respectively. In these words, $t_0^{-+}$ is lower and $t_0$ is upper. 

Let $c_1,\dots,c_n$ $(c_1<\dots <c_n)$ be the critical heights of $h|_{\Delta_2}$ contained in 
$[t_0^{-+},t_0]$. Note that the property `upper' or `lower' is unchanged at any height of 
$[t_0^{-+},t_0]\setminus \{ c_1,\dots ,c_n\}$. Hence there exists a critical height $c_i$ such that a height 
$t$ is changed from lower to upper at $c_i$. 
The graph $\Lambda(t)$ is changed as in Figure \ref{fig:B04} around the critical height $c_i$. 

\begin{figure}[htb]\begin{center}
  {\unitlength=1cm
  \begin{picture}(6,5.4)
   \put(.9,.3){\includegraphics[keepaspectratio]{%
             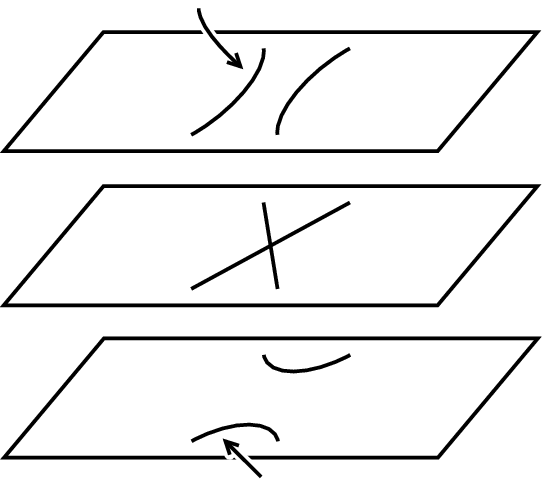}}
   \put(2.3,5.3){upper}
   \put(3.6,0){lower}
   \put(2.5,2.2){$w$}
   \put(-.3,4.3){$t=c_i^+$}
   \put(-.3,2.7){$t=c_i$}
   \put(-.3,1.1){$t=c_i^-$}
  \end{picture}}
  \caption{}
  \label{fig:B04}
\end{center}\end{figure}

Let $c_i^+$ ($c_i^-$ resp.) be a regular height just above (below resp.) $c_i$. We note that the lower disk 
for $\Lambda (c_i^-)$ and the upper disk for $\Lambda (c_i^+)$ in Figure \ref{fig:B04} are contained in the 
same component of $\Delta_2$, say $D$. We take parallel copies, say $D'$ and $D''$, of $D$ such that 
$D'$ is obtained by pushing $D$ into one side and that $D''$ is obtained by pushing $D$ into the other side 
(\textit{cf}. Figure \ref{fig:B05}). 
Then we may suppose that there is an upper (a lower resp.) simple outermost edge for a fat-vertex 
in $D'$ ($D''$ resp.). Hence we can apply the arguments of the proof of Lemma \ref{B3} 
to modify $\Sigma$ so that the modified graph contains an unknotted cycle. 

\begin{figure}[htb]\begin{center}
  {\unitlength=1cm
  \begin{picture}(4.3,3.1)
   \put(0.4,0.5){\includegraphics[keepaspectratio]{%
   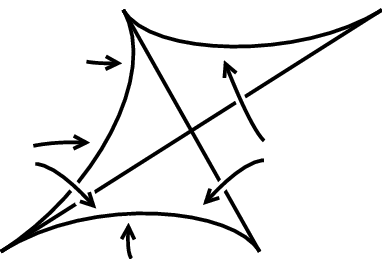}}
   \put(0.2,2.5){upper}
   \put(0.3,1.5){$\sigma''$}
   \put(0.0,.4){$w$}
   \put(1.5,0.1){lower}
   \put(3.2,1.5){$\sigma'$}
  \end{picture}}
  \caption{}
  \label{fig:B05}
\end{center}\end{figure}

\end{proof}

Let $v^-$ be the $y$-vertex of $\Sigma$ at the height $t_0^-$ and 
$t_0^{--}$ a regular height just below $t_0^-$. Let $v^{--}$ be the intersection point 
of the descending edges from $v^-$ in $\Sigma$ and $P(t_0^{--})$, 
and let $D_{v^{--}}$ be the fat-vertex of $\Lambda(t_0^{--})$ correponding to $v^{--}$. 

\begin{lem}\label{B5}
Every simple outermost edge for any fat-vertex of $\Lambda(t_0^{--})$ 
is lower, or $\Sigma$ is modified by edge slides so that the modified graph contains an unknotted cycle. 
\end{lem}

\begin{proof}
Suppose that there is a fat-vertex $D_w$ of $\Lambda(t_0^{--})$ such that $\Lambda(t_0^{--})$ contains 
an upper simple outermost edge $\gamma$ for $D_w$. 
Let $\sigma$ be the edge of $\Sigma$ with $h^1_{\sigma}\supset D_w$. 
Let $\delta_{\gamma}$ be the outermost disk for $(D_w,\gamma)$. 
Let $\bar{\gamma}$ ($\bar{\gamma}'$ resp.) be a union of the components obtained by cutting $\sigma$ 
by the two fat-vertices of $\Lambda(t_0)$ incident to $\gamma$ ($\gamma'$ resp.) 
such that a $1$-handle correponding to each component intersects 
$\partial \delta_{\gamma}\setminus \gamma$ 
($\partial \delta_{\gamma'}\setminus \gamma'$ resp.). 
Then $\bar{\gamma}$ ($\bar{\gamma}'$ resp.) satisfies one of the conditions $(1)$, $(2)$ and $(3)$ 
in the proof of Lemma \ref{B3}. 

\medskip
\noindent\textit{Case A.}\ \ \ \ $D_{v^{--}}\ne D_w$. 

\medskip
Then we have the following three cases. In each case, we can slide (a component of) $\bar{\gamma}$ to 
$\gamma$ along the disk $\delta_{\gamma}$. Moreover, we can isotope $\Sigma$ to reduce 
$(W_{\Sigma},n_{\Sigma})$, a contradiction. 

\medskip
\noindent\textit{Case A}-(1). \ \ \ \ 
$\bar{\gamma}$ satisfies the condition $(1)$. 

\medskip
Then there are two possibilities: (i) $v^{--}\not\in \bar{\gamma}$ and (ii) $v^{--}\in \bar{\gamma}$. 
In each case, see Figure \ref{fig:C01}. 

\begin{figure}[htb]\begin{center}
  {\unitlength=1cm
  \begin{picture}(12,11)
   \put(0.2,6.3){\includegraphics[keepaspectratio,width=5cm]{%
             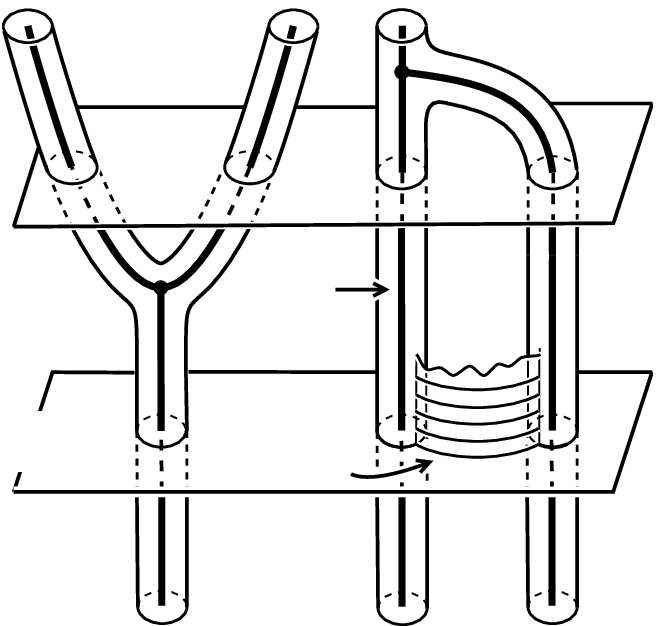}}
   \put(7,6.3){\includegraphics[keepaspectratio,width=5cm]{%
             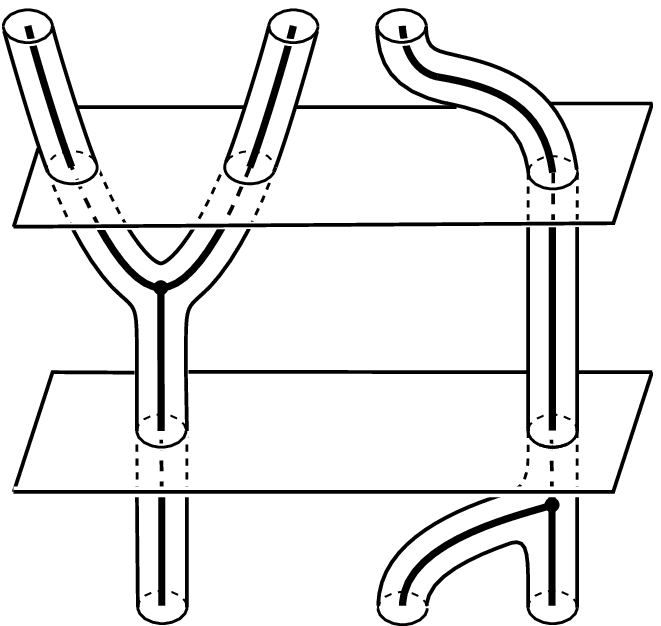}}
   \put(0.2,0.3){\includegraphics[keepaspectratio,width=5cm]{%
             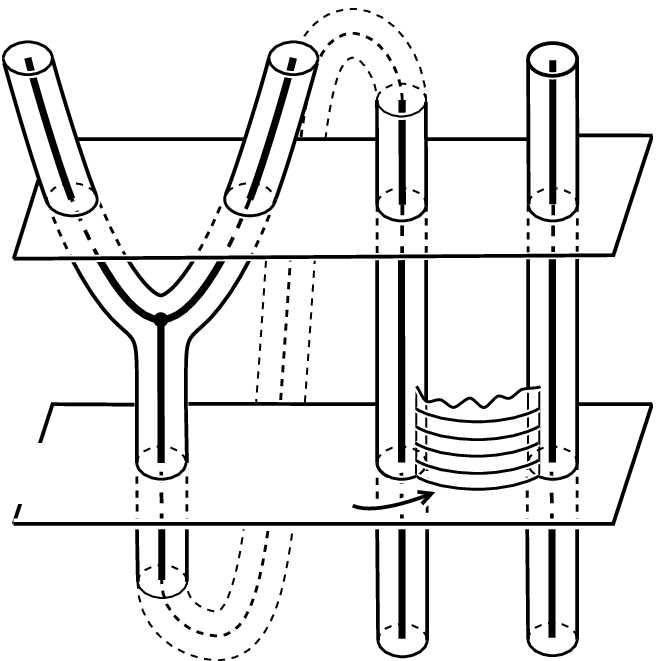}}
   \put(7,0.3){\includegraphics[keepaspectratio,width=5cm]{%
             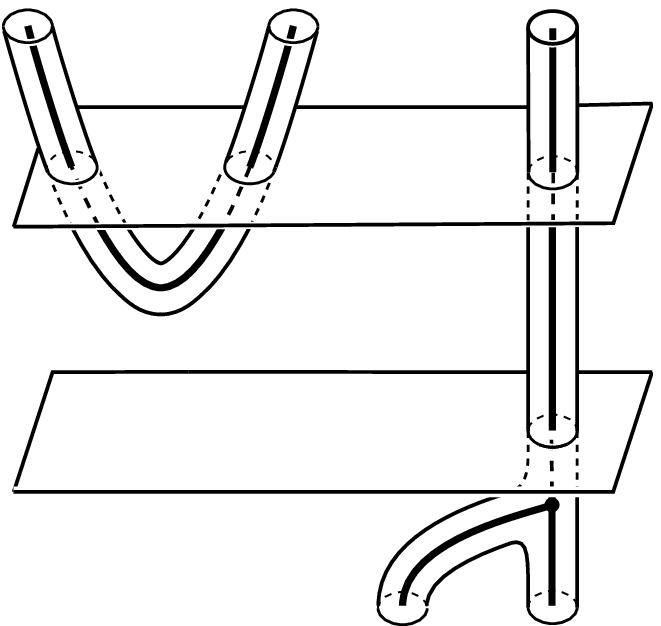}}
   \put(5,6){(i)~~$v^{--}\not\in \bar{\gamma}$}
   \put(5,0){(ii)~~$v^{--}\in \bar{\gamma}$}
   \put(-.2,9){$t=t_0$}
   \put(-.2,6.9){$t=t_0^{--}$}
   \put(-.2,3){$t=t_0$}
   \put(-.2,.9){$t=t_0^{--}$}
   \put(5.9,8.4){$\Longrightarrow$}
   \put(5.9,2.4){$\Longrightarrow$}
   \put(0.4,7.6){$D_{v^{--}}$}
   \put(2.5,7.8){$D_w$}
   \put(2.56,7.4){$\gamma$}
   \put(2.5,8.8){$\bar{\gamma}$}
   \put(0.4,1.6){$D_{v^{--}}$}
   \put(2.5,1.9){$D_w$}
   \put(2.56,1.5){$\gamma$}
  \end{picture}}
  \caption{}
  \label{fig:C01}
\end{center}\end{figure}

\medskip
\noindent\textit{Case A}-(2). \ \ \ \ 
$\bar{\gamma}$ satisfies the condition $(2)$. 

\medskip
See Figure \ref{fig:C02}. 

\begin{figure}[htb]\begin{center}
  {\unitlength=1cm
  \begin{picture}(12,4.7)
   \put(0.2,0.0){\includegraphics[keepaspectratio,width=5cm]{%
             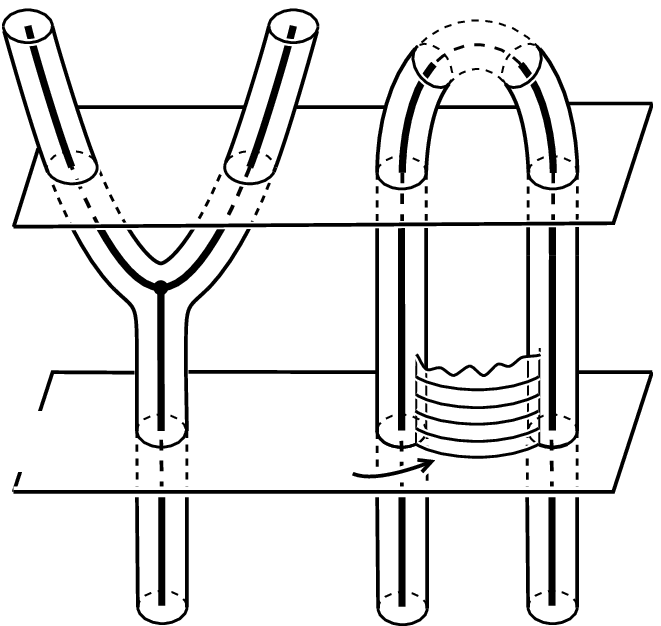}}
   \put(7,0.0){\includegraphics[keepaspectratio,width=5cm]{%
             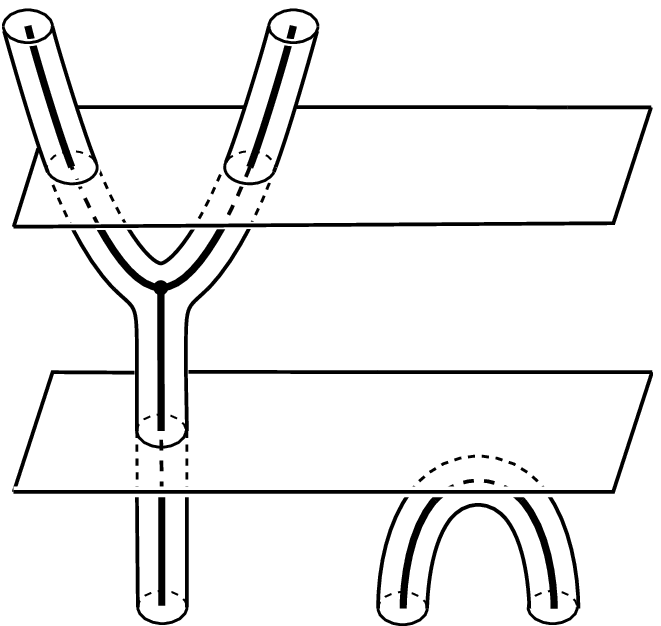}}
   \put(-.2,2.7){$t=t_0$}
   \put(-.2,.6){$t=t_0^{--}$}
   \put(5.9,2.4){$\Longrightarrow$}
   \put(0.4,1.3){$D_{v^{--}}$}
   \put(2.5,1.6){$D_w$}
   \put(2.6,1.2){$\gamma$}
  \end{picture}}
  \caption{}
  \label{fig:C02}
\end{center}\end{figure}

\medskip
\noindent\textit{Case A}-(3). \ \ \ \ 
$\bar{\gamma}$ satisfies the condition $(3)$. 

\medskip
Then there are two possibilities: (i) $v^{--}\not\in \bar{\gamma}$ and (ii) $v^{--}\in \bar{\gamma}$. 
In each case, see Figure \ref{fig:C03}. 

\begin{figure}[htb]\begin{center}
  {\unitlength=1cm
  \begin{picture}(12,4.7)
   \put(0.2,0.0){\includegraphics[keepaspectratio,width=5cm]{%
             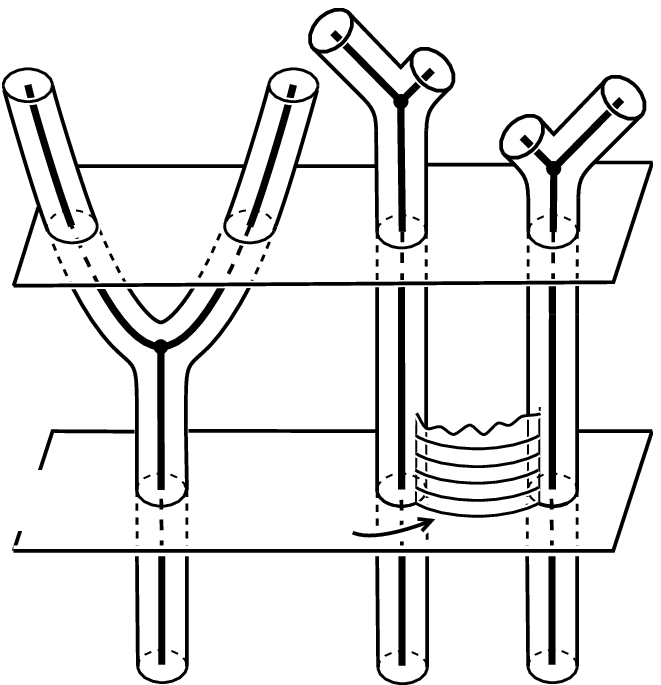}}
   \put(7,0.0){\includegraphics[keepaspectratio,width=5cm]{%
             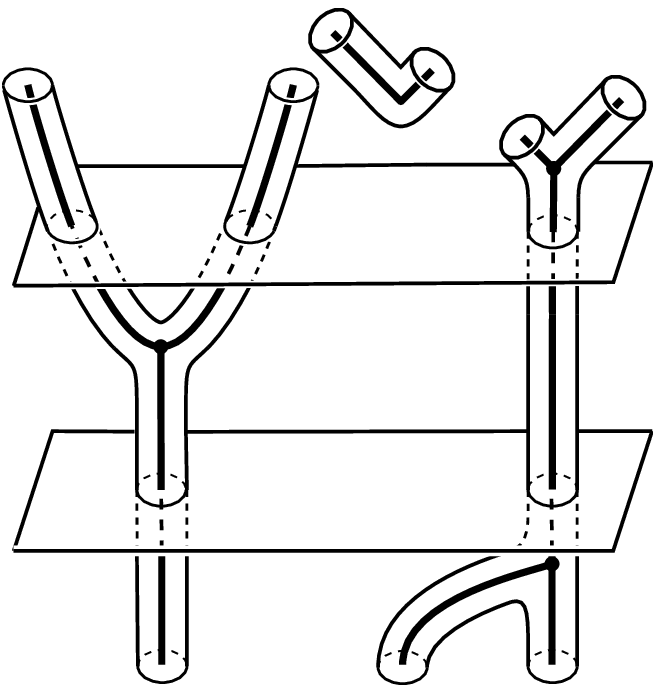}}
   \put(-.2,2.7){$t=t_0$}
   \put(-.2,.6){$t=t_0^{--}$}
   \put(5.9,2.4){$\Longrightarrow$}
   \put(0.4,1.3){$D_{v^{--}}$}
   \put(2.5,1.6){$D_w$}
   \put(2.6,1.2){$\gamma$}
  \end{picture}}
  \caption{}
  \label{fig:C03}
\end{center}\end{figure}

\medskip
\noindent\textit{Case B.}\ \ \ \ $D_{v^{--}}= D_{w}$. 

\medskip
Since $\delta_{\gamma}$ is upper, we see that $\bar{\gamma}$ does not satisfy the condition $(2)$. 
Hence we have the following. 

\medskip
\noindent\textit{Case B}-(1). \ \ \ \ 
$\bar{\gamma}$ satisfies the condition $(1)$. 

\medskip
Since $\gamma$ is upper, we see that the $y$-vertex of $\Sigma$ at the height $t_0^-$ 
is an endpoint of $\bar{\gamma}$, i.e., 
$\bar{\gamma}$ is the short vertical arc joining $v^-$ to $v^{--}$. Then we can slide $\bar{\gamma}$ 
to $\gamma$ along the disk $\delta_{\gamma}$ to obtain a new graph $\Sigma'$. Note that 
$(W_{\Sigma'},n_{\Sigma'})=(W_{\Sigma},n_{\Sigma})$ (\textit{cf}. Figure \ref{fig:D02}). 
However, the critical point for $\Sigma'$ corresponding to $v^-$ is a minimum. Hence we can apply 
the arguments in the proof of Lemma \ref{B4_5} to show that there is an unknotted cycle in $\Sigma'$. 

\begin{figure}[htb]\begin{center}
  {\unitlength=1cm
  \begin{picture}(10,5)
   \put(0.2,0.0){\includegraphics[keepaspectratio,width=4cm]{%
             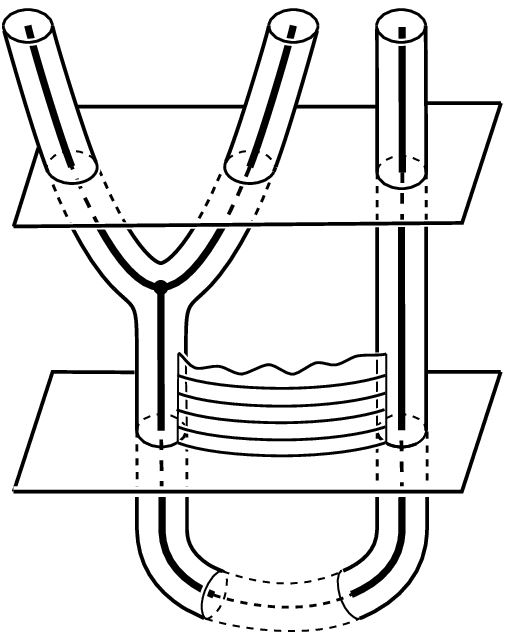}}
   \put(6,0.0){\includegraphics[keepaspectratio,width=4cm]{%
             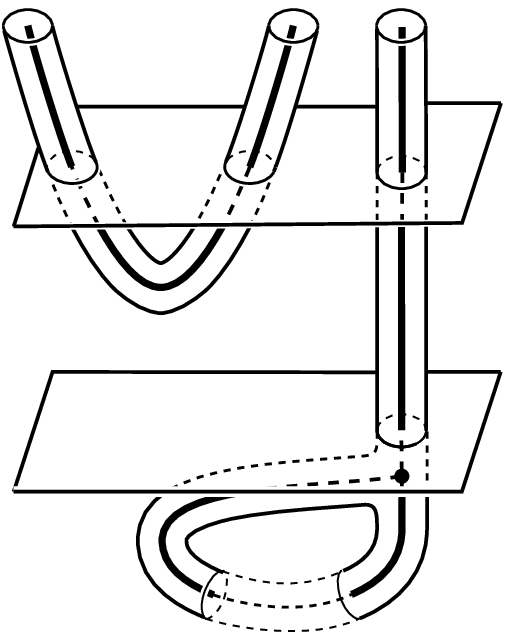}}
   \put(-.2,2.9){$t=t_0$}
   \put(-.2,.8){$t=t_0^{--}$}
   \put(4.9,2.4){$\Longrightarrow$}
   \put(2.7,1.2){$\gamma$}
  \end{picture}}
  \caption{}
  \label{fig:D02}
\end{center}\end{figure}

\medskip
\noindent\textit{Case B}-(3). \ \ \ \ 
$\bar{\gamma}$ satisfies the condition $(3)$. 

\medskip
Let $\bar{\gamma}_1$ and $\bar{\gamma}_2$ be the components of $\bar{\gamma}$ with 
$\partial \bar{\gamma}_1\ni v^-$. Then we can slide $\bar{\gamma}_2$ into $\gamma$ along the disk 
$\delta_{\gamma}$. Moreover, we can isotope $\Sigma$ to reduce $(W_{\Sigma},n_{\Sigma})$, 
a contradiction (\textit{cf}. Figure \ref{fig:D03}). 

\begin{figure}[htb]\begin{center}
  {\unitlength=1cm
  \begin{picture}(10,5.5)
   \put(0.2,0.0){\includegraphics[keepaspectratio,width=4cm]{%
             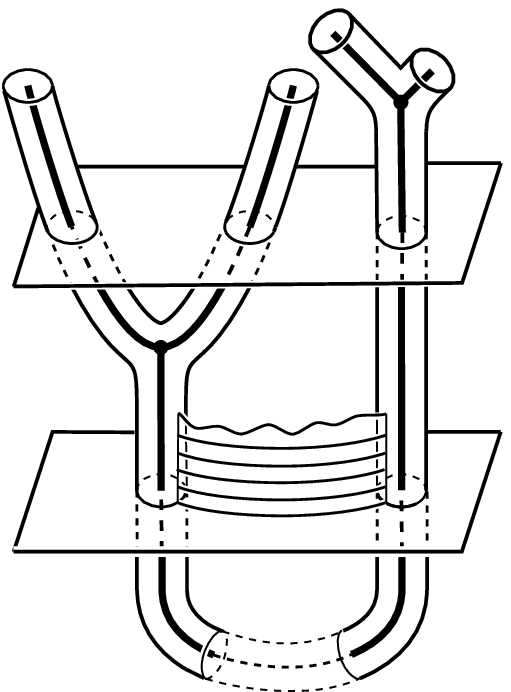}}
   \put(6,0.0){\includegraphics[keepaspectratio,width=4cm]{%
             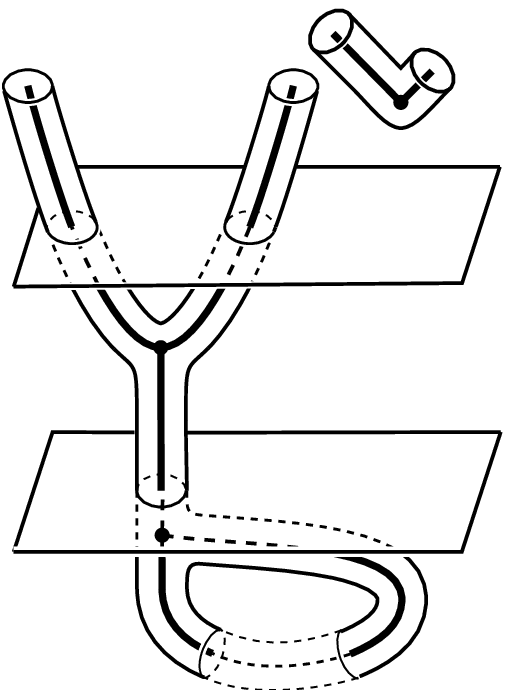}}
   \put(-.2,2.9){$t=t_0$}
   \put(-.2,.8){$t=t_0^{--}$}
   \put(4.9,2.4){$\Longrightarrow$}
   \put(2.7,1.2){$\gamma$}
  \end{picture}}
  \caption{}
  \label{fig:D03}
\end{center}\end{figure}

\end{proof}

\begin{lem}\label{B6}
Every simple outermost edge for any fat-vertex of $\Lambda(t_0^{--})$ is incident to $D_{v^{--}}$, 
or $\Sigma$ is modified so that there is an unknotted cycle. 
\end{lem}

\begin{proof}
Suppose that $\Lambda(t_0^{--})$ contains a simple outermost edge $\gamma$ for $D_w$ 
and is not incident to $D_{v^{--}}$. Then it follows from Lemma \ref{B5} that $\gamma$ is lower. 
This means that $\Lambda(t_0)$ contains a lower simple edge, because an edge disjoint from 
$D_{v^{--}}$ is not affected at all in $[t_0^{--},t_0]$. 
This contradicts Lemma \ref{B4}. 

\end{proof}

We now prove Proposition \ref{B}. 

\begin{proof}[Proof of Proposition \ref{B}]
We first prove the following. 

\medskip
\noindent\textit{Claim.}\ \ \ \ 
For any fat-vertex $D_w(\ne D_{v^{--}})$ of $\Lambda(t_0^{--})$, there are no loops of $\Lambda(t_0^{--})$ 
based on $D_w$, or $\Sigma$ is modified by edge slides so that the modified graph contains an unknotted cycle. 

\medskip
\noindent\textit{Proof.}\ \ \ \ 
Suppose that there is a fat-vertex $D_w(\ne D_{v^{--}})$ of $\Lambda(t_0^{--})$ such that there is a loop 
$\alpha$ of $\Lambda(t_0^{--})$ based on $D_w$. Then $\alpha$ separates $\mathrm{cl}(P(t_0^{--})\setminus D_w)$ 
into two disks $E_1$ and $E_2$ with $D_{v^{--}}\subset E_2$. By retaking $D_{w}$ and $\alpha$, if necessary, 
we may suppose that there are no loop components of $\Lambda(t_0^{--})$ in $\mathrm{int}(E_1)$. 
It follows from Lemma \ref{B1} that there is a fat-vertex $D_{w'}$ of $\Lambda(t_0^{--})$ in $\mathrm{int}(E_1)$. 
Then every outermost edge for $D_{w'}$ of $\Lambda(t_0^{--})$ is simple. 
Hence it follows from Lemma \ref{B6} that $\Sigma$ contains an unknotted cycle and therefore we have the claim. 

\medskip
Then we have the following cases. 

\medskip
\noindent\textit{Case A.}\ \ \ \ 
The descending edges of $\Sigma$ from the maximum or $\lambda$-vertex $v^+$ at the height 
$t_0^+$ are equal to the ascending edges from $v^-$ (\textit{cf}. Figure \ref{fig:28}). 

\begin{figure}[htb]\begin{center}
  {\unitlength=1cm
  \begin{picture}(4.3,4.6)
   \put(0,0){\includegraphics[keepaspectratio]{%
   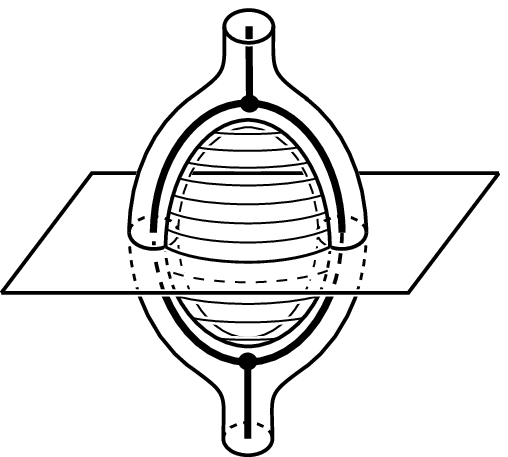}}
  \end{picture}}
  \caption{}
  \label{fig:28}
\end{center}\end{figure}

\medskip
Then we can immediately see that there is an unknotted cycle . 

\medskip
\noindent\textit{Case B.}\ \ \ \ 
Exactly one of the descending edges from $v^+$ is equal to one of the ascending edges from 
$v^-$ (\textit{cf}. Figure \ref{fig:29}). 

\begin{figure}[htb]\begin{center}
  {\unitlength=1cm
  \begin{picture}(4.3,7.1)
   \put(0,0){\includegraphics[keepaspectratio]{%
   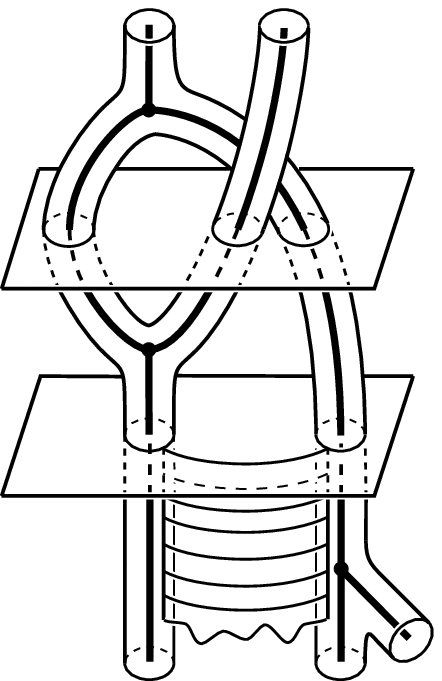}}
  \end{picture}}
  \caption{}
  \label{fig:29}
\end{center}\end{figure}

\medskip
Let $\sigma'$ be the other edge disjoint from $v^-$, and let $w^{--}$ be the first intersection point 
of $P(t_0^{--})$ and the edge $\sigma'$. 
Let $\gamma$ be an outermost edge for $D_{w^{--}}$ of $\Lambda(t_0^{--})$. 
By the claim above, we see that $\gamma$ is simple. 
It follows from Lemma \ref{B5} that we may suppose that $\gamma$ is lower. It also follows from 
Lemma \ref{B6} that we may suppose that the endpoints of $\gamma$ are $v^{--}$ and $w^{--}$. 
Let $\delta_{\gamma}$ be the outermost disk for $(D_w,\gamma)$. Set $\bar{\gamma}=\sigma'\cap \delta_{\gamma}$. 
Since the subarc of $\sigma'$ whose endpoints are $v^{--}$ and $w^{--}$ is monotonous and 
$\gamma$ is lower, we see that $\bar{\gamma}$ cannot satisfy the condition $(3)$ in the proof 
of Lemma \ref{B4}. Hence $\bar{\gamma}$ satisfies the condition $(1)$ or $(2)$. In each case, 
we can slide $\bar{\gamma}$ to $\gamma$ along the disk $\delta_{\gamma}$ to obtain a new graph 
with an unknotted cycle. 

\medskip
\noindent\textit{Case C.}\ \ \ \ 
Any descending edge of $\Sigma$ from $v^+$ is disjoint from an ascending edge from $v^-$. 

\medskip
It follows from \ref{B5}, \ref{B6} and the claim that $\Lambda(t^{--})$ 
contains a lower simple outermost edge $\gamma_i$ $(i=1,2)$ 
for $D_{w_i}$ which is adjacent to $D_{w_i}$ and $D_{v^{--}}$. Let 
$\delta_{\gamma_i}$ be the outermost disk for $(D_{w_i},\gamma_i)$. 
Set $\bar{\gamma}_i=\sigma_i\cap \delta_{\gamma_i}$. Since the subarc of $\sigma_i$ whose endpoints 
$v^+$ and $w_i$ are monotonous and $\delta_{\gamma_i}$ is lower, we see that $\gamma_i$ cannot 
satisfy the condition $(3)$. Then we have the following. 

\medskip
\noindent\textit{Case C}-(1). \ \ \ \ 
Both $\bar{\gamma}_1$ and $\bar{\gamma}_2$ satisfy the condition $(1)$. 

\medskip
If $\sigma_1=\sigma_2$, then we can slide $\bar{\gamma}_1$ to $\gamma_1$ along the disk 
$\delta_{\gamma_1}$. We can further isotope $\Sigma$ to reduce $(W_{\Sigma},n_{\Sigma})$, a contradiction 
(\textit{cf}. Figure \ref{fig:73}). Hence $\sigma_1\ne \sigma_2$. 

\begin{figure}[htb]\begin{center}
  {\unitlength=1cm
  \begin{picture}(9.6,5.0)
   \put(0,-.2){\includegraphics[keepaspectratio,width=3.6cm]{%
   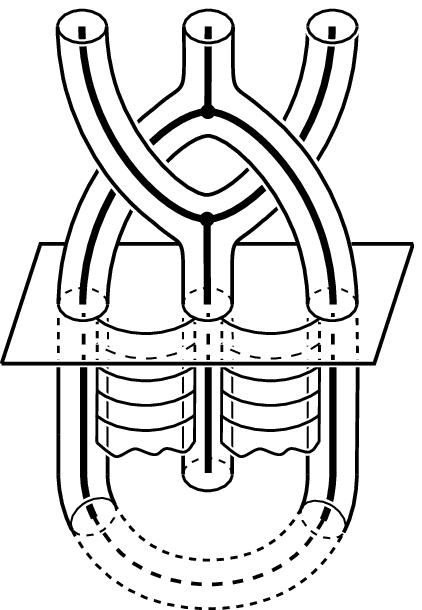}}
   \put(6.0,.8){\includegraphics[keepaspectratio,width=3.6cm]{%
   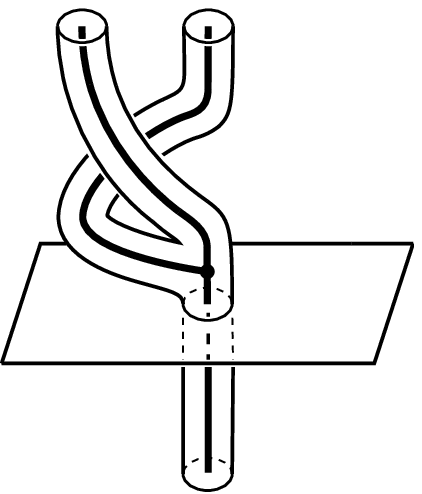}}
   \put(4.5,2.3){$\Longrightarrow$}
  \end{picture}}
  \caption{}
  \label{fig:73}
\end{center}\end{figure}

Then we can slide $\bar{\gamma}_1\cup \bar{\gamma}_2$ to $\gamma_1\cup \gamma_2$ along 
$\delta_{\gamma_1}\cup \delta_{\gamma_2}$ so that a new graph contains an unknotted cycle 
(\textit{cf}. Figure \ref{fig:30}). 

\begin{figure}[htb]\begin{center}
  {\unitlength=1cm
  \begin{picture}(4.2,5.7)
   \put(0,0){\includegraphics[keepaspectratio]{%
   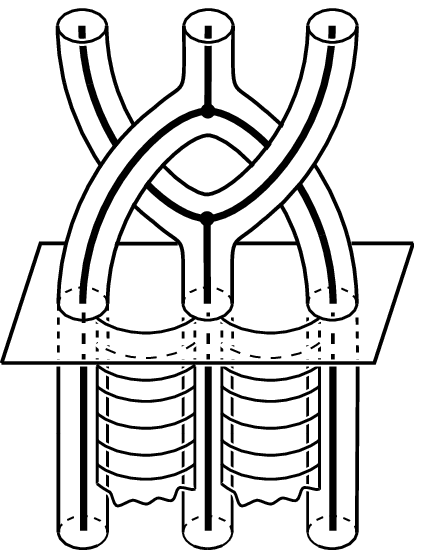}}
  \end{picture}}
  \caption{}
  \label{fig:30}
\end{center}\end{figure}

\medskip
\noindent\textit{Case C}-(2). \ \ \ \ 
Either $\bar{\gamma}_1$ or $\bar{\gamma}_2$, say $\bar{\gamma}_1$, satisfies the condition $(2)$. 

\medskip
Since $\bar{\gamma}_1$ satisfies the condition $(2)$, we see that the endpoints of $\sigma_1$ are 
$v^+$ and $v^-$. Hence $w_2\not\in \sigma_1$. This implies that $\bar{\gamma}_2$ satisfies the condition 
$(1)$. Then we first slide $\bar{\gamma}_2$ to $\gamma_2$ along $\delta_{\gamma_2}$. 
We can further slide $\bar{\gamma}_1$ to $\gamma_1$ along $\delta_{\gamma_1}$ so that 
a new graph contains an unknotted cycle. 

This completes the proof of Proposition \ref{B}. 
\end{proof}

\subsection{Applications of Haken's theorem and Waldhausen's theorem}

\begin{cor}\label{S}
Let $M$ be a compact $3$-manifold and $(C_1,C_2;S)$ a reducible Heegaard splitting. 
Then $M$ is reducible or $(C_1,C_2;S)$ is stabilized. 
\end{cor}

\begin{proof}
Suppose that $M$ is irreducible. Let $P$ be a $2$-sphere such that $P\cap S$ is an 
essential loop. Since $M$ is irreducible, we see that $P$ bounds a $3$-ball in $M$. Hence we can regard 
$M$ as a connected sum of $S^3$ and $M$. By Theorem \ref{S6}, the induced 
Heegaard splitting of $S^3$ is stabilized. Hence this cancelling pair of disks shows that 
$(C_1,C_2;S)$ is stabilized. 
\end{proof}

\begin{cor}\label{S7}
Any Heegaard splitting of a handlebody is standard, i.e, is obtained from 
a trivial splitting by stabilization. 
\end{cor}

\begin{exercise}
\rm
Show Corollary \ref{S7}. 
\end{exercise}

\begin{thm}\label{S8}
Let $M$ be a closed $3$-manifold. Let $(C_1,C_2;S)$ and $(C'_1,C'_2;S')$ be Heegaard 
splittings of $M$. Then there is a Heegaard splitting which is obtained by stabilization 
of both $(C_1,C_2;S)$ and $(C'_1,C'_2;S')$. 
\end{thm}

\begin{proof}
Let $\Sigma_{C_1}$ and $\Sigma_{C_1'}$ be spines of $C_1$ and $C'_1$ respectively. By 
an isotopy, we may assume that $\Sigma_{C_1}\cap \Sigma_{C'_1}=\emptyset$ and 
$C_1\cap C'_1=\emptyset$. 
Set $M'=\mathrm{cl}(M\setminus (C_1\cup C'_1))$, $\partial_1 M'=\partial C_1$ and 
$\partial_2 M'=\partial C_2$. Let $(\bar{C}_1,\bar{C}_2;\bar{S})$ be a Heegaard splitting of 
$(M';\partial_1 M', \partial_2 M')$. Set $C^{\ast}_1=C_1\cup \bar{C}_1$ and 
$C^{\ast}_2=C_2\cup \bar{C}_2$. Then it is easy to see that $(C^{\ast}_1,C^{\ast}_2;\bar{S})$ 
is a Heegaard splitting of $M$. Note that 
$C_2'=C_1\cup M'=C_1\cup (\bar{C}_1\cup \bar{C}_2)=(C_1\cup \bar{C}_1)\cup \bar{C}_2$. 
Here, we note that $(C^{\ast}_1,\bar{C}_2;\bar{S})$ is a Heegaard splitting of 
$C_2'$. It follows from Corollary \ref{S7} that 
$(C^{\ast}_1,\bar{C}_2;\bar{S})$ is obtained from a trivial splitting of $C_2'$ by stabilization. 
This implies that $(C^{\ast}_1,C^{\ast}_2;\bar{S})$ is obtained from $(C'_1,C'_2;S')$ 
by stabilization. On the argument above, by replacing $C_1$ to $C_1'$, we see that 
$(C^{\ast}_1,C^{\ast}_2;\bar{S})$ is also obtained from $(C_1,C_2;S)$ by stabilization. 
\end{proof}

\begin{rem}
\rm
The stabilization problem is one of the most important themes on Heegaard theory. But we do not 
give any more here. 
For the detail, for example, see \cite{Lei}, \cite{Reidemeister}, \cite{Rubinstein}, 
\cite{Sedgwick} and \cite{Singer}. 
\end{rem}

\section{Generalized Heegaard splittings}

\subsection{Definitions}

\begin{defi}
\rm
A \textit{$0$-fork} is a connected $1$-complex obtained 
by joining a point $p$ to a point $g$ whose $1$-simplices are oriented toward $g$ and away from $p$. 
For $n\ge 1$, an \textit{$n$-fork} is a connected $1$-complex obtained 
by joining a point $p$ to each of distinct $n$ points $t_i$  $(i=1,...,n)$ and to a point $g$ 
whose $1$-simplices are oriented toward $g$ and away from $t_i$. 
We call $p$ a \textit{root}, $t_i$ a \textit{tine} and $g$ a \textit{grip}. 
\end{defi}

\begin{rem}
\rm
An $n$-fork corresponds to a compression body $C$ such that each of $t_i$ $(i=1,2, ..., n)$ corresponds to 
a component of $\partial_- C$ and $g$ correponds to $\partial_+ C$ (\textit{cf}. Figure \ref{fig:31}). 

\begin{figure}[htb]\begin{center}
  {\unitlength=1cm
  \begin{picture}(8.4,1.6)
   \put(.7,0){\includegraphics[keepaspectratio]{%
   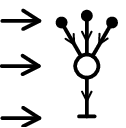}}
   \put(3.3,0){\includegraphics[keepaspectratio]{%
   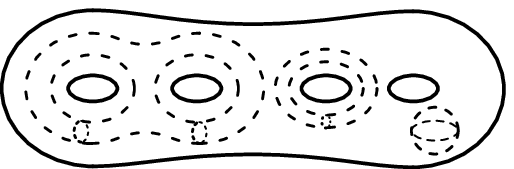}}
   \put(0,1.0){tine}
   \put(0,.54){root}
   \put(0,0.0){grip}
  \end{picture}}
  \caption{}
  \label{fig:31}
\end{center}\end{figure}

\end{rem}

\begin{defi}
\rm
Let $\mathcal{A}$ ($\mathcal{B}$ resp.) be a collection of finite forks, $T_{\mathcal{A}}$ 
($T_{\mathcal{B}}$ resp.) a collection of tines of $\mathcal{A}$ ($\mathcal{B}$ resp.) 
and $G_{\mathcal{A}}$ ($G_{\mathcal{B}}$ resp.) a collection of grips of $\mathcal{A}$ 
($\mathcal{B}$ resp.). We suppose that there are bijections $\mathcal{T}:T_{\mathcal{A}} 
\to T_{\mathcal{B}}$ and $\mathcal{G}:G_{\mathcal{A}} \to G_{\mathcal{B}}$. 
A \textit{fork complex} $\mathscr{F}$ is an oriented connected $1$-complex $\mathcal{A}\cup 
(-\mathcal{B})/\{\mathcal{T},\mathcal{G}\}$, where $-\mathcal{B}$ denotes the $1$-complex obtained 
by taking the opposite orientation of each $1$-simplex and the equivalence relation $/\{\mathcal{T},
\mathcal{G}\}$ is given by $t\sim \mathcal{T}(t)$ for any $t\in T_{\mathcal{A}}$ 
and $g\sim \mathcal{G}(g)$ for any $g\in G_{\mathcal{A}}$. We define: 

\bigskip
\begin{center}
$\partial_1 \mathscr{F}=\{($tines of $\mathcal{A})\setminus T_{\mathcal{A}}\}\cup \{($grips of $\mathcal{B})
\setminus G_{\mathcal{B}}\}$ and 
\end{center}

\begin{center}
$\partial_2 \mathscr{F}=\{($tines of $\mathcal{B})\setminus T_{\mathcal{B}}\}
\cup \{($grips of $\mathcal{A})\setminus G_{\mathcal{A}}\}$. 
\end{center}
\end{defi}

\begin{defi}
\rm
A fork complex is $\textit{exact}$ if there exists $e\in \mathrm{Hom} (C_0(\mathscr{F}),\mathbb{R})$ such that 
\begin{enumerate}
\item $e(v_1)=0$ for any $v_1\in \partial_1 \mathscr{F}$, 
\item $(\delta e)(e_\mathcal{A})>0$ for any $1$-simplex $e_\mathcal{A}$ in $\mathcal{A}$ with the standard 
orientation, $(\delta e)(e_\mathcal{A})<0$ for any $1$-simplex $e_\mathcal{B}$ in $\mathcal{B}$ with the 
standard orientation, where $\delta$ denotes the coboundary operator 
$\mathrm{Hom} (C_0(\mathscr{F}),\mathbb{R})\to \mathrm{Hom} (C_1(\mathscr{F}),\mathbb{R})$ and 
\item $e(v_2)=1$ for any $v_2\in \partial_2 \mathscr{F}$. 
\end{enumerate}
\end{defi}

\begin{rem}\rm
Geometrically speaking, $\mathscr{F}$ is exact if and only if we can put $\mathscr{F}$ in $\mathbb{R}^3$ 
so that 
\begin{enumerate}
\item $\partial_1 \mathscr{F}$ lies in the plane of height $0$, 
\item for any path $\alpha$ in $\mathscr{F}$ from a point in $\partial_1 \mathscr{F}$ to a point 
in $\partial_2 \mathscr{F}$, $h|_{\alpha}$ is monotonically increasing, where $h$ is the height 
function of $\mathbb{R}^3$ and 
\item $\partial_2 \mathscr{F}$ lies in the plane of height $1$ (\textit{cf}. Figure \ref{fig:31.2}). 
\end{enumerate}

\begin{figure}[htb]\begin{center}
  {\unitlength=1cm
  \begin{picture}(4,3.5)
   \put(0,0){\includegraphics[keepaspectratio]{%
   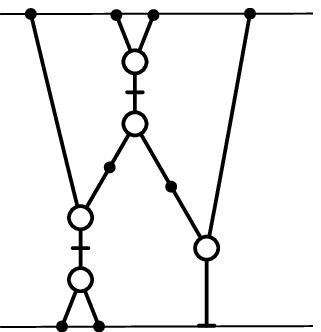}}
   \put(3.5,3.1){$1$}
   \put(3.5,0.0){$0$}
  \end{picture}}
  \caption{}
  \label{fig:31.2}
\end{center}\end{figure}

\end{rem}

In the following, we regard fork complexes as geometric objects, i.e., $1$-dimensional polyhedra. 

\begin{defi}
\rm
A \textit{fork of} $\mathscr{F}$ is the image of a fork in $\mathcal{A}\cup \mathcal{B}$ in $\mathscr{F}$. 
A \textit{grip} (\textit{root} and \textit{tine} resp.) \textit{of} $\mathscr{F}$ is the image of a 
grip (root and tine resp.) in $\mathcal{A}\cup \mathcal{B}$ in $\mathscr{F}$. 
\end{defi}

\begin{defi}
\rm
Let $M$ be a compact orientable $3$-manifold, and let $(\partial_1 M,\partial_2 M)$ be 
a partition of boundary components of $M$. 
A \textit{generalized Heegaard splitting} of $(M;\partial_1 M,\partial_2 M)$ is a pair of 
an exact fork complex $\mathscr{F}$ and a proper map $\rho:(M;\partial_1 M,\partial_2 M)\to 
(\mathscr{F};\partial_1 \mathscr{F},\partial_2 \mathscr{F})$ which satisfies 
the following. 

\begin{enumerate}
\item The map $\rho$ is transverse to $\mathscr{F}-\{$the roots of $\mathscr{F}\}$. 
\item For each fork $\mathcal{F}\subset \mathscr{F}$, we have the following (\textit{cf}. Figure \ref{fig:31.3}). 
\begin{enumerate}
\item 
If $\mathcal{F}$ is a $0$-fork, then $\rho^{-1}(\mathcal{F})$ is a handlebody 
$V_{\mathcal{F}}$ such that (1) $\rho^{-1}(g)=\partial V_{\mathcal{F}}$ 
and (2) $\rho^{-1}(p)$ is a $1$-complex which is a spine of $V_{\mathcal{F}}$, 
where $g$ is the grip of $\mathcal{F}$. 
\item
If $\mathcal{F}$ is an $n$-fork with $n \ge 1$, then $\rho^{-1}(\mathcal{F})$ is a connected compression body 
$V_{\mathcal{F}}$ such that (1) $\rho^{-1}(g)=\partial_+ V_{\mathcal{F}}$, (2) for each tine $t_i$, 
$\rho^{-1}(t_i)$ is a connected component of $\partial_- V_{\mathcal{F}}$ and 
$\rho^{-1}(t_i)\neq \rho^{-1}(t_j)$ for $i\neq j$ and (3) $\rho^{-1}(p)$ is a $1$-complex which is 
a deformation retract of $V_{\mathcal{F}}$, where $g$ is the grip of $\mathcal{F}$, 
$p$ is the root of $\mathcal{F}$ and $\{t_i\}_{1\le i\le n}$ is the set of the tines of $\mathcal{F}$. 
\end{enumerate}
\end{enumerate}

\begin{figure}[htb]\begin{center}
  {\unitlength=1cm
  \begin{picture}(6.7,4.2)
   \put(0,2.6){\includegraphics[keepaspectratio]{%
               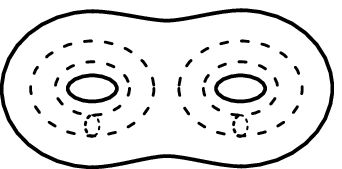}}
   \put(0,0){\includegraphics[keepaspectratio]{%
               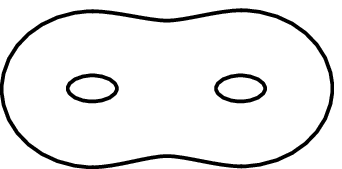}}
   \put(4.0,2.6){\includegraphics[keepaspectratio,scale=1.2]{%
               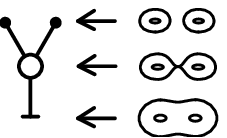}}
   \put(4.2,.1){\includegraphics[keepaspectratio,scale=1.2]{%
               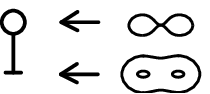}}
  \end{picture}}
  \caption{}
  \label{fig:31.3}
\end{center}\end{figure}

\end{defi}

\begin{rem}
\rm
Let $g$ be a grip of $\mathscr{F}$ which is contained in the interior of $\mathscr{F}$. 
Let $\mathcal{F}_1$ and $\mathcal{F}_2$ be the forks of $\mathscr{F}$ which are adjacent to $g$. 
Then $(\rho^{-1}(\mathcal{F}_1),\rho^{-1}(\mathcal{F}_2);\rho^{-1}(g))$ is a Heegaard splitting of 
$\rho^{-1}(\mathcal{F}_1\cup \mathcal{F}_2)$. 
\end{rem}

\begin{defi}
\rm
A generalized Heegaard splitting $(\mathscr{F},\rho)$ is said to be \textit{strongly irreducible} if 
$(1)$ for each tine $t$, $\rho^{-1}(t)$ is incompressible, and $(2)$ for each grip $g$ with 
two forks attached to $g$, say $\mathcal{F}_1$ and $\mathcal{F}_2$, 
$(\rho^{-1}(\mathcal{F}_1),\rho^{-1}(\mathcal{F}_2);\rho^{-1}(g))$ is strongly irreducible. 
\end{defi}

Let $\mathcal{M}$ be the set of finite multisets of $\mathbb{Z}_{\ge 0}=\{0,1,2, ...\}$. 
We define a total order $<$ on $\mathcal{M}$ as follows. 
For $M_1$ and $M_2\in \mathcal{M}$, we first arrange the elements of $M_i$ $(i=1,2)$ in 
non-increasing order respectively. Then we compare the arranged tuples of non-negative integers 
by lexicographic order. 

\begin{example}\label{E1}
\rm
$(1)$\ \ If $M_1=\{5,4,1,1\}$ and $M_2=\{5,3,2,2,2,1\}$, then $M_2<M_1$. \\
$(2)$\ \ If $M_1=\{3,1,0,0\}$ and $M_2=\{3,1,0,0,0\}$, then $M_1<M_2$. 
\end{example}

\begin{defi}
\rm
Let $(\mathscr{F},\rho)$ be a generalized Heegaard splitting of $(M;\partial_1 M,\partial_2 M)$. 
We define \textit{the width of} $(\mathscr{F},\rho)$ to be the multiset 

\bigskip
\begin{center}
$w(\mathscr{F},\rho)=\{\textit{genus}(\rho^{-1}(g_1)),\dots , \textit{genus}(\rho^{-1}(g_m))\}$, 
\end{center}

\medskip\noindent
where $\{g_1,\dots ,g_m\}$ is the set of the grips of $\mathscr{F}$. 
We say that $(\mathscr{F},\rho)$ is \textit{thin} 
if $w(\mathscr{F},\rho)$ is minimal among all generalized Heegaard splittings of 
$(M;\partial_1 M,\partial_2 M)$. 
\end{defi}

\begin{example}\label{thinB3}
The thin generalized Heegaard splittings of the $3$-ball $B^3$ are two fork complexes illustrated in 
Figure \ref{fig:79}, 
where $\rho^{-1}(\mathcal{F}_1)$ is a $3$-ball and $\rho^{-1}(\mathcal{F}_2)\cong S^2\times [0,1]$. 

\begin{figure}[htb]\begin{center}
  {\unitlength=1cm
  \begin{picture}(4.2,.8)
   \put(0,.4){\includegraphics[keepaspectratio]{
   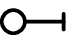}}
   \put(2.4,.4){\includegraphics[keepaspectratio]{
   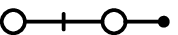}}
   \put(.2,0){$\mathcal{F}_1$}
   \put(2.6,0){$\mathcal{F}_1$}
   \put(3.4,0){$\mathcal{F}_2$}
  \end{picture}}
  \caption{}
  \label{fig:79}
\end{center}\end{figure}

\end{example}

\subsection{Properties of thin generalized Heegaard splittings}

In this subsection, let $(\mathscr{F},\rho)$ be a thin generalized Heegaard 
splitting of $(M;\partial_1 M,\partial_2 M)$. 

\begin{obs}\label{S9}
Let $t$ be a tine of $\mathscr{F}$. Then any $2$-sphere component of $\rho^{-1}(t)$ is essential 
in $M$ unless $M$ is a $3$-ball. 
\end{obs}

\begin{proof}
Suppose that there is a tine $t$ such that $\rho^{-1}(t)$ is a $2$-sphere, say $P$, which 
bounds a $3$-ball $B$ in $M$. 
Let $\mathscr{F}_B$ be the subcomplex of $\mathscr{F}$ with $\rho^{-1}(\mathscr{F}_B)=B$. 
If $\mathscr{F}_B=\mathscr{F}$, then we see that $M$ is a $3$-ball. 
Otherwise, there is a fork $\mathcal{F}'$ with $t\in \mathcal{F}'$ and 
$\mathcal{F}'\not\subset \mathscr{F}_B$. Let $e_t$ be the $1$-simplex in $\mathcal{F}'$ joining $t$ 
to the root of $\mathcal{F}'$. Set $\mathscr{F}^{\ast}=\mathscr{F}\setminus (\mathscr{F}_B\cup e_t)$. 
Note that $\rho^{-1}(\mathcal{F}'\cup \mathscr{F}_B)$ ($=\rho^{-1}(\mathcal{F}')\cup B$) is a 
compression body $V^{\ast}$. 
Then it is easy to see that we can modify $\rho$ in $V^{\ast}$ to obtain 
$\rho^{\ast}:M\to \mathscr{F}^{\ast}$ such that $(\rho^{\ast})^{-1}(\mathcal{F}'\setminus e_t)$ 
is the compression body $V^{\ast}$ (\textit{cf}. Figure \ref{fig:31.31}). 

\begin{figure}[htb]\begin{center}
  {\unitlength=1cm
  \begin{picture}(6.3,4.2)
   \put(0,3.1){\includegraphics[keepaspectratio]{%
               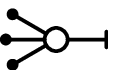}}
   \put(2.3,2.6){\includegraphics[keepaspectratio]{%
               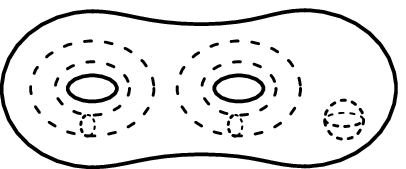}}
   \put(0,0.5){\includegraphics[keepaspectratio]{%
               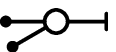}}
   \put(2.3,0){\includegraphics[keepaspectratio]{%
               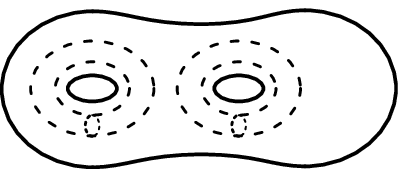}}
   \put(-.2,3.7){$t$}
   \put(1.6,3.6){$\rho$}
   \put(1.6,1.0){$\rho^*$}
   \put(1.4,.7){$\longleftarrow$}
   \put(1.4,3.34){$\longleftarrow$}
   \put(4.3,2.0){$\downarrow$}
  \end{picture}}
  \caption{}
  \label{fig:31.31}
\end{center}\end{figure}

Moreover, the generalized Heegaard structure on $(\mathscr{F},\rho)$ 
(e.g. $\mathcal{A},\mathcal{B}$ decomposition etc) is naturally inherited to $(\mathscr{F}^{\ast},\rho^{\ast})$. 
Then we clearly have $w(\mathscr{F}^{\ast},\rho^{\ast})<w(\mathscr{F},\rho)$, contradicting the assumption 
that $(\mathscr{F},\rho)$ is thin. 

\end{proof}

\begin{lem}\label{S10}
Suppose that there is a fork $\mathcal{F}$ such that $\rho^{-1}(t)$ is trivial. Let $t$ be the 
tine of $\mathcal{F}$. Then $\rho^{-1}(t)$ is a component of $\partial M$ and one of the following holds. 

\begin{enumerate}
\item $M$ is a $3$-ball. 
\item $M\cong \rho^{-1}(t)\times [0,1]$. 
\item $\rho^{-1}(t)$ is compressible in $M$.
\end{enumerate}
\end{lem}

\begin{proof}
We first prove that $\rho^{-1}(t)$ is a boundary component of $M$. Suppose that $\rho^{-1}(t)$ 
is not a boundary component of $M$. Let $g$ be the grip of $\mathcal{F}$. If $\rho^{-1}(g)$ is a 
boundary component of $M$, then we can reduce the width by removing $\mathcal{F}$, a contradiction 
(\textit{cf}. Figure \ref{fig:31.5}). 

\begin{figure}[htb]\begin{center}
  {\unitlength=1cm
  \begin{picture}(5.8,2.7)
   \put(0,0.2){\includegraphics[keepaspectratio]{%
             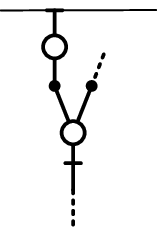}}
   \put(4.0,0.2){\includegraphics[keepaspectratio]{%
             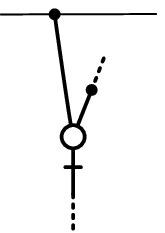}}
   \put(2.7,1.3){$\Longrightarrow$}
   \put(1.8,2.3){$1$}
   \put(5.8,2.3){$1$}
   \put(-.2,1.9){$F$}
  \end{picture}}
  \caption{}
  \label{fig:31.5}
\end{center}\end{figure}

Hence $\mathcal{F}$ is contained in the interior of $\mathscr{F}$. Note that $\mathcal{F}$ is a $1$-fork. 
Let $\mathcal{F}_1$ be the fork attached to $g$ and $\mathcal{F}_2$ the fork attached to $t$. Note that 
since $\rho^{-1}(\mathcal{F})$ is a trivial compression body, we see that 
$\rho^{-1}(\mathcal{F}_1\cup \mathcal{F}\cup \mathcal{F}_2)$ is also a compression body. Hence we can 
replace $\mathcal{F}_1\cup \mathcal{F}\cup \mathcal{F}_2$ in $\mathscr{F}$ to a new fork 
so that we can obtain a new fork complex, say $\mathscr{F}^{\ast}$. Moreover, 
we can modify $\rho:M\to \mathscr{F}$ to obtain 
$\rho^{\ast}:M\to \mathscr{F}^{\ast}$ so that $(\mathscr{F}^{\ast},\rho^{\ast})$ is a generalized Heegaard 
splitting of $(M;\partial_1 M,\partial_2 M)$ with $w(\mathscr{F}^{\ast},\rho^{\ast})<w(\mathscr{F},\rho)$, 
a contradiction (\textit{cf}. Figure \ref{fig:31.6}). Hence $\rho^{-1}(t)$ is a boundary component of $M$. 

\begin{figure}[htb]\begin{center}
  {\unitlength=1cm
  \begin{picture}(5.8,3.4)
   \put(0.2,0.2){\includegraphics[keepaspectratio]{%
             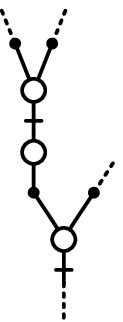}}
   \put(4.2,0.2){\includegraphics[keepaspectratio]{%
             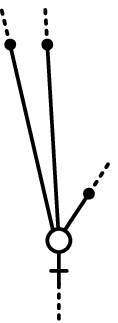}}
   \put(2.7,1.3){$\Longrightarrow$}
   \put(-.2,2.5){$F_1$}
   \put(-.2,1.8){$F$}
   \put(-.2,0.9){$F_2$}
  \end{picture}}
  \caption{}
  \label{fig:31.6}
\end{center}\end{figure}

We next show that one of the conclusions $(1)$-$(3)$ of Lemma \ref{S10} holds. Suppose that 
both conclusions $(1)$ and $(2)$ of Lemma \ref{S10} do not hold, i.e., $M$ is not a $3$-ball and 
$M\not\cong \rho^{-1}(t)\times [0,1]$. Then there is a fork $\mathcal{F}'(\ne \mathcal{F})$ attached to $g$. 
Moreover, since $(\mathscr{F},\rho)$ is thin and $M\not\cong \rho^{-1}(t)\times [0,1]$, 
we see that $\rho^{-1}(\mathcal{F}')$ is a non-trivial compression body. 
Also, since $M$ is not a $3$-ball, $\rho^{-1}(t)$ is not a $2$-sphere. Hence we see that $\rho^{-1}(t)$ 
is compressible in $\rho^{-1}(\mathcal{F})\cup \rho^{-1}(\mathcal{F}')$. This implies that the conclusion $(3)$ 
of Lemma \ref{S10} holds. 
\end{proof}

\begin{prop}\label{S11}
Let $\mathcal{F}_1$ and $\mathcal{F}_2$ be forks of 
$\mathscr{F}$ which have the same grip $g$ of $\mathscr{F}$. Then $(\rho^{-1}(\mathcal{F}_1),
\rho^{-1}(\mathcal{F}_2);\rho^{-1}(g))$ is strongly irreducible. 
\end{prop}

\begin{proof}
Set $A_g=\rho^{-1}(\mathcal{F}_1)$, $B_g=\rho^{-1}(\mathcal{F}_2)$, $S_g=\rho^{-1}(g)$, 
$M_g=A_g\cup B_g$, $\partial_1 M_g=\partial_- A_g$ and $\partial_2 M_g=\partial_- B_g$. 
Then $(A_g,B_g;S_g)$ is a Heegaard splitting of $(M_g;\partial_1 M_g,\partial_2 M_g)$. 
Suppose that $(A_g,B_g;S_g)$ is weakly reducible. Let $D_A$ and $D_B$ 
be meridian disks of $A_g$ and $B_g$ respectively which satisfy $\partial D_A\cap \partial D_B=\emptyset$. 
Let $\Delta_A$ ($\Delta_B$ resp.) be a complete meridian system of $A_g$ ($B_g$ resp.) such that 
$D_A$ ($D_B$ resp.) is a component of $\Delta_A$ ($\Delta_B$ resp.) (\textit{cf}. (6) of Remark \ref{R1}). 
Note that $A_g$ is obtained from $\partial_- A_g\times [0,1]$ and $0$-handles $\mathcal{H}^0$ by 
attaching $1$-handles $\mathcal{H}^1$ corresponding to $\Delta_A$ (\textit{cf}. (3) of Remark \ref{R1}) 
and that $B_g$ is obtained from $S_g\times [0,1]$ by attaching $2$-handles $\mathcal{H}^2$ corresponding 
to $\Delta_B$ and $3$-handles $\mathcal{H}^3$ (\textit{cf}. Definition \ref{D1}). Hence we see that 
$M_g$ admits the following decomposition (\textit{cf}. Remark \ref{R2}): 

\bigskip
\begin{center}
$M_g=(\partial_1 M_g\times [0,1])\cup \mathcal{H}^0 \cup \mathcal{H}^1 \cup \mathcal{H}^2\cup \mathcal{H}^3$. 
\end{center}

\medskip\noindent
Let $h^1$ be the component of $\mathcal{H}^1$ corresponding to $D_A$ and 
$h^2$ the component of $\mathcal{H}^2$ corresponding to $D_B$. Then $M_g$ admits the following 
decomposition: 

\bigskip
\begin{center}
$M_g=
(\partial_1 M_g\times [0,1])\cup \mathcal{H}^0 \cup (\mathcal{H}^1\setminus h^1)\cup h^2\cup h^1\cup 
(\mathcal{H}^2\setminus h^2)\cup \mathcal{H}^3$. 
\end{center}

\medskip\noindent
Set $A_g'=(\partial_1 M_g\times [0,1])\cup \mathcal{H}^0 \cup (\mathcal{H}^1\setminus h^1)$. 
We divide the proof into the following two cases. 

\medskip
\noindent\textit{Case 1.}\ \ \ \ 
$\partial D_A$ or $\partial D_B$ is non-separating in $S_g$. 

\medskip
Suppose first that $\partial D_A$ is non-separating in $S_g$. 
Then $A_g'$ is a compression body (\textit{cf}. (6) of Remark \ref{R1}). 
Since $A_g'=(\partial_1 M_g\times [0,1])\cup \mathcal{H}^0 \cup (\mathcal{H}^1\setminus h^1)$, we obtain 

\bigskip
\begin{center}
$M_g = A_g'\cup h^2\cup h^1\cup (\mathcal{H}^2\setminus h^2)\cup \mathcal{H}^3$.
\end{center}

\medskip\noindent
Note that the attaching region of the $2$-handle $h^2$ is contained in $\partial_+ A_g'$. 
Hence we have: 

\bigskip
\begin{center}
$M_g \cong 
A_g'\cup \bigl( (\partial_+ A_g'\times [0,1])\cup h^2 \bigr)
\cup h^1\cup (\mathcal{H}^2\setminus h^2)\cup \mathcal{H}^3$.
\end{center}

\medskip\noindent
Set $B_g'=(\partial_+ A_g'\times [0,1])\cup h^2$. Then $B_g'$ is also a compression body and we have: 

\bigskip
\begin{center}
$M_g \cong  A_g'\cup B_g'\cup h^1\cup (\mathcal{H}^2\setminus h^2)\cup \mathcal{H}^3$.
\end{center}

\medskip\noindent
Note that $\partial_- B_g'$ is homeomorphic to the surface obtained from $S_g$ by performing surgery 
along $\partial D_A\cup \partial D_B$. Then we have the following subcases. 

\medskip
\noindent\textit{Case 1.1.}\ \ \ \ 
$\partial D_B$ is non-separating in $S_g$ and 
$\partial D_A\cup \partial D_B$ is non-separating in $S_g$. 

\medskip
Then $\partial_- B_g'$ is connected. Note that 

\begin{eqnarray*}
M_g 
&\cong & A_g'\cup B_g'\cup h^1\cup (\mathcal{H}^2\setminus h^2)\cup \mathcal{H}^3\\
&\cong & A_g'\cup B_g'\cup \bigl( (\partial_- B_g'\times [0,1])\cup h^1\bigr) 
\cup (\mathcal{H}^2\setminus h^2)\cup \mathcal{H}^3.
\end{eqnarray*}

\noindent
Set $A_g''=(\partial_- B_g'\times [0,1])\cup h^1$. 
Then $A_g''$ is also a compression body and we have: 

\begin{eqnarray*}
M_g 
&\cong & A_g'\cup B_g'\cup A_g''\cup (\mathcal{H}^2\setminus h^2)\cup \mathcal{H}^3\\
&\cong & A_g'\cup B_g'\cup A_g''\cup \bigl( (\partial_+ A_g''\times [0,1]\cup 
(\mathcal{H}^2\setminus h^2)\cup \mathcal{H}^3 \bigr). 
\end{eqnarray*}

\noindent
Set $B_g''=\partial_+ A_g''\times [0,1]\cup (\mathcal{H}^2\setminus h^2)\cup \mathcal{H}^3$. 
Note that $B_g'\cap A_g''=\partial_- B_g'=\partial_- A_g''$. This shows that each handle of 
$\mathcal{H}^2\setminus h^2$ and $\mathcal{H}^3$ is adjacent to $A_g''$ along $\partial_+ A_g''$. 
This implies that $B_g''$ is also a compression body. Hence we have: 

\bigskip
\begin{center}
$M_g \cong (A_g'\cup B_g')\cup (A_g''\cup B_g'')$. 
\end{center}

\medskip\noindent
Then we can substitute $\mathcal{F}_1\cup \mathcal{F}_2$ in $\mathscr{F}$ for 
$\mathcal{F}_1'\cup \mathcal{F}_2'\cup \mathcal{F}_1''\cup \mathcal{F}_2''$, where 
$\mathcal{F}_1'$, $\mathcal{F}_2'$, $\mathcal{F}_1''$ and $\mathcal{F}_2''$ are forks 
corresponding to $A_g'$, $B_g'$, $A_g''$ and $B_g''$ respectively. 
Set $\mathscr{F}^{\ast}=(\mathscr{F}\setminus (\mathcal{F}_1\cup \mathcal{F}_2))\cup 
(\mathcal{F}_1'\cup \mathcal{F}_2'\cup \mathcal{F}_1''\cup \mathcal{F}_2'')$. 
Then we can modify $\rho:M\to \mathscr{F}$ in $M_g$ to obtain 
$\rho^{\ast}:M\to \mathscr{F}^{\ast}$ such that $(\rho^{\ast})^{-1}(\mathcal{F}_1')=A_g'$, 
$(\rho^{\ast})^{-1}(\mathcal{F}_2')=B_g'$, $(\rho^{\ast})^{-1}(\mathcal{F}_1'')=A_g''$ and 
$(\rho^{\ast})^{-1}(\mathcal{F}_2'')=B_g''$. It is easy to see that $w(\mathscr{F}^{\ast},\rho^{\ast})
<w(\mathscr{F},\rho)$, a contradiction (\textit{cf}. Figure \ref{fig:32}). 

\begin{figure}[htb]\begin{center}
  {\unitlength=1cm
  \begin{picture}(10.5,8)
   \put(0.4,4.9){\includegraphics[keepaspectratio]{%
   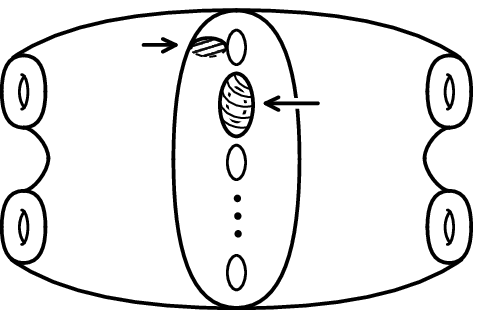}}
   \put(7.4,6){\includegraphics[keepaspectratio]{%
   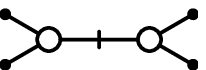}}
   \put(0,0){\includegraphics[keepaspectratio]{%
   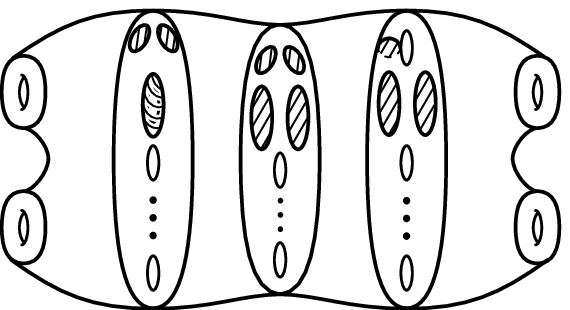}}
   \put(6.4,1.1){\includegraphics[keepaspectratio]{%
   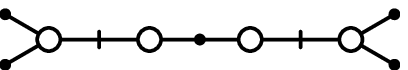}}
   \put(1.3,7.5){$D_A$}
   \put(1.3,5.5){$A_g$}
   \put(3.8,6.9){$D_B$}
   \put(4,5.5){$B_g$}
   \put(2.7,3.8){{\Large $\Downarrow$}}
   \put(8.25,3.8){{\Large $\Downarrow$}}
%
   \put(.6,.4){$A'_g$}
   \put(2,.4){$B'_g$}
   \put(3.3,.4){$A''_g$}
   \put(4.7,.4){$B''_g$}
  \end{picture}}
  \caption{}
  \label{fig:32}
\end{center}\end{figure}

\medskip
\noindent\textit{Case 1.2.}\ \ \ \ 
$\partial D_B$ is non-separating in $S_g$ and 
$\partial D_A\cup \partial D_B$ is separating in $S_g$. 

\medskip
Then $\partial_- B_g'$ consists of two components, say $G_1$ and $G_2$. 

\begin{eqnarray*}
M_g 
&\cong & A_g'\cup B_g'\cup h^1\cup (\mathcal{H}^2\setminus h^2)\cup \mathcal{H}^3\\
&\cong & A_g'\cup B_g'\cup (\partial_- B_g'\times [0,1])\cup h^1\cup (\mathcal{H}^2\setminus h^2)\cup 
\mathcal{H}^3\\
&\cong & A_g'\cup B_g'\cup ((G_1\cup G_2)\times [0,1])\cup h^1\cup (\mathcal{H}^2\setminus h^2)\cup 
\mathcal{H}^3.
\end{eqnarray*}

\noindent
Set $A_g''=((G_1\cup G_2)\times [0,1])\cup h^1$. Since $\partial D_B$ is 
non-separating in $S_g$, we see that $h^1$ joins $G_1$ to $G_2$. 
Hence $A_g''$ is a compression body and we have: 

\begin{eqnarray*}
M_g 
&\cong & A_g'\cup B_g'\cup A_g''\cup (\mathcal{H}^2\setminus h^2)\cup \mathcal{H}^3\\
&\cong & A_g'\cup B_g'\cup A_g''\cup (\partial_+ A_g''\times [0,1])\cup 
(\mathcal{H}^2\setminus h^2)\cup \mathcal{H}^3.
\end{eqnarray*}

\noindent
Set $B_g''=\partial_+ A_g''\times [0,1]\cup (\mathcal{H}^2\setminus h^2)\cup \mathcal{H}^3$. 
Note that $B_g'\cap A_g''=\partial_- B_g'=\partial_- A_g''$. This shows that each handle of 
$\mathcal{H}^2\setminus h^2$ and $\mathcal{H}^3$ is adjacent to $A_g''$ along $\partial_+ A_g''$. 
This implies that $B_g''$ is also a compression body. Hence we have: 

\bigskip
\begin{center}
$M_g \cong  (A_g'\cup B_g')\cup (A_g''\cup B_g'')$. 
\end{center}

\medskip\noindent
According to this decomposition, we can modify the fork complex $(\mathscr{F},\rho)$ as in Figure \ref{fig:33} 
or Figure \ref{fig:33.1}. 
It is easy to see that for a new complex $(\mathscr{F}^{\ast},\rho^{\ast})$, we have 
$w(\mathscr{F}^{\ast},\rho^{\ast})<w(\mathscr{F},\rho)$, a contradiction. 

\begin{figure}[htb]\begin{center}
  {\unitlength=1cm
  \begin{picture}(10.5,8)
   \put(0.4,4.9){\includegraphics[keepaspectratio]{%
   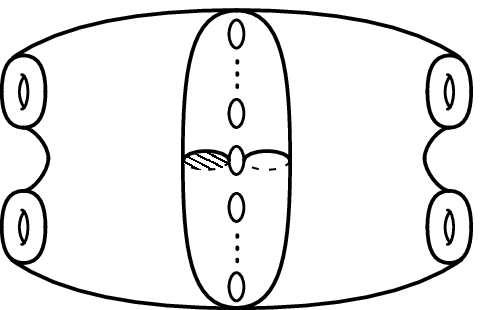}}
   \put(7.4,6){\includegraphics[keepaspectratio]{%
   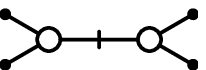}}
   \put(0,0){\includegraphics[keepaspectratio]{%
   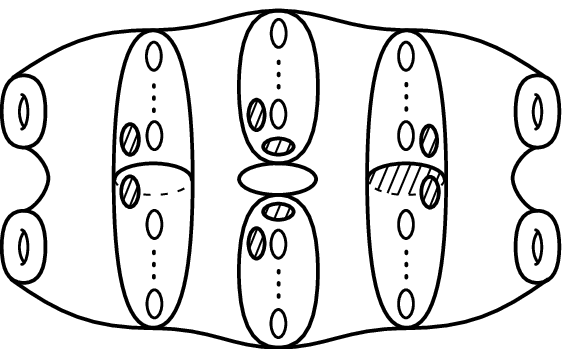}}
   \put(6.4,1.1){\includegraphics[keepaspectratio]{%
   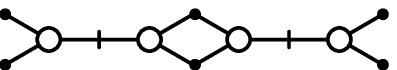}}
   \put(1.6,6.3){$D_A$}
   \put(1.3,5.5){$A_g$}
   \put(3.46,6.3){$D_B$}
   \put(4,5.5){$B_g$}
   \put(2.7,3.8){{\Large $\Downarrow$}}
   \put(8.25,3.8){{\Large $\Downarrow$}}
%
   \put(.6,.7){$A'_g$}
   \put(2,.4){$B'_g$}
   \put(3.3,.4){$A''_g$}
   \put(4.7,.7){$B''_g$}
  \end{picture}}	       
  \caption{\sl The case of irreducible splittings}
  \label{fig:33}
\end{center}\end{figure}

\begin{figure}[htb]\begin{center}
  {\unitlength=1cm
  \begin{picture}(10.5,8)
   \put(0.4,4.9){\includegraphics[keepaspectratio]{%
   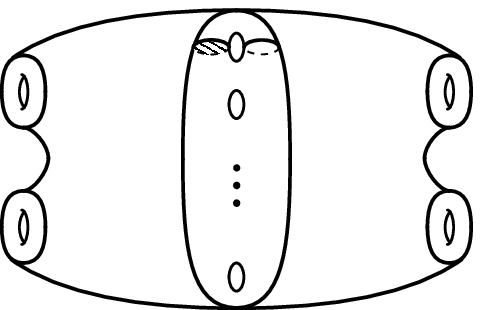}}
   \put(7.4,6){\includegraphics[keepaspectratio]{%
   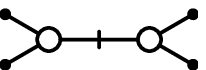}}
   \put(0,0){\includegraphics[keepaspectratio]{%
   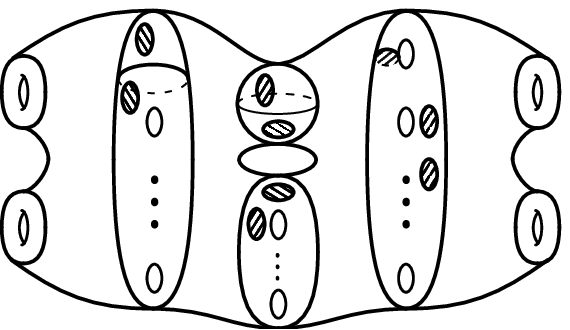}}
   \put(6.4,1.1){\includegraphics[keepaspectratio]{%
   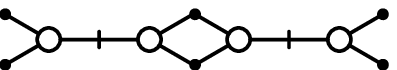}}
   \put(1.7,7.5){$D_A$}
   \put(1.3,5.5){$A_g$}
   \put(3.4,7.5){$D_B$}
   \put(4,5.5){$B_g$}
   \put(2.7,3.8){{\Large $\Downarrow$}}
   \put(8.25,3.8){{\Large $\Downarrow$}}
%
   \put(.6,.7){$A'_g$}
   \put(2,.4){$B'_g$}
   \put(3.3,.4){$A''_g$}
   \put(4.7,.7){$B''_g$}
  \end{picture}}	       
  \caption{\sl The case of reducible splittings}
  \label{fig:33.1}
\end{center}\end{figure}

\medskip
\noindent\textit{Case 1.3.}\ \ \ \ 
$\partial D_B$ is separating in $S_g$ 
(hence $\partial D_A\cup \partial D_B$ is separating in $S_g$). 

\medskip
Then $\partial_- B_g'$ consists of two components, say $\bar{G}_1$ and $\bar{G}_2$. Since $\partial D_B$ is 
separating in $S_g$, we see that $h^1$ joins $\bar{G}_1$ or $\bar{G}_2$, say $\bar{G}_1$, to itself. 
Let $\mathcal{H}_1^2$ ($\mathcal{H}_2^2$ resp.) be the components of $\mathcal{H}^2\setminus h^2$ 
adjacent to $\bar{G}_1$ ($\bar{G}_2$ resp.). 
Let $\mathcal{H}_1^3$ ($\mathcal{H}_2^3$ resp.) be the components of $\mathcal{H}^3$ 
adjacent to $\bar{G}_1$ ($\bar{G}_2$ resp.). Set 
$\bar{B}_g'=B_g'\cup \mathcal{H}_2^2\cup \mathcal{H}_2^3$. 
Then $\bar{B}_g'$ is a compression body with $\partial_+ \bar{B}_g'=\partial_+ B_g'$. 
Set $A_g''=(\bar{G}_1\times [0,1])\cup h^1$ and $B_g''=(\partial_+ A_g''\times [0,1])\cup 
\mathcal{H}_1^2\cup \mathcal{H}_1^3$. Then each of $A_g''$ and $B_g''$ is a compression body. 
Hence we have: 

\begin{eqnarray*}
M_g 
&\cong & A_g'\cup B_g'\cup h^1\cup (\mathcal{H}^2_1\cup \mathcal{H}^2_2)\cup 
(\mathcal{H}^3_1\cup \mathcal{H}^3_2)\\
&\cong & A_g'\cup (B_g'\cup \mathcal{H}_2^2\cup \mathcal{H}_2^3)
\cup h^1\cup \mathcal{H}^2_1\cup \mathcal{H}^3_1\\
&\cong & A_g'\cup \bar{B}_g'\cup (\bar{G}_1\times [0,1])\cup h^1 \cup 
\mathcal{H}^2_1\cup \mathcal{H}^3_1\\
&\cong & A_g'\cup \bar{B}_g'\cup A_g''\cup (\partial_+ A_g''\times [0,1])\cup 
\mathcal{H}^2_1\cup \mathcal{H}^3_1\\
&\cong &  (A_g'\cup \bar{B}_g')\cup (A_g''\cup B_g'').
\end{eqnarray*}

\noindent
According to this decomposition, we can modify the fork complex $(\mathscr{F},\rho)$ as in Figure \ref{fig:34}. 
It is easy to see that for a new complex $(\mathscr{F}^{\ast},\rho^{\ast})$, we have 
$w(\mathscr{F}^{\ast},\rho^{\ast})<w(\mathscr{F},\rho)$, a contradiction. 
Therefore if $\partial D_A$ is non-separating, we have the desired conclusion. 

\begin{figure}[htb]\begin{center}
  {\unitlength=1cm
  \begin{picture}(10.5,9)
   \put(0.4,5.8){\includegraphics[keepaspectratio]{%
   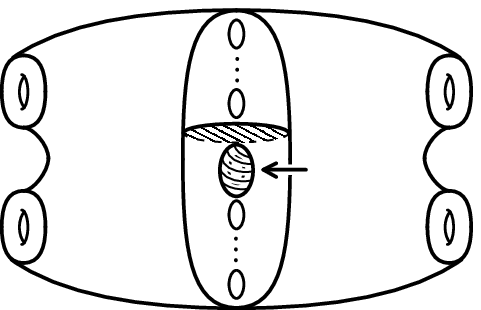}}
   \put(7.4,6.9){\includegraphics[keepaspectratio]{%
   Fig33p2.eps}}
   \put(0,0){\includegraphics[keepaspectratio]{%
   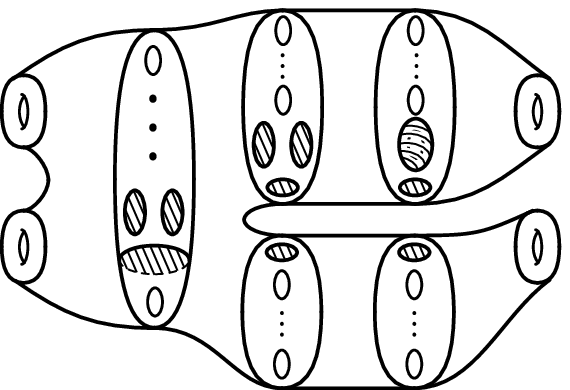}}
   \put(6.4,1.6){\includegraphics[keepaspectratio]{%
   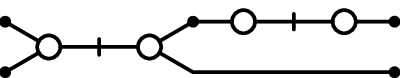}}
   \put(1.7,7.5){$D_A$}
   \put(1.3,6.4){$A_g$}
   \put(3.6,7.16){$D_B$}
   \put(4,6.4){$B_g$}
   \put(2.7,4.7){{\Large $\Downarrow$}}
   \put(8.25,4.7){{\Large $\Downarrow$}}
%
%
   \put(.6,1.2){$A'_g$}
   \put(2.05,1.26){$\bar{B}'_g$}
   \put(3.3,2.2){$A''_g$}
   \put(4.7,2.3){$B''_g$}
  \end{picture}}
  \caption{}
  \label{fig:34}
\end{center}\end{figure}

\medskip
Suppose next that $\partial D_B$ is non-separating in $S_g$. 
Then we start with the dual handle decomposition 

\bigskip
\begin{center}
$M_g=(\partial_2 M_g\times [0,1])\cup \bar{\mathcal{H}}^0 
\cup \bar{\mathcal{H}}^1 \cup \bar{\mathcal{H}}^2 \cup \bar{\mathcal{H}}^3$ 
\end{center}

\medskip\noindent
and apply the above arguments which gives a contradiction. 

\medskip
\noindent\textit{Case 2.}\ \ \ \ 
Each of $\partial D_A$ and $\partial D_B$ is separating in $S$. 

\medskip
Then $A_g'$ consists of two compression bodies, say $\bar{A}_g'$ and $\tilde{A}_g'$ 
(\textit{cf}. (6) of Remark \ref{R1}). 
We may suppose that $h^2$ is attached to $\partial_+ \bar{A}_g'$. Set 
$B_g'=(\partial_+ \bar{A}_g'\times [0,1])\cup h^2$. Since $D_B$ is separating in $S$, 
we see that $\partial_- B_g'$ consists of two components, say $G_1$ and $G_2$. 
Note that $D_A$ is also separating in $S_g$. Hence we may suppose that $h^1\cap G_2\ne \emptyset$ and 
$h^1\cap G_1=\emptyset$. 
Let $\mathcal{H}_1^2$ be the components of $\mathcal{H}^2\setminus h^2$ adjacent to $G_1$ and 
$\mathcal{H}_1^3$ be the components of $\mathcal{H}^3$ adjacent to $G_1$. 
Set $\mathcal{H}_2^2=\mathcal{H}^2\setminus (h^2\cup \mathcal{H}_1^2)$, 
$\mathcal{H}_2^3=\mathcal{H}^3\setminus \mathcal{H}_1^3$ and 
$\bar{B}_g'=B_g'\cup \mathcal{H}^2_1\cup \mathcal{H}^3_1$. 
Then $\bar{B}_g'$ is a compression body. 
Set $A_g^{\ast}=(G_2\times [0,1])\cup \tilde{A}_g'\cup h_1$ and 
$B_g''=(\partial_+ A_g''\times [0,1])\cup \mathcal{H}^2_2\cup \mathcal{H}^3_2$. 
Note that each of $A_g''$ and $B_g''$ is a compression body. 
Set $A_g''=\tilde{A}_g'\cup A_g^{\ast}$. 
Note also that $A_g''$ is a compression body (\textit{cf}. (5) of Remark \ref{R1}). 
Hence we have: 

\begin{eqnarray*}
M_g 
&\cong & A_g'\cup h^2\cup h^1\cup (\mathcal{H}^2_1\cup \mathcal{H}^2_2)\cup 
(\mathcal{H}^3_1\cup \mathcal{H}^3_2)\\
&\cong & (\bar{A}_g'\cup \tilde{A}_g')\cup h^2\cup h^1\cup (\mathcal{H}^2_1\cup \mathcal{H}^2_2)\cup 
(\mathcal{H}^3_1\cup \mathcal{H}^3_2)\\
&\cong & \bar{A}_g'\cup \bigl( (\partial_+ \bar{A}_g'\times [0,1])\cup h^2\bigr) \cup \tilde{A}_g' 
\cup h^1\cup (\mathcal{H}^2_1\cup \mathcal{H}^2_2)\cup 
(\mathcal{H}^3_1\cup \mathcal{H}^3_2)\\
&\cong & \bar{A}_g'\cup B_g'\cup ((G_1\cup G_2)\times [0,1])\cup \tilde{A}_g' 
\cup h^1\cup (\mathcal{H}^2_1\cup \mathcal{H}^2_2)\cup 
(\mathcal{H}^3_1\cup \mathcal{H}^3_2)\\
&\cong & \bar{A}_g'\cup B_g'\cup ((G_1\cup G_2)\times [0,1])\cup \tilde{A}_g'
\cup h^1\cup (\mathcal{H}^2_1\cup \mathcal{H}^2_2)\cup 
(\mathcal{H}^3_1\cup \mathcal{H}^3_2)\\
&\cong & \bar{A}_g'\cup \bigl( B_g'\cup 
(G_1\times [0,1])\cup \mathcal{H}^2_1\cup \mathcal{H}^3_1 \bigr) \cup \tilde{A}_g'\\
&&
\cup \bigl( (G_2\times [0,1])\cup \tilde{A}_g'\cup h_1\bigr) \cup 
\mathcal{H}^2_2\cup \mathcal{H}^3_2\\
&\cong & \bar{A}_g'\cup \bar{B}_g'\cup (\tilde{A}_g'\cup A_g^{\ast})\cup \bigl( 
(\partial_+ A_g^{\ast}\times [0,1])\cup \mathcal{H}^2_2\cup \mathcal{H}^3_2 \bigr)\\
&\cong & (\bar{A}_g'\cup \bar{B}_g')\cup (A_g''\cup B_g'').
\end{eqnarray*}

According to this decomposition, we can modify the fork complex $(\mathscr{F},\rho)$ as in 
Figure \ref{fig:35} or Figure \ref{fig:35.1}. 
It is easy to see that for a new complex $(\mathscr{F}^{\ast},\rho^{\ast})$, we have 
$w(\mathscr{F}^{\ast},\rho^{\ast})<w(\mathscr{F},\rho)$, a contradiction. 

\begin{figure}[htb]\begin{center}
  {\unitlength=1cm
  \begin{picture}(10.5,10)
   \put(0.4,6.9){\includegraphics[keepaspectratio]{%
   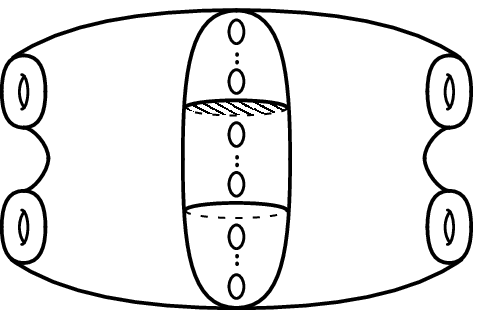}}
   \put(7.4,8){\includegraphics[keepaspectratio]{%
   Fig33p2.eps}}
   \put(0,0){\includegraphics[keepaspectratio]{%
   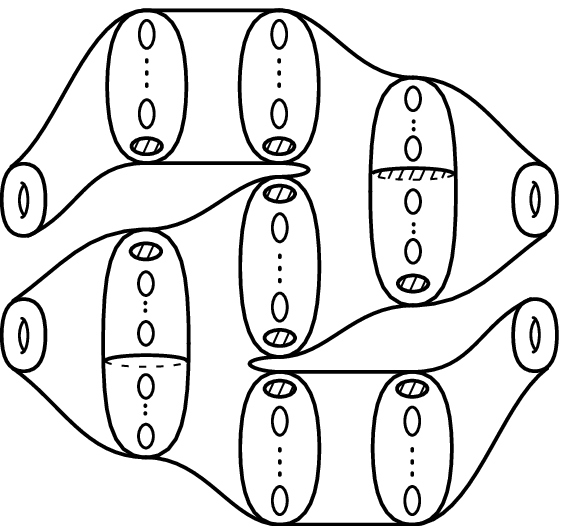}}
   \put(6.4,2.1){\includegraphics[keepaspectratio]{%
   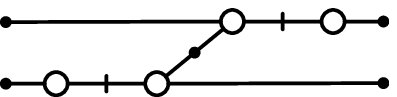}}
   \put(1.6,8.9){$D_A$}
   \put(1.3,7.5){$A_g$}
   \put(3.46,7.9){$D_B$}
   \put(4,7.5){$B_g$}
   \put(2.76,6.7){$g$}
   \put(2.7,5.9){{\Large $\Downarrow$}}
   \put(8.25,5.9){{\Large $\Downarrow$}}
   \put(1.46,5.45){$k$}
   \put(4.2,4.8){$l$}
   \put(1.1,.3){$g-k$}
   \put(3.8,-.3){$g-l$}
   \put(.5,1.7){$\bar{A}'_{g}$}
   \put(2,1.1){$\bar{B}'_{g}$}
   \put(3.3,4.0){$\bar{A}''_{g}$}
   \put(4.7,3.1){$\bar{B}''_{g}$}
  \end{picture}}
  \caption{\sl The case of irreducible splittings}
  \label{fig:35}
\end{center}\end{figure}

\begin{figure}[htb]\begin{center}
  {\unitlength=1cm
  \begin{picture}(10.5,9)
   \put(0.4,5.9){\includegraphics[keepaspectratio]{%
   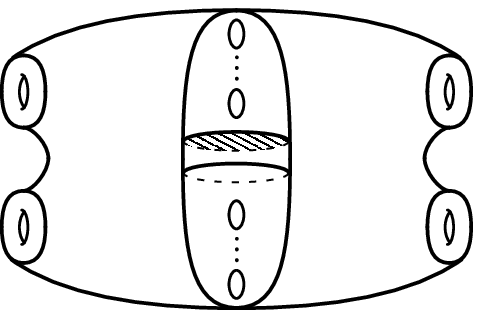}}
   \put(7.4,7){\includegraphics[keepaspectratio]{%
   Fig33p2.eps}}
   \put(0,0){\includegraphics[keepaspectratio]{%
   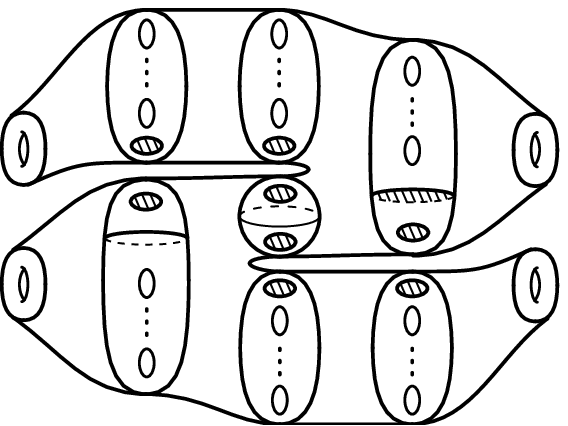}}
   \put(6.4,1.6){\includegraphics[keepaspectratio]{%
   Fig35p4.eps}}
   \put(1.7,7.6){$D_A$}
   \put(1.3,6.5){$A_g$}
   \put(3.4,7.2){$D_B$}
   \put(4,6.5){$B_g$}
   \put(2.76,5.66){$g$}
   \put(2.7,4.9){{\Large $\Downarrow$}}
   \put(8.25,4.9){{\Large $\Downarrow$}}
   \put(.5,1.3){$\bar{A}'_{g}$}
   \put(2,.6){$\bar{B}'_{g}$}
   \put(3.3,3.5){$\bar{A}''_{g}$}
   \put(4.7,2.7){$\bar{B}''_{g}$}
  \end{picture}}	       
  \caption{\sl The case of reducible splittings}
  \label{fig:35.1}
\end{center}\end{figure}

\end{proof}

\begin{lem}\label{S12}
Any component $\rho^{-1}(t)$ is incompressible in $M$ unless $M$ is $\partial$-compressible, 
where $t$ is a tine of 
$\mathscr{F}$. 
\end{lem}

\begin{proof}
Suppose that $\rho^{-1}(t)$ is compressible in $M$ for a tine $t$ of $\mathscr{F}$. Let 
$D$ be a compressing disk of $\rho^{-1}(t)$. Let $\mathcal{T}$ be the union of the tines of 
$\mathscr{F}$. 
By an innermost disk argument, we may assume that $D\cap \rho^{-1}(\mathcal{T})=\partial D$. 
Let $\mathcal{F}_1$ be the fork containing $\rho (\eta(\partial D;D))$. 
Note that $\rho^{-1}(t)$ is incompressible in $\rho^{-1}(\mathcal{F}_1)$ (\textit{cf}. (4) of Remark \ref{R1}). 
Hence there is a fork $\mathcal{F}_2(\ne \mathcal{F}_1)$ attached to the grip, say $g$, of $\mathcal{F}_1$. 
Since $D\cap \rho^{-1}(\mathcal{T})=\partial D$, we have $D\subset \rho^{-1}(\mathcal{F}_1\cup \mathcal{F}_2)$. 
Hence $D$ is a $\partial$-compressing 
disk of $M'=\rho^{-1}(\mathcal{F}_1)\cup \rho^{-1}(\mathcal{F}_2)$. Hence 
it follows from (2) of Theorem \ref{S5} and 
Lemma \ref{S4} that the Heegaard splitting $(\rho^{-1}(\mathcal{F}_1),
\rho^{-1}(\mathcal{F}_2);\rho^{-1}(g))$ is either weakly reducible or trivial. 
It also follows from Proposition \ref{S11} that $(\rho^{-1}(\mathcal{F}_1),
\rho^{-1}(\mathcal{F}_2);\rho^{-1}(g))$ is strongly irreducible and hence the splitting must be trivial. 
Since $\rho^{-1}(\mathcal{F}_1)$ contains $\eta(\partial D;D)$, we see that 
$\rho^{-1}(\mathcal{F}_1)$ is a trivial compression body and that $t$ is the tine of $\mathcal{F}_1$. 
Hence by Lemma \ref{S10}, we have one of the following: $(1)$ $M$ is a $3$-ball, 
$(2)$ $M\cong \rho^{-1}(t)\times [0,1]$ and $(3)$ $M$ is $\partial$-compressible. We suppose that 
$M$ does not satisfy the condition $(3)$, i.e., $M$ is $\partial$-incompressible. 
If $M$ satisfies the condition $(1)$, i.e., $M$ is a $3$-ball, then it follows from Example \ref{thinB3} 
that $t$ is the only tine of $\mathscr{F}$ and that $\rho^{-1}(t)$ is incompressible in $M$. 
This contradicts that we suppose that $\rho^{-1}(t)$ is compressible in $M$. 
If $M$ satisfies the condition $(2)$, i.e., $M\cong \rho^{-1}(t)\times [0,1]$, then $\rho^{-1}(t)$ is 
incompressible in $M$, a contradiction. 
\end{proof}

As a direct consequence of Proposition \ref{S11} and Lemma \ref{S12}, we have the following. 

\begin{cor}
$(\mathscr{F},\rho)$ is strongly irreducible unless $M$ is $\partial$-compressible. 
\end{cor}

\begin{rem}
\rm
There are strongly irreducible splittings which are not thin. 
In fact, there are strongly irreducible Heegaard splittings which are not minimal genus 
(\textit{cf}. \cite{CG2} and \cite{Kobayashi}). 
\end{rem}

\begin{lem}\label{S13}
Suppose that $\mathscr{F}$ contains a tine. There exists a tine $t$ of $\mathscr{F}$ such that 
$\rho^{-1}(t)$ is a $2$-sphere if and only if $M$ is reducible or is a $3$-ball. 
\end{lem}

\begin{proof}
The ``only if part'' is immediate from  Observation \ref{S9}. 
Hence we will give a proof of the ``if part''. 

Suppose that $M$ is reducible or is a $3$-ball. 
If $M$ is a $3$-ball, then it follows from Example \ref{thinB3} that there is exactly one tine, say $t$, 
of $\mathscr{F}$ and $\rho^{-1}(t)=\partial M$ is a $2$-sphere. Hence in the remainder of the proof,  we suppose 
that $M$ is reducible. 
Let $\mathcal{T}$ be the union of the tines of $\mathscr{F}$. 
Let $P$ be a reducing $2$-sphere such that $|P\cap \rho^{-1}(\mathcal{T})|$ is minimal among such 
all reducing $2$-spheres. By an innermost disk argument, we see that $P\cap \rho^{-1}(\mathcal{T})=\emptyset$. 
Let $\mathcal{F}_1$ be a fork of $\mathscr{F}$ with $\rho^{-1}(\mathcal{F}_1)\cap P\ne \emptyset$. 

Suppose first that there are no forks of $\mathscr{F}$ attaching to the grip of $\mathcal{F}_1$. 
Then this implies that $P$ is an essential $2$-sphere in $\rho^{-1}(t)$. 

Suppose next that there is a fork of $\mathscr{F}$, say $\mathcal{F}_2$, other than $\mathcal{F}_1$ 
which attaches to the grip, say $g$, of $\mathcal{F}_1$. Note that 
$\rho^{-1}(\mathcal{F}_1)\cup 
\rho^{-1}(\mathcal{F}_2)$ contains $P$. It follows from (1) of Theorem \ref{S5} that 
$\rho^{-1}(\mathcal{F}_1)$ or $\rho^{-1}(\mathcal{F}_2)$ is reducible, or $(\rho^{-1}(\mathcal{F}_1),
\rho^{-1}(\mathcal{F}_2);\rho^{-1}(g))$ is reducible. 
The latter condition, however, contradicts Proposition \ref{S11}. 
Hence we may assume that $\rho^{-1}(\mathcal{F}_1)$ is reducible, that is, there is a $2$-sphere component $P_0$ 
of $\partial_- (\rho^{-1}(\mathcal{F}_1))$ (\textit{cf}. (1) of Remark \ref{R1}). 
This implies that there is a tine $t$ with $\rho^{-1}(t)=P_0$. 
\end{proof}

\begin{lem}\label{S14}
If some $\rho^{-1}(g)$ is a torus, where $g$ is a grip of $\mathscr{F}$, then one of the 
following holds. \\
$(1)$\ \ $M$ is reducible.\\
$(2)$\ \ $M$ is $($a torus$)\times [0,1]$.\\
$(3)$\ \ $M$ is a solid torus.\\
$(4)$\ \ $M$ is a lens space.
\end{lem}

\begin{proof}
Suppose that $M$ does not satisfy the conclusion (1) of Lemma \ref{S14}, i.e., $M$ is irreducible. 
Note that $\rho^{-1}(g)$ may be a boundary component 
of $M$. Let $\mathcal{F}$ be a fork such that the grip of $\mathcal{F}$ is $g$. Set 
$V=\rho^{-1}(\mathcal{F})$. 

If $V$ is trivial, then $M$ is either $T^2\times [0,1]$ or a 
solid torus by Lemma \ref{S10}. Hence conclusion $(2)$ or $(3)$ of Lemma \ref{S14} holds. 

If $V$ is non-trivial, then we see that $V$ is a solid torus by Observation \ref{S9} and 
Example \ref{thinB3}. Suppose further that the conclusion $(3)$ does not hold., i.e., 
$M$ is not a solid torus. Then there is a fork $\mathcal{F}'(\ne \mathcal{F})$ attached to $g$. 
Set $V'=\rho^{-1}(\mathcal{F}')$. If $V'$ is trivial, then it follows from  Lemma \ref{S10} that 
$M$ is a solid torus, a contradiction. 
If $V'$ is non-trivial, then we see that $V'$ is a solid torus by Observation \ref{S9} and 
Example \ref{thinB3}. Hence $M$ is a lens space and we have the conclusion $(4)$ of Lemma \ref{S14}. 
\end{proof}

\subsection{Examples of generalized Heegaard splittings}

In this section, we use some theorems without proofs to obtain generalized Heegaard 
splittings and associated fork complexes. 
Let $F_g$ be a connected closed orientable surface of genus $g$. 

\bigskip\noindent
$\bullet \ \ \ \  M=F_g\times [0,1]$. 

\medskip
Set $M=F_g\times [0,1]$, 
$A=F_g\times [0,1/2]$, $B=F_g\times [1/2,1]$ and $S=F_g\times \{1/2\}$. 
Clearly, $(A,B;S)$ is a Heegaard splitting of $M$, and we call this Heegaard splitting 
the \textit{trivial Heegaard splitting of type I}. Let $p$ be a point in $F_g$. 
Set 

\bigskip
\begin{center}
$A'=\eta((F_g\times \{0\})\cup (p\times [0,1])\cup (F_g\times \{1\});M)$, 
\end{center}

\medskip\noindent
$B'=\mathrm{cl}(M\setminus A')$ and $S'=A'\cap B'$. Then $(A',B';S')$ is also a Heegaard splitting of 
$M$, and we call this splitting the \textit{trivial Heegaard splitting of type II}. 
The proof of the next observation is left to the reader. 

\begin{obs}
Both these Heegaard splittings are strongly irreducible. 
\end{obs}

In fact, Scharlemann and Thompson proved the following. 

\begin{thm}[\cite{Scharlemann2} 2.11 Main Theorem]\label{S15}
Any irreducible Heegaard splitting of $F_g\times [0,1]$ is trivial of type I or II. 
\end{thm}

\begin{figure}[htb]\begin{center}
  {\unitlength=1cm
  \begin{picture}(4,1.8)
   \put(1.8,0){\includegraphics[keepaspectratio]{%
   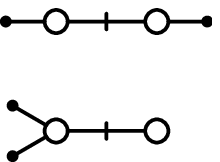}}
   \put(0,1.4){Type I}
   \put(0,.26){Type II}
  \end{picture}}
  \caption{}
  \label{fig:36}
\end{center}\end{figure}

We remark that the fork complexes associated to these Heegaard splittings are 
illustrated in Figure \ref{fig:36}. 

\bigskip\noindent
$\bullet \ \ \ \ M=F_g\times S^1$. 

\medskip
Note that $S^1$ is regarded as $[0,1]/\{0\}\sim \{1\}$. 
Let $p$ and $q$ be distinct points in $F_g$. Set 

\medskip
\begin{center}
$A=\mathrm{cl}((F_g\times [0,1/2])\setminus \eta(p\times [0,1/2];F_g\times [0,1/2]))
\cup \eta(q\times [1/2,1];F_g\times [1/2,1])$ 
\end{center}

\medskip\noindent
and 

\begin{eqnarray*}
B 
&=& \mathrm{cl}(M\setminus A)\\
&=& \mathrm{cl}((F_g\times [1/2,1])\setminus \eta(q\times [1/2,1];F_g\times [1/2,1]))\cup 
\eta(p\times [0,1/2];F_g\times [0,1/2]). 
\end{eqnarray*}

\begin{figure}[htb]\begin{center}
  {\unitlength=1cm
  \begin{picture}(7.6,5.8)
   \put(0,0){\includegraphics[keepaspectratio]{%
   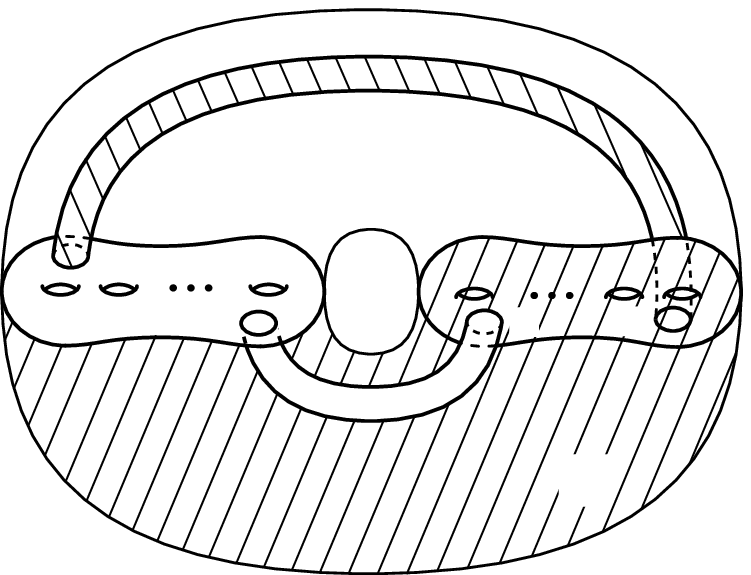}}
   \put(.26,3.1){$p$}
   \put(2.2,2.6){$q$}
   \put(5.24,2.6){$q$}
   \put(6.35,2.6){$p$}
   \put(5.8,.9){$A$}
   \put(5.8,3.9){$B$}
  \end{picture}}
  \caption{}
  \label{fig:37}
\end{center}\end{figure}

Note that $A$ and $B$ are handlebodies. Set $S=\partial A\cap \partial B$. 
Then $(A,B;S)$ is a Heegaard splitting of $M=F_g\times S^1$ and is 
called the \textit{trivial Heegaard 
splitting} of $M=F_g\times S^1$ (\textit{cf}. Figure \ref{fig:37}). 

\begin{exercise}
Show that this trivial Heegaard splitting is weakly reducible. 
\end{exercise}

\begin{thm}[\cite{Schultens} Theorem 5.7]\label{S16}
Any irreducible Heegaard splitting of $F_g\times S^1$ is the trivial splitting. 
\end{thm}

\bigskip\noindent
$\bullet \ \ \ \  T^3=T^2\times S^1$. 

\begin{figure}[htb]\begin{center}
  {\unitlength=1cm
  \begin{picture}(10,5.6)
   \put(0,.8){\includegraphics[keepaspectratio]{%
   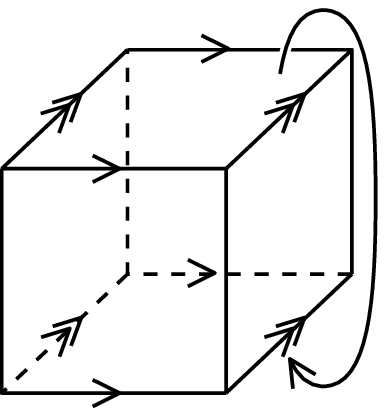}}
   \put(6,.8){\includegraphics[keepaspectratio]{%
   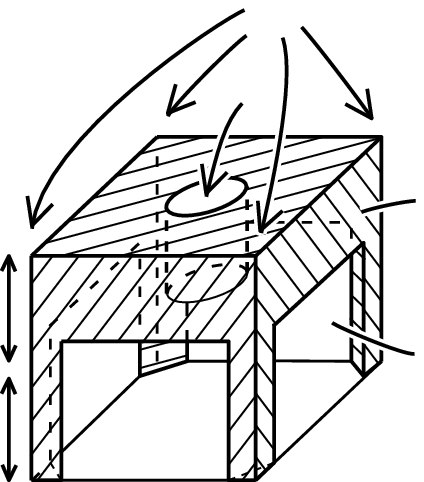}}
   \put(5.6,2.5){$I_1$}
   \put(5.6,1.3){$I_2$}
   \put(8.7,5.6){$p$}
   \put(8.5,4.84){$q$}
   \put(10.3,3.6){$A$}
   \put(10.3,2.0){$B$}

   \put(1.6,0){(a)}
   \put(7.8,0){(b)}

  \end{picture}}
  \caption{}
  \label{fig:38}
\end{center}\end{figure}

\medskip
It is known that $T^3$ is obtained from a cube $[0,1]\times [0,1]\times [0,1]$ 
by attaching corresponding edges and faces as in Figure \ref{fig:38} (a). Set 

\bigskip
\begin{center}
$A=\mathrm{cl}((T^2\times [0,1/2])\setminus \eta(p\times [0,1/2];T^2\times [0,1/2]))
\cup \eta(q\times [1/2,1];T^2\times [1/2,1])$ 
\end{center}

\medskip\noindent
and 

\begin{eqnarray*}
B 
&=& \mathrm{cl}(T^3\setminus A)\\
&=& \mathrm{cl}((T^2\times [1/2,1])\setminus \eta(q\times [1/2,1];T^2\times [1/2,1]))\cup 
\eta(p\times [0,1/2];T^2\times [0,1/2]). 
\end{eqnarray*}

Then we see that $A$ and $B$ are genus two handlebodies and that 
it follows from Theorem \ref{S16} that $(A,B;S)$ is the Heegaard splitting of $T^3$, 
where $S=\partial A=\partial B$ (\textit{cf}. Figure \ref{fig:38} (b)). 
Set $h^1=\eta(q\times [1/2,1];T^2\times [1/2,1])$ and 
$h^2=\eta(p\times [0,1/2];T^2\times [0,1/2])$. 
Note that $h^1$ ($h^2$ resp.) can be regarded as a $1$-handle ($2$-handle resp.) in a 
handle decomposition of $T^3$ obtained from the Heegaard splitting $(A,B;S)$. 
Since $h^1\cap h^2=\emptyset$, we can perform a weak reduction to obtain a generalized 
Heegaard splitting. We give a concrete description of the generalized Heegaard splitting in 
the following. 
First, set $A_1=\mathrm{cl}(T^3\times [0,1/2]\setminus h^1)$ and 
$B_2=\mathrm{cl}(T^3\times [1/2,1]\setminus h^2)$. 
That is, $A_1$ is obtained from $A$ by removing the $1$-handle $h^1$ and 
$B_2$ is obtained from $B$ by removing the $2$-handle $h^2$. Then we have:  

\begin{eqnarray*}
T^3 
&=& A\cup B\\
&=& A_1\cup h^1\cup h^2\cup B_2\\
&\cong & A_1\cup (\partial A_1\times [0,1])\cup h^1\cup h^2\cup B_2\\
&=& A_1\cup \bigl( (\partial A_1\times [0,1])\cup h^2\bigr) \cup h^1\cup B_2.
\end{eqnarray*}

Set $B_1=(\partial A_1\times [0,1])\cup h^2$. Then $B_1$ is a compression body 
such that $\partial_+ B_1=\partial A_1$ and $\partial_- B$ consists of two tori. 
Hence we have: 

\begin{eqnarray*}
T^3 
&\cong & A_1\cup B_1\cup h^1\cup B_2\\
&\cong & A_1\cup B_1\cup \bigl( (\partial_- B_1\times [0,1])\cup h^1\bigr) \cup B_2.
\end{eqnarray*}

Set $A_2=(\partial_- B_1\times [0,1])\cup h^1$. Then $A_2$ is a compression body 
such that $\partial_+ A_2=\partial B_2$ and $\partial_- A_2=\partial_- B_1$. 
Hence we have: 

\bigskip
\begin{center}
$T^3=(A_1\cup B_1)\cup (A_2\cup B_2)$. 
\end{center}

\begin{figure}[htb]\begin{center}
  {\unitlength=1cm
  \begin{picture}(9.3,6.7)
   \put(0,3.4){\includegraphics[keepaspectratio]{%
   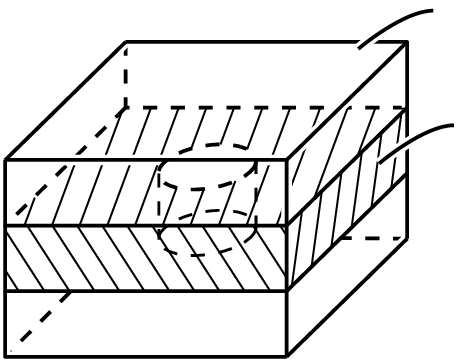}}
   \put(6,3.4){\includegraphics[keepaspectratio]{%
   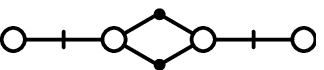}}
   \put(0,0){\includegraphics[keepaspectratio]{%
   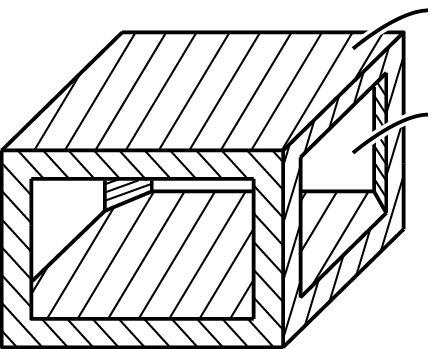}}
   \put(4.5,6.9){$B_1$}
   \put(4.7,5.7){$A_1$}
   \put(4.5,3.4){$A_2$}
   \put(4.5,2.27){$B_2$}
   \put(5.9,3.1){$A_1$}
   \put(7,3.1){$B_1$}
   \put(7.9,3.1){$A_2$}
   \put(8.9,3.1){$B_2$}
  \end{picture}}
  \caption{}
  \label{fig:39}
\end{center}\end{figure}

\medskip\noindent
This together with the fork complex as in Figure \ref{fig:39} gives a generalized Heegaard splitting. 

\begin{exercise}
\rm
Show that this is the only fork complex associated with a generalized Heegaard splitting 
of $T^3$ via weak reduction. 
\end{exercise}

\begin{rem}
\rm
The inverse procedure of weak reduction is called \textit{amalgamation}. 
\end{rem}

\begin{exercise}
\rm
Show that the above generalized Heegaard splitting of $T^3$ is strongly irreducible. 
\end{exercise}

\bigskip\noindent
$\bullet \ \ \ \  M=F_2\times S^1$. 

\begin{figure}[htb]\begin{center}
  {\unitlength=1cm
  \begin{picture}(11.6,5.6)
   \put(0,0.8){\includegraphics[keepaspectratio]{%
   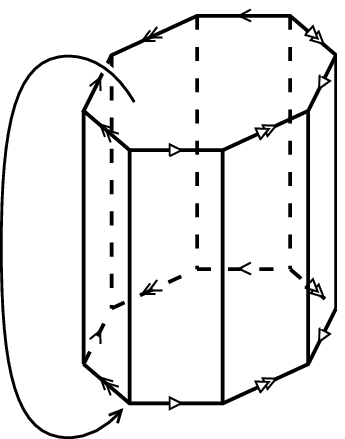}}
   \put(4.6,0.8){\includegraphics[keepaspectratio]{%
   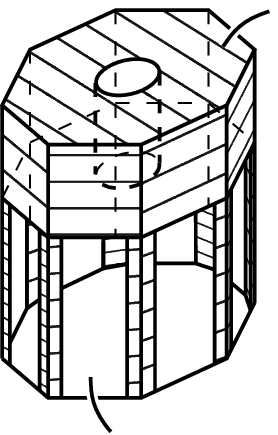}}
   \put(8.6,0.8){\includegraphics[keepaspectratio]{%
   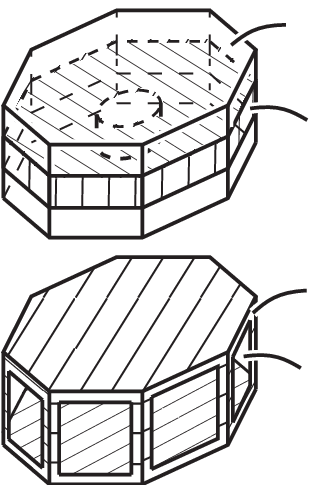}}
   \put(.65,4.3){$a$}
   \put(1.3,5){$b$}
   \put(2.5,5.3){$a$}
   \put(3.3,5){$d$}
   \put(3.46,4.3){$c$}
   \put(2.6,3.5){$d$}
   \put(1.7,3.4){$c$}
   \put(.9,3.5){$b$}
   \put(.91,2.04){$a$}
   \put(1.5,2.5){$b$}
   \put(2.5,2.7){$a$}
   \put(3.2,2.5){$d$}
   \put(3.46,1.8){$c$}
   \put(2.6,1){$d$}
   \put(1.7,.9){$c$}
   \put(.9,1.1){$b$}
   \put(5.8,.7){$B$}
   \put(7.44,5){$A$}
   \put(11.6,5.4){$B_1$}
   \put(11.7,4.3){$A_1$}
   \put(11.74,2.7){$A_2$}
   \put(11.74,1.83){$B_2$}
   \put(1.85,0){(a)}
   \put(5.7,0){(b)}
   \put(9.7,0){(c)}
  \end{picture}}
  \caption{}
  \label{fig:40}
\end{center}\end{figure}

\medskip
$M$ is obtained from a (an octagon)
$\times [0,1]$ by attaching corresponding edges and faces as in Figure \ref{fig:40} (a). 
Set 

\bigskip
\begin{center}
$A=\mathrm{cl}((F_g\times [0,1/2])\setminus \eta(p\times [0,1/2];F_g\times [0,1/2]))
\cup \eta(q\times [1/2,1];F_g\times [1/2,1])$ 
\end{center}

\medskip\noindent
and 

\begin{eqnarray*}
B 
&=& \mathrm{cl}(M\setminus A)\\
&=& \mathrm{cl}((F_g\times [1/2,1])\setminus \eta(q\times [1/2,1];F_g\times [1/2,1]))\cup 
\eta(p\times [0,1/2];F_g\times [0,1/2]). 
\end{eqnarray*}

\noindent
Then it follows from Theorem \ref{S16} that we obtain the Heegaard splitting $M=A\cup B$ (
see Figure \ref{fig:40} (b)). 
As descibed in case of $M=T^3$, we can perform a weak reduction and 
we obtain the same fork complex as that illustrated in Figure \ref{fig:39}. In this case, each of $A_1$ 
and $B_2$ is a handlebody of genus four and each of $A_2$ and $B_1$ is a compression body 
with $\partial_+ A_2=\partial B_2$, $\partial A_1=\partial_+ B_1$ and $\partial_- A_2=
\partial_- B_1$. 

For the Heegaard splitting $A\cup B$ of $M$, we can find another weak reduction as 
follows. 
Recall that $M=P_8\times [0,1]/\sim$, where $P_8$ is an octagon (\textit{cf}. Figure \ref{fig:40}). 
Then there is a handle decomposition 

\bigskip
\begin{center}
$M=h^0\cup h^1_a\cup h^1_b\cup h^1_c\cup h^1_d\cup h^1_e\cup 
h^2_a\cup h^2_b\cup h^2_c\cup h^2_d\cup h^2_e\cup h^3$, 
\end{center}

\medskip\noindent
where a $0$-handle $h^0$ corrsponds to a vertex of $P_8$, 
a $1$-handle $h^1_a$ ($h^1_b$, $h^1_c$, $h^1_d$ and $h^1_e$ resp.) corrsponds to 
$a$ ($b$, $c$, $d$ and $e$ resp.) in $P_8$, 
a $2$-handle $h^2_a$ ($h^2_b$, $h^2_c$ and $h^2_d$ resp.) corrsponds to 
the face bounded by $eae^{-1}a^{-1}$ ($ebe^{-1}b^{-1}$, $ece^{-1}c^{-1}$ and $ede^{-1}d^{-1}$ resp.) 
in $\partial P_8\times [0,1]$, 
a $2$-handle $h^2_e$ corrsponds to the face bounded by $aba^{-1}b^{-1}cdc^{-1}d^{-1}$ in $P_8$ and 
a $3$-handle $h^3$ corrsponds to the vertex in the interior of $P_8\times [0,1]$. 
Set $A_1=h^0\cup h^1_a\cup h^1_e$. Then $A_1$ is a genus two handlebody and we have: 

\begin{eqnarray*}
M 
&=& A_1\cup h^1_b\cup h^1_c\cup h^1_d\cup 
h^2_a\cup h^2_b\cup h^2_c\cup h^2_d\cup h^2_e\cup h^3\\
&\cong & A_1\cup (\partial A_1\times [0,1]) \cup h^1_b\cup h^1_c\cup h^1_d\cup 
h^2_a\cup h^2_b\cup h^2_c\cup h^2_d\cup h^2_e\cup h^3\\
&=& A_1\cup \bigl( (\partial A_1\times [0,1]) \cup h^2_a\bigr) \cup h^1_b\cup h^1_c\cup h^1_d\cup 
\cup h^2_b\cup h^2_c\cup h^2_d\cup h^2_e\cup h^3 
\end{eqnarray*}

\noindent
Set $B_1=(\partial A_1\times [0,1]) \cup h^2_a$. Then $B_1$ is a compression body such that 
$\partial_+ B_1=\partial A_1$ and $\partial_- B_1$ consists of two tori. Then we have: 

\begin{eqnarray*}
M 
&\cong & A_1\cup B_1 \cup h^1_b\cup h^1_c\cup h^1_d\cup 
\cup h^2_b\cup h^2_c\cup h^2_d\cup h^2_e\cup h^3\\
&\cong & A_1\cup B_1 \cup (\partial_- B_1\times [0,1]) \cup h^1_b\cup h^1_c\cup h^1_d\cup 
\cup h^2_b\cup h^2_c\cup h^2_d\cup h^2_e\cup h^3\\
&=& A_1\cup B_1 \cup \bigl( (\partial_- B_1\times [0,1])\cup h^1_b \bigr) \cup h^1_c\cup h^1_d\cup 
\cup h^2_b\cup h^2_c\cup h^2_d\cup h^2_e\cup h^3. 
\end{eqnarray*}

\noindent
Set $A_2=(\partial_- B_1\times [0,1]) \cup h^1_b$. Then $A_2$ is a compression body such that 
$\partial_+ A_2$ is a closed surface of genus two and $\partial_- A_2=\partial_- B_1$. Then we have: 

\begin{eqnarray*}
M 
&\cong& A_1\cup B_1 \cup A_2 \cup h^1_c\cup h^1_d\cup 
\cup h^2_b\cup h^2_c\cup h^2_d\cup h^2_e\cup h^3\\
&\cong& A_1\cup B_1 \cup A_2 \cup (\partial_+ A_2\times [0,1])\cup h^1_c\cup h^1_d\cup 
\cup h^2_b\cup h^2_c\cup h^2_d\cup h^2_e\cup h^3\\
&=& A_1\cup B_1 \cup A_2 \cup \bigl( (\partial_+ A_2\times [0,1])\cup h^2_b \bigr) 
\cup h^1_c\cup h^1_d\cup h^2_c\cup h^2_d\cup h^2_e\cup h^3. 
\end{eqnarray*}

\noindent
Set $B_2=(\partial_+ A_2\times [0,1])\cup h^2_b$. Then $B_2$ is a compression body such that 
$\partial_+ B_2=\partial_+ A_2$ and $\partial_- B_2$ consists of a torus. Then we have: 

\begin{eqnarray*}
M 
&\cong& A_1\cup B_1 \cup A_2 \cup B_2 \cup h^1_c\cup h^1_d\cup h^2_c\cup h^2_d\cup h^2_e\cup h^3\\
&\cong& A_1\cup B_1 \cup A_2 \cup B_2 \cup (\partial_- B_2\times [0,1])
\cup h^1_c\cup h^1_d\cup h^2_c\cup h^2_d\cup h^2_e\cup h^3\\
&=& A_1\cup B_1 \cup A_2 \cup B_2 \bigl( (\partial_- B_2\times [0,1])\cup h^1_c\bigr) 
\cup h^1_d\cup h^2_c\cup h^2_d\cup h^2_e\cup h^3. 
\end{eqnarray*}

\noindent
Set $A_3=(\partial_- B_2\times [0,1]) \cup h^1_c$. Then $A_3$ is a compression body such that 
$\partial_+ A_3$ is a closed surface of genus two and $\partial_- A_3=\partial_- B_2$. Then we have: 

\begin{eqnarray*}
M 
&\cong& A_1\cup B_1 \cup A_2 \cup B_2 \cup A_3\cup h^1_d\cup h^2_c\cup h^2_d\cup h^2_e\cup h^3\\
&\cong& A_1\cup B_1 \cup A_2 \cup B_2 \cup A_3\cup (\partial_+ A_3\times [0,1])
\cup h^1_d\cup h^2_c\cup h^2_d\cup h^2_e\cup h^3\\
&=& A_1\cup B_1 \cup A_2 \cup B_2 \cup A_3\cup \bigl( (\partial_+ A_3\times [0,1])\cup h^2_c\bigr) 
\cup h^1_d\cup h^2_d\cup h^2_e\cup h^3. 
\end{eqnarray*}

\noindent
Set $B_3=(\partial_+ A_3\times [0,1])\cup h^2_c$. Then $B_3$ is a compression body such that 
$\partial_+ B_3=\partial_+ A_3$ and $\partial_- B_3$ consists of two tori. Then we have: 

\begin{eqnarray*}
M 
&\cong& A_1\cup B_1 \cup A_2 \cup B_2 \cup A_3\cup B_3\cup h^1_d\cup h^2_d\cup h^2_e\cup h^3\\
&\cong& A_1\cup B_1 \cup A_2 \cup B_2 \cup A_3\cup B_3\cup (\partial_- B_3\times [0,1])
\cup h^1_d\cup h^2_d\cup h^2_e\cup h^3\\
&=& A_1\cup B_1 \cup A_2 \cup B_2 \cup A_3\cup B_3\cup 
\bigl( (\partial_- B_3\times [0,1])\cup h^1_d \bigr)
\cup h^2_d\cup h^2_e\cup h^3. 
\end{eqnarray*}

\noindent
Set $A_4=(\partial_- B_4\times [0,1]) \cup h^1_d$. Then $A_4$ is a compression body such that 
$\partial_+ A_4$ is a closed surface of genus two and $\partial_- A_4=\partial_- B_3$. Then we have: 

\begin{eqnarray*}
M 
&\cong& A_1\cup B_1 \cup A_2 \cup B_2 \cup A_3\cup B_3\cup A_4\cup h^2_d\cup h^2_e\cup h^3\\
&\cong& A_1\cup B_1 \cup A_2 \cup B_2 \cup A_3\cup B_3\cup A_4\cup (\partial_+ A_4\times [0,1])
\cup h^2_d\cup h^2_e\cup h^3.
\end{eqnarray*}

\noindent
Set $B_4=(\partial_+ A_4\times [0,1])\cup h^2_d\cup h^2_e\cup h^3$. 
Then $B_4$ is a genus two handlebody such that 
$\partial B_4=\partial_+ A_4$. Therefore we have the following decomposition. 

\bigskip
\begin{center}
$T^3=(A_1\cup B_1)\cup (A_2\cup B_2)\cup (A_3\cup B_3)\cup (A_4\cup B_4)$. 
\end{center}

\begin{figure}[htb]\begin{center}
  {\unitlength=1cm
  \begin{picture}(7.4,2.5)
   \put(0,.8){\includegraphics[keepaspectratio]{%
   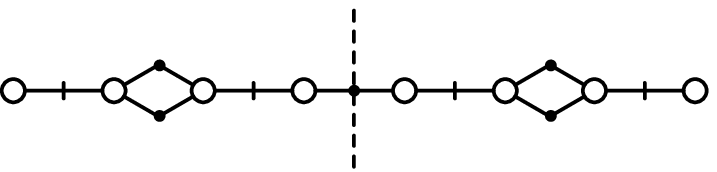}}
   \put(0,2.1){$A_1$}
   \put(1,2.1){$B_1$}
   \put(1.9,2.1){$A_2$}
   \put(2.9,2.1){$B_2$}
   \put(3.9,2.1){$A_3$}
   \put(4.94,2.1){$B_3$}
   \put(5.8,2.1){$A_4$}
   \put(6.8,2.1){$B_4$}
   \put(0,1.1){$\underbrace{\hspace*{3.2cm}}$}
   \put(0,.4){(a torus with}
   \put(.5,0){a single hole)$\times S^1$}
   \put(4,1.1){$\underbrace{\hspace*{3.2cm}}$}
   \put(4,.4){(a torus with}
   \put(4.5,0){a single hole)$\times S^1$}
  \end{picture}}
  \caption{}
  \label{fig:41}
\end{center}\end{figure}

\medskip\noindent
This together with the fork complex as in Figure \ref{fig:41} gives a generalized Heegaard splitting. 
We remark that $(A_1\cup B_1)\cup (A_2\cup B_2)$ ($(A_3\cup B_3)\cup (A_4\cup B_4)$ resp.) composes 
a $($a torus with a single hole$)\times S^1$. 

\begin{figure}[htb]\begin{center}
  {\unitlength=1cm
  \begin{picture}(7.3,4)
   \put(0,2.2){\includegraphics[keepaspectratio]{%
   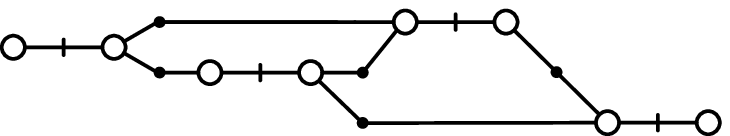}}
   \put(0,.2){\includegraphics[keepaspectratio]{%
   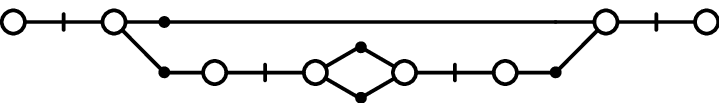}}   
   \put(0,2.6){$A_1$}
   \put(1,2.6){$B_1$}
   \put(2,2.4){$A_3$}
   \put(2.9,2.4){$B_3$}
   \put(3.9,3.7){$A_2$}
   \put(4.94,3.7){$B_2$}
   \put(6,1.9){$A_4$}
   \put(7,1.9){$B_4$}
   \put(0,.6){$A_1$}
   \put(.9,.6){$B_1$}
   \put(1.9,.1){$A_3$}
   \put(2.9,.1){$B_3$}
   \put(3.9,.1){$A_2$}
   \put(4.94,.1){$B_2$}
   \put(6,.6){$A_4$}
   \put(7,.6){$B_4$}
  \end{picture}}
  \caption{}
  \label{fig:42}
\end{center}\end{figure}

By changing the attaching order of $h^1_b$, $h^1_c$ and $h^1_d$, we can obtain 
two more strongly irreducible generalized Heegaard splittings via 
weak reduction (\textit{cf}. Figure \ref{fig:42}). 

\begin{exercise}
\rm
Show that these are the only fork complexes associated with a generalized Heegaard 
splitting of $F_2\times S^1$ via weak reduction. 
\end{exercise}

\begin{rem}
\rm
The fork complexes associated with distinct weak reduction of a Heegaard splitting 
need not homotopic. 
\end{rem}

\bigskip
{\sc Department of Mathematics, Graduate School of Science, Osaka University, 
Machikaneyama 1-16, Toyonaka, Osaka 560-0043, Japan}

\textit{E-mail address}: \texttt{saito@gaia.math.wani.osaka-u.ac.jp}

\bigskip
{\sc Mathematics Department, University of California, Santa Barbara, CA93106, USA}

\textit{E-mail address}: \texttt{mgscharl@math.ucsb.edu}

\bigskip
{\sc Department of Mathematics, One Shields Avenue, University of California, Davis, CA 95616, USA}

\textit{E-mail address}: \texttt{jcs@math.ucdavis.edu}

\end{document}